
\ifx\shlhetal\undefinedcontrolsequence\let\shlhetal\relax\fi
\def\fmtname{AmS-TeX}

\def\fmtversion{2.1}
\catcode`\@=11
\ifx\amstexloaded@\relax\catcode`\@=\active
  \endinput\else\let\amstexloaded@\relax\fi
\newlinechar=`\^^J
\def\W@{\immediate\write\sixt@@n}
\def\CR@{\W@{^^J\fmtname - Version \fmtversion^^J%
COPYRIGHT 1985, 1990, 1991 - AMERICAN MATHEMATICAL SOCIETY^^J%
Use of this macro package is not restricted provided^^J%
each use is acknowledged upon publication.^^J}}
\CR@ \everyjob{\CR@}
\message{Loading definitions for}
\message{misc utility macros,}
\toksdef\toks@@=2
\long\def\rightappend@#1\to#2{\toks@{\\{#1}}\toks@@
 =\expandafter{#2}\xdef#2{\the\toks@@\the\toks@}\toks@{}\toks@@{}}
\def\alloclist@{}
\newif\ifalloc@
\def\showallocations{{\def\\{\immediate\write\m@ne}\alloclist@}\alloc@true}
\def\alloc@#1#2#3#4#5{\global\advance\count1#1by\@ne
 \ch@ck#1#4#2\allocationnumber=\count1#1
 \global#3#5=\allocationnumber
 \edef\next@{\string#5=\string#2\the\allocationnumber}%
 \expandafter\rightappend@\next@\to\alloclist@}
\newcount\count@@
\newcount\count@@@
\def\FN@{\futurelet\next}
\def\DN@{\def\next@}
\def\DNii@{\def\nextii@}
\def\RIfM@{\relax\ifmmode}
\def\RIfMIfI@{\relax\ifmmode\ifinner}
\def\setboxz@h{\setbox\z@\hbox}
\def\wdz@{\wd\z@}
\def\boxz@{\box\z@}
\def\setbox@ne{\setbox\@ne}
\def\wd@ne{\wd\@ne}
\def\iterate{\body\expandafter\iterate\else\fi}
\def\err@#1{\errmessage{AmS-TeX error: #1}}
\newhelp\defaulthelp@{Sorry, I already gave what help I could...^^J
Maybe you should try asking a human?^^J
An error might have occurred before I noticed any problems.^^J
``If all else fails, read the instructions.''}
\def\Err@{\errhelp\defaulthelp@\err@}
\def\eat@#1{}
\def\in@#1#2{\def\in@@##1#1##2##3\in@@{\ifx\in@##2\in@false\else\in@true\fi}%
 \in@@#2#1\in@\in@@}
\newif\ifin@
\def\space@.{\futurelet\space@\relax}
\space@. %
\newhelp\athelp@
{Only certain combinations beginning with @ make sense to me.^^J
Perhaps you wanted \string\@\space for a printed @?^^J
I've ignored the character or group after @.}
{\catcode`\~=\active 
 \lccode`\~=`\@ \lowercase{\gdef~{\FN@\at@}}}
\def\at@{\let\next@\at@@
 \ifcat\noexpand\next a\else\ifcat\noexpand\next0\else
 \ifcat\noexpand\next\relax\else
   \let\next\at@@@\fi\fi\fi
 \next@}
\def\at@@#1{\expandafter
 \ifx\csname\space @\string#1\endcsname\relax
  \expandafter\at@@@ \else
  \csname\space @\string#1\expandafter\endcsname\fi}
\def\at@@@#1{\errhelp\athelp@ \err@{\Invalid@@ @}}
\def\atdef@#1{\expandafter\def\csname\space @\string#1\endcsname}
\newhelp\defahelp@{If you typed \string\define\space cs instead of
\string\define\string\cs\space^^J
I've substituted an inaccessible control sequence so that your^^J
definition will be completed without mixing me up too badly.^^J
If you typed \string\define{\string\cs} the inaccessible control sequence^^J
was defined to be \string\cs, and the rest of your^^J
definition appears as input.}
\newhelp\defbhelp@{I've ignored your definition, because it might^^J
conflict with other uses that are important to me.}
\def\define{\FN@\define@}
\def\define@{\ifcat\noexpand\next\relax
 \expandafter\define@@\else\errhelp\defahelp@                               
 \err@{\string\define\space must be followed by a control
 sequence}\expandafter\def\expandafter\nextii@\fi}                          
\def\undefined@@@@@@@@@@{}
\def\preloaded@@@@@@@@@@{}
\def\next@@@@@@@@@@{}
\def\define@@#1{\ifx#1\relax\errhelp\defbhelp@                              
 \err@{\string#1\space is already defined}\DN@{\DNii@}\else
 \expandafter\ifx\csname\expandafter\eat@\string                            
 #1@@@@@@@@@@\endcsname\undefined@@@@@@@@@@\errhelp\defbhelp@
 \err@{\string#1\space can't be defined}\DN@{\DNii@}\else
 \expandafter\ifx\csname\expandafter\eat@\string#1\endcsname\relax          
 \global\let#1\undefined\DN@{\def#1}\else\errhelp\defbhelp@
 \err@{\string#1\space is already defined}\DN@{\DNii@}\fi
 \fi\fi\next@}

\def\predefine#1#2{\let#1#2}
\def\undefine#1{\let#1\undefined}
\message{page layout,}
\newdimen\captionwidth@
\captionwidth@\hsize
\advance\captionwidth@-1.5in
\def\pagewidth#1{\hsize#1\relax
 \captionwidth@\hsize\advance\captionwidth@-1.5in}
\def\pageheight#1{\vsize#1\relax}
\def\hcorrection#1{\advance\hoffset#1\relax}
\def\vcorrection#1{\advance\voffset#1\relax}
\message{accents/punctuation,}

\let\graveaccent\`
\let\acuteaccent\'
\let\tildeaccent\~
\let\hataccent\^
\let\underscore\_
\let\B\=
\let\D\.
\let\ic@\/
\def\/{\unskip\ic@}
\def\textfonti{\the\textfont\@ne}
\def\t#1#2{{\edef\next@{\the\font}\textfonti\accent"7F \next@#1#2}}
\def~{\unskip\nobreak\ \ignorespaces}
\def\.{.\spacefactor\@m}
\atdef@;{\leavevmode\null;}
\atdef@:{\leavevmode\null:}
\atdef@?{\leavevmode\null?}
\edef\@{\string @}
\def\relaxnext@{\let\next\relax}
\atdef@-{\relaxnext@\leavevmode
 \DN@{\ifx\next-\DN@-{\FN@\nextii@}\else
  \DN@{\leavevmode\hbox{-}}\fi\next@}%
 \DNii@{\ifx\next-\DN@-{\leavevmode\hbox{---}}\else
  \DN@{\leavevmode\hbox{--}}\fi\next@}%
 \FN@\next@}
\def\srdr@{\kern.16667em}
\def\drsr@{\kern.02778em}
\def\sldl@{\drsr@}
\def\dlsl@{\srdr@}
\atdef@"{\unskip\relaxnext@
 \DN@{\ifx\next\space@\DN@. {\FN@\nextii@}\else
  \DN@.{\FN@\nextii@}\fi\next@.}%
 \DNii@{\ifx\next`\DN@`{\FN@\nextiii@}\else
  \ifx\next\lq\DN@\lq{\FN@\nextiii@}\else
  \DN@####1{\FN@\nextiv@}\fi\fi\next@}%
 \def\nextiii@{\ifx\next`\DN@`{\sldl@``}\else\ifx\next\lq
  \DN@\lq{\sldl@``}\else\DN@{\dlsl@`}\fi\fi\next@}%
 \def\nextiv@{\ifx\next'\DN@'{\srdr@''}\else
  \ifx\next\rq\DN@\rq{\srdr@''}\else\DN@{\drsr@'}\fi\fi\next@}%
 \FN@\next@}

\def\textfontii{\the\textfont\tw@}
\def\lbrace@{\delimiter"4266308 }
\def\rbrace@{\delimiter"5267309 }
\def\{{\RIfM@\lbrace@\else{\textfontii f}\spacefactor\@m\fi}
\def\}{\RIfM@\rbrace@\else
 \let\@sf\empty\ifhmode\edef\@sf{\spacefactor\the\spacefactor}\fi
 {\textfontii g}\@sf\relax\fi}
\let\lbrace\{
\let\rbrace\}
\def\AmSTeX{{\textfontii A\kern-.1667em%
  \lower.5ex\hbox{M}\kern-.125emS}-\TeX}
\message{line and page breaks,}
\def\vmodeerr@#1{\Err@{\string#1\space not allowed between paragraphs}}
\def\mathmodeerr@#1{\Err@{\string#1\space not allowed in math mode}}
\def\linebreak{\RIfM@\mathmodeerr@\linebreak\else
 \ifhmode\unskip\unkern\break\else\vmodeerr@\linebreak\fi\fi}

\newskip\saveskip@
\def\allowlinebreak{\RIfM@\mathmodeerr@\allowlinebreak\else
 \ifhmode\saveskip@\lastskip\unskip
 \allowbreak\ifdim\saveskip@>\z@\hskip\saveskip@\fi
 \else\vmodeerr@\allowlinebreak\fi\fi}
\def\nolinebreak{\RIfM@\mathmodeerr@\nolinebreak\else
 \ifhmode\saveskip@\lastskip\unskip
 \nobreak\ifdim\saveskip@>\z@\hskip\saveskip@\fi
 \else\vmodeerr@\nolinebreak\fi\fi}
\def\newline{\relaxnext@
 \DN@{\RIfM@\expandafter\mathmodeerr@\expandafter\newline\else
  \ifhmode\ifx\next\par\else
  \expandafter\unskip\expandafter\null\expandafter\hfill\expandafter\break\fi
  \else
  \expandafter\vmodeerr@\expandafter\newline\fi\fi}%
 \FN@\next@}
\def\dmatherr@#1{\Err@{\string#1\space not allowed in display math mode}}
\def\nondmatherr@#1{\Err@{\string#1\space not allowed in non-display math
 mode}}
\def\onlydmatherr@#1{\Err@{\string#1\space allowed only in display math mode}}
\def\nonmatherr@#1{\Err@{\string#1\space allowed only in math mode}}
\def\mathbreak{\RIfMIfI@\break\else
 \dmatherr@\mathbreak\fi\else\nonmatherr@\mathbreak\fi}
\def\nomathbreak{\RIfMIfI@\nobreak\else
 \dmatherr@\nomathbreak\fi\else\nonmatherr@\nomathbreak\fi}
\def\allowmathbreak{\RIfMIfI@\allowbreak\else
 \dmatherr@\allowmathbreak\fi\else\nonmatherr@\allowmathbreak\fi}
\def\pagebreak{\RIfM@
 \ifinner\nondmatherr@\pagebreak\else\postdisplaypenalty-\@M\fi
 \else\ifvmode\removelastskip\break\else\vadjust{\break}\fi\fi}
\def\nopagebreak{\RIfM@
 \ifinner\nondmatherr@\nopagebreak\else\postdisplaypenalty\@M\fi
 \else\ifvmode\nobreak\else\vadjust{\nobreak}\fi\fi}
\def\nonvmodeerr@#1{\Err@{\string#1\space not allowed within a paragraph
 or in math}}
\def\vnonvmode@#1#2{\relaxnext@\DNii@{\ifx\next\par\DN@{#1}\else
 \DN@{#2}\fi\next@}%
 \ifvmode\DN@{#1}\else
 \DN@{\FN@\nextii@}\fi\next@}
\def\newpage{\vnonvmode@{\vfill\break}{\nonvmodeerr@\newpage}}
\def\smallpagebreak{\vnonvmode@\smallbreak{\nonvmodeerr@\smallpagebreak}}
\def\medpagebreak{\vnonvmode@\medbreak{\nonvmodeerr@\medpagebreak}}
\def\bigpagebreak{\vnonvmode@\bigbreak{\nonvmodeerr@\bigpagebreak}}
\def\NoBlackBoxes{\global\overfullrule\z@}
\def\BlackBoxes{\global\overfullrule5\p@}
\def\Invalid@#1{\def#1{\Err@{\Invalid@@\string#1}}}
\def\Invalid@@{Invalid use of }
\message{figures,}
\Invalid@\caption
\Invalid@\captionwidth
\newdimen\smallcaptionwidth@
\def\topspace{\mid@false\ins@}
\def\midspace{\mid@true\ins@}
\newif\ifmid@
\def\captionfont@{}
\def\ins@#1{\relaxnext@\allowbreak
 \smallcaptionwidth@\captionwidth@\gdef\thespace@{#1}%
 \DN@{\ifx\next\space@\DN@. {\FN@\nextii@}\else
  \DN@.{\FN@\nextii@}\fi\next@.}%
 \DNii@{\ifx\next\caption\DN@\caption{\FN@\nextiii@}%
  \else\let\next@\nextiv@\fi\next@}%
 \def\nextiv@{\vnonvmode@
  {\ifmid@\expandafter\midinsert\else\expandafter\topinsert\fi
   \vbox to\thespace@{}\endinsert}
  {\ifmid@\nonvmodeerr@\midspace\else\nonvmodeerr@\topspace\fi}}%
 \def\nextiii@{\ifx\next\captionwidth\expandafter\nextv@
  \else\expandafter\nextvi@\fi}%
 \def\nextv@\captionwidth##1##2{\smallcaptionwidth@##1\relax\nextvi@{##2}}%
 \def\nextvi@##1{\def\thecaption@{\captionfont@##1}%
  \DN@{\ifx\next\space@\DN@. {\FN@\nextvii@}\else
   \DN@.{\FN@\nextvii@}\fi\next@.}%
  \FN@\next@}%
 \def\nextvii@{\vnonvmode@
  {\ifmid@\expandafter\midinsert\else
  \expandafter\topinsert\fi\vbox to\thespace@{}\nobreak\smallskip
  \setboxz@h{\noindent\ignorespaces\thecaption@\unskip}%
  \ifdim\wdz@>\smallcaptionwidth@\centerline{\vbox{\hsize\smallcaptionwidth@
   \noindent\ignorespaces\thecaption@\unskip}}%
  \else\centerline{\boxz@}\fi\endinsert}
  {\ifmid@\nonvmodeerr@\midspace
  \else\nonvmodeerr@\topspace\fi}}%
 \FN@\next@}
\message{comments,}
\def\newcodes@{\catcode`\\12\catcode`\{12\catcode`\}12\catcode`\#12%
 \catcode`\%12\relax}
\def\oldcodes@{\catcode`\\0\catcode`\{1\catcode`\}2\catcode`\#6%
 \catcode`\%14\relax}
\def\comment{\newcodes@\endlinechar=10 \comment@}
{\lccode`\0=`\\
\lowercase{\gdef\comment@#1^^J{\comment@@#10endcomment\comment@@@}%
\gdef\comment@@#10endcomment{\FN@\comment@@@}%
\gdef\comment@@@#1\comment@@@{\ifx\next\comment@@@\let\next\comment@
 \else\def\next{\oldcodes@\endlinechar=`\^^M\relax}%
 \fi\next}}}
\def\pr@m@s{\ifx'\next\DN@##1{\prim@s}\else\let\next@\egroup\fi\next@}
\def\prime{{\null\prime@\null}}
\mathchardef\prime@="0230
\let\dsize\displaystyle

\let\ssize\scriptstyle

\message{math spacing,}
\def\,{\RIfM@\mskip\thinmuskip\relax\else\kern.16667em\fi}
\def\!{\RIfM@\mskip-\thinmuskip\relax\else\kern-.16667em\fi}
\let\thinspace\,
\let\negthinspace\!
\def\medspace{\RIfM@\mskip\medmuskip\relax\else\kern.222222em\fi}
\def\negmedspace{\RIfM@\mskip-\medmuskip\relax\else\kern-.222222em\fi}
\def\thickspace{\RIfM@\mskip\thickmuskip\relax\else\kern.27777em\fi}
\let\;\thickspace
\def\negthickspace{\RIfM@\mskip-\thickmuskip\relax\else
 \kern-.27777em\fi}
\atdef@,{\RIfM@\mskip.1\thinmuskip\else\leavevmode\null,\fi}
\atdef@!{\RIfM@\mskip-.1\thinmuskip\else\leavevmode\null!\fi}
\atdef@.{\RIfM@&&\else\leavevmode.\spacefactor3000 \fi}
\def\and{\DOTSB\;\mathchar"3026 \;}

\message{fractions,}
\def\frac#1#2{{#1\over#2}}

\newdimen\ex@
\ex@.2326ex
\Invalid@\thickness
\def\thickfrac{\relaxnext@
 \DN@{\ifx\next\thickness\let\next@\nextii@\else
 \DN@{\nextii@\thickness1}\fi\next@}%
 \DNii@\thickness##1##2##3{{##2\above##1\ex@##3}}%
 \FN@\next@}

\def\thickfracwithdelims#1#2{\relaxnext@\def\ldelim@{#1}\def\rdelim@{#2}%
 \DN@{\ifx\next\thickness\let\next@\nextii@\else
 \DN@{\nextii@\thickness1}\fi\next@}%
 \DNii@\thickness##1##2##3{{##2\abovewithdelims
 \ldelim@\rdelim@##1\ex@##3}}%
 \FN@\next@}
\def\binom#1#2{{#1\choose#2}}

\def\:{\nobreak\hskip.1111em\mathpunct{}\nonscript\mkern-\thinmuskip{:}\hskip
 .3333emplus.0555em\relax}
\def\snug{\unskip\kern-\mathsurround}
\message{smash commands,}
\def\topsmash{\top@true\bot@false\smash@}
\def\botsmash{\top@false\bot@true\smash@}
\newif\iftop@
\newif\ifbot@
\def\smash{\top@true\bot@true\smash@}
\def\smash@{\RIfM@\expandafter\mathpalette\expandafter\mathsm@sh\else
 \expandafter\makesm@sh\fi}
\def\finsm@sh{\iftop@\ht\z@\z@\fi\ifbot@\dp\z@\z@\fi\leavevmode\boxz@}
\message{large operator symbols,}
\def\LimitsOnSums{\global\let\slimits@\displaylimits}
\def\NoLimitsOnSums{\global\let\slimits@\nolimits}
\LimitsOnSums
\mathchardef\coprod@="1360       \def\coprod{\DOTSB\coprod@\slimits@}
\mathchardef\bigvee@="1357       \def\bigvee{\DOTSB\bigvee@\slimits@}
\mathchardef\bigwedge@="1356     \def\bigwedge{\DOTSB\bigwedge@\slimits@}
\mathchardef\biguplus@="1355     \def\biguplus{\DOTSB\biguplus@\slimits@}
\mathchardef\bigcap@="1354       \def\bigcap{\DOTSB\bigcap@\slimits@}
\mathchardef\bigcup@="1353       \def\bigcup{\DOTSB\bigcup@\slimits@}
\mathchardef\prod@="1351         \def\prod{\DOTSB\prod@\slimits@}
\mathchardef\sum@="1350          \def\sum{\DOTSB\sum@\slimits@}
\mathchardef\bigotimes@="134E    \def\bigotimes{\DOTSB\bigotimes@\slimits@}
\mathchardef\bigoplus@="134C     \def\bigoplus{\DOTSB\bigoplus@\slimits@}
\mathchardef\bigodot@="134A      \def\bigodot{\DOTSB\bigodot@\slimits@}
\mathchardef\bigsqcup@="1346     \def\bigsqcup{\DOTSB\bigsqcup@\slimits@}
\message{integrals,}
\def\LimitsOnInts{\global\let\ilimits@\displaylimits}
\def\NoLimitsOnInts{\global\let\ilimits@\nolimits}
\NoLimitsOnInts
\def\int{\DOTSI\intop\ilimits@}
\def\oint{\DOTSI\ointop\ilimits@}
\def\intic@{\mathchoice{\hskip.5em}{\hskip.4em}{\hskip.4em}{\hskip.4em}}
\def\negintic@{\mathchoice
 {\hskip-.5em}{\hskip-.4em}{\hskip-.4em}{\hskip-.4em}}
\def\intkern@{\mathchoice{\!\!\!}{\!\!}{\!\!}{\!\!}}
\def\intdots@{\mathchoice{\plaincdots@}
 {{\cdotp}\mkern1.5mu{\cdotp}\mkern1.5mu{\cdotp}}
 {{\cdotp}\mkern1mu{\cdotp}\mkern1mu{\cdotp}}
 {{\cdotp}\mkern1mu{\cdotp}\mkern1mu{\cdotp}}}
\newcount\intno@
\def\iint{\DOTSI\intno@\tw@\FN@\ints@}
\def\iiint{\DOTSI\intno@\thr@@\FN@\ints@}
\def\iiiint{\DOTSI\intno@4 \FN@\ints@}
\def\idotsint{\DOTSI\intno@\z@\FN@\ints@}
\def\ints@{\findlimits@\ints@@}
\newif\iflimtoken@
\newif\iflimits@
\def\findlimits@{\limtoken@true\ifx\next\limits\limits@true
 \else\ifx\next\nolimits\limits@false\else
 \limtoken@false\ifx\ilimits@\nolimits\limits@false\else
 \ifinner\limits@false\else\limits@true\fi\fi\fi\fi}
\def\multint@{\int\ifnum\intno@=\z@\intdots@                                
 \else\intkern@\fi                                                          
 \ifnum\intno@>\tw@\int\intkern@\fi                                         
 \ifnum\intno@>\thr@@\int\intkern@\fi                                       
 \int}                                                                      
\def\multintlimits@{\intop\ifnum\intno@=\z@\intdots@\else\intkern@\fi
 \ifnum\intno@>\tw@\intop\intkern@\fi
 \ifnum\intno@>\thr@@\intop\intkern@\fi\intop}
\def\ints@@{\iflimtoken@                                                    
 \def\ints@@@{\iflimits@\negintic@\mathop{\intic@\multintlimits@}\limits    
  \else\multint@\nolimits\fi                                                
  \eat@}                                                                    
 \else                                                                      
 \def\ints@@@{\iflimits@\negintic@
  \mathop{\intic@\multintlimits@}\limits\else
  \multint@\nolimits\fi}\fi\ints@@@}
\def\LimitsOnNames{\global\let\nlimits@\displaylimits}
\def\NoLimitsOnNames{\global\let\nlimits@\nolimits@}
\LimitsOnNames
\def\nolimits@{\relaxnext@
 \DN@{\ifx\next\limits\DN@\limits{\nolimits}\else
  \let\next@\nolimits\fi\next@}%
 \FN@\next@}
\message{operator names,}
\def\newmcodes@{\mathcode`\'"27\mathcode`\*"2A\mathcode`\."613A%
 \mathcode`\-"2D\mathcode`\/"2F\mathcode`\:"603A }
\def\operatorname#1{\mathop{\newmcodes@\kern\z@\fam\z@#1}\nolimits@}
\def\operatornamewithlimits#1{\mathop{\newmcodes@\kern\z@\fam\z@#1}\nlimits@}
\def\qopname@#1{\mathop{\fam\z@#1}\nolimits@}
\def\qopnamewl@#1{\mathop{\fam\z@#1}\nlimits@}
\def\arccos{\qopname@{arccos}}
\def\arcsin{\qopname@{arcsin}}
\def\arctan{\qopname@{arctan}}
\def\arg{\qopname@{arg}}
\def\cos{\qopname@{cos}}
\def\cosh{\qopname@{cosh}}
\def\cot{\qopname@{cot}}
\def\coth{\qopname@{coth}}
\def\csc{\qopname@{csc}}
\def\deg{\qopname@{deg}}
\def\det{\qopnamewl@{det}}
\def\dim{\qopname@{dim}}
\def\exp{\qopname@{exp}}
\def\gcd{\qopnamewl@{gcd}}
\def\hom{\qopname@{hom}}
\def\inf{\qopnamewl@{inf}}
\def\injlim{\qopnamewl@{inj\,lim}}
\def\ker{\qopname@{ker}}
\def\lg{\qopname@{lg}}
\def\lim{\qopnamewl@{lim}}
\def\liminf{\qopnamewl@{lim\,inf}}
\def\limsup{\qopnamewl@{lim\,sup}}
\def\ln{\qopname@{ln}}
\def\log{\qopname@{log}}
\def\max{\qopnamewl@{max}}
\def\min{\qopnamewl@{min}}
\def\Pr{\qopnamewl@{Pr}}
\def\projlim{\qopnamewl@{proj\,lim}}
\def\sec{\qopname@{sec}}
\def\sin{\qopname@{sin}}
\def\sinh{\qopname@{sinh}}
\def\sup{\qopnamewl@{sup}}
\def\tan{\qopname@{tan}}
\def\tanh{\qopname@{tanh}}
\def\varinjlim{\mathop{\vtop{\ialign{##\crcr
 \hfil\rm lim\hfil\crcr\noalign{\nointerlineskip}\rightarrowfill\crcr
 \noalign{\nointerlineskip\kern-\ex@}\crcr}}}}
\def\varprojlim{\mathop{\vtop{\ialign{##\crcr
 \hfil\rm lim\hfil\crcr\noalign{\nointerlineskip}\leftarrowfill\crcr
 \noalign{\nointerlineskip\kern-\ex@}\crcr}}}}
\def\varliminf{\mathop{\underline{\vrule height\z@ depth.2exwidth\z@
 \hbox{\rm lim}}}}

\newdimen\buffer@
\buffer@\fontdimen13 \tenex
\newdimen\buffer
\buffer\buffer@

\def\ResetBuffer{\fontdimen13 \tenex\buffer@\global\buffer\buffer@}
\def\shave#1{\mathop{\hbox{$\m@th\fontdimen13 \tenex\z@                     
 \displaystyle{#1}$}}\fontdimen13 \tenex\buffer}

\message{multilevel sub/superscripts,}
\Invalid@\\
\def\Let@{\relax\iffalse{\fi\let\\=\cr\iffalse}\fi}
\Invalid@\vspace
\def\vspace@{\def\vspace##1{\crcr\noalign{\vskip##1\relax}}}
\def\multilimits@{\bgroup\vspace@\Let@
 \baselineskip\fontdimen10 \scriptfont\tw@
 \advance\baselineskip\fontdimen12 \scriptfont\tw@
 \lineskip\thr@@\fontdimen8 \scriptfont\thr@@
 \lineskiplimit\lineskip
 \vbox\bgroup\ialign\bgroup\hfil$\m@th\scriptstyle{##}$\hfil\crcr}
\def\Sb{_\multilimits@}
\def\endSb{\crcr\egroup\egroup\egroup}
\def\Sp{^\multilimits@}

\def\spreadlines#1{\RIfMIfI@\onlydmatherr@\spreadlines\else
 \openup#1\relax\fi\else\onlydmatherr@\spreadlines\fi}
\def\Mathstrut@{\copy\Mathstrutbox@}
\newbox\Mathstrutbox@
\setbox\Mathstrutbox@\null
\setboxz@h{$\m@th($}
\ht\Mathstrutbox@\ht\z@
\dp\Mathstrutbox@\dp\z@
\message{matrices,}
\newdimen\spreadmlines@
\def\spreadmatrixlines#1{\RIfMIfI@
 \onlydmatherr@\spreadmatrixlines\else
 \spreadmlines@#1\relax\fi\else\onlydmatherr@\spreadmatrixlines\fi}
\def\matrix{\null\,\vcenter\bgroup\Let@\vspace@
 \normalbaselines\openup\spreadmlines@\ialign
 \bgroup\hfil$\m@th##$\hfil&&\quad\hfil$\m@th##$\hfil\crcr
 \Mathstrut@\crcr\noalign{\kern-\baselineskip}}
\def\endmatrix{\crcr\Mathstrut@\crcr\noalign{\kern-\baselineskip}\egroup
 \egroup\,}
\def\format{\crcr\egroup\iffalse{\fi\ifnum`}=0 \fi\format@}
\newtoks\hashtoks@
\hashtoks@{#}
\def\format@#1\\{\def\preamble@{#1}%
 \def\l{$\m@th\the\hashtoks@$\hfil}%
 \def\c{\hfil$\m@th\the\hashtoks@$\hfil}%
 \def\r{\hfil$\m@th\the\hashtoks@$}%
 \edef\preamble@@{\preamble@}\ifnum`{=0 \fi\iffalse}\fi
 \ialign\bgroup\span\preamble@@\crcr}
\def\smallmatrix{\null\,\vcenter\bgroup\vspace@\Let@
 \baselineskip9\ex@\lineskip\ex@
 \ialign\bgroup\hfil$\m@th\scriptstyle{##}$\hfil&&\thickspace\hfil
 $\m@th\scriptstyle{##}$\hfil\crcr}
\def\endsmallmatrix{\crcr\egroup\egroup\,}

\newmuskip\dotsspace@
\dotsspace@1.5mu
\def\strip@#1 {#1}
\def\spacehdots#1\for#2{\multispan{#2}\xleaders
 \hbox{$\m@th\mkern\strip@#1 \dotsspace@.\mkern\strip@#1 \dotsspace@$}\hfill}
\def\hdotsfor#1{\spacehdots\@ne\for{#1}}
\def\multispan@#1{\omit\mscount#1\unskip\loop\ifnum\mscount>\@ne\sp@n\repeat}
\def\spaceinnerhdots#1\for#2\after#3{\multispan@{\strip@#2 }#3\xleaders
 \hbox{$\m@th\mkern\strip@#1 \dotsspace@.\mkern\strip@#1 \dotsspace@$}\hfill}
\def\innerhdotsfor#1\after#2{\spaceinnerhdots\@ne\for#1\after{#2}}
\def\cases{\bgroup\spreadmlines@\jot\left\{\,\matrix\format\l&\quad\l\\}
\def\endcases{\endmatrix\right.\egroup}
\message{multiline displays,}
\newif\ifinany@
\newif\ifinalign@
\newif\ifingather@
\def\strut@{\copy\strutbox@}
\newbox\strutbox@
\setbox\strutbox@\hbox{\vrule height8\p@ depth3\p@ width\z@}
\def\topaligned{\null\,\vtop\aligned@}
\def\botaligned{\null\,\vbox\aligned@}
\def\aligned{\null\,\vcenter\aligned@}
\def\aligned@{\bgroup\vspace@\Let@
 \ifinany@\else\openup\jot\fi\ialign
 \bgroup\hfil\strut@$\m@th\displaystyle{##}$&
 $\m@th\displaystyle{{}##}$\hfil\crcr}
\def\endaligned{\crcr\egroup\egroup}

\def\alignedat#1{\null\,\vcenter\bgroup\doat@{#1}\vspace@\Let@
 \ifinany@\else\openup\jot\fi\ialign\bgroup\span\preamble@@\crcr}
\newcount\atcount@
\def\doat@#1{\toks@{\hfil\strut@$\m@th
 \displaystyle{\the\hashtoks@}$&$\m@th\displaystyle
 {{}\the\hashtoks@}$\hfil}
 \atcount@#1\relax\advance\atcount@\m@ne                                    
 \loop\ifnum\atcount@>\z@\toks@=\expandafter{\the\toks@&\hfil$\m@th
 \displaystyle{\the\hashtoks@}$&$\m@th
 \displaystyle{{}\the\hashtoks@}$\hfil}\advance
  \atcount@\m@ne\repeat                                                     
 \xdef\preamble@{\the\toks@}\xdef\preamble@@{\preamble@}}

\def\gathered{\null\,\vcenter\bgroup\vspace@\Let@
 \ifinany@\else\openup\jot\fi\ialign
 \bgroup\hfil\strut@$\m@th\displaystyle{##}$\hfil\crcr}
\def\endgathered{\crcr\egroup\egroup}
\newif\iftagsleft@
\def\TagsOnLeft{\global\tagsleft@true}
\def\TagsOnRight{\global\tagsleft@false}
\TagsOnLeft
\newif\ifmathtags@
\def\TagsAsMath{\global\mathtags@true}
\def\TagsAsText{\global\mathtags@false}
\TagsAsText
\def\tagform@#1{\hbox{\rm(\ignorespaces#1\unskip)}}
\def\thetag{\leavevmode\tagform@}
\def\tag#1$${\iftagsleft@\leqno\else\eqno\fi                                
 \maketag@#1\maketag@                                                       
 $$}                                                                        
\def\maketag@{\FN@\maketag@@}
\def\maketag@@{\ifx\next"\expandafter\maketag@@@\else\expandafter\maketag@@@@
 \fi}
\def\maketag@@@"#1"#2\maketag@{\hbox{\rm#1}}                                
\def\maketag@@@@#1\maketag@{\ifmathtags@\tagform@{$\m@th#1$}\else
 \tagform@{#1}\fi}
\interdisplaylinepenalty\@M
\def\allowdisplaybreaks{\RIfMIfI@
 \onlydmatherr@\allowdisplaybreaks\else
 \interdisplaylinepenalty\z@\fi\else\onlydmatherr@\allowdisplaybreaks\fi}
\Invalid@\allowdisplaybreak
\Invalid@\displaybreak
\Invalid@\intertext
\def\allowdisplaybreak@{\def\allowdisplaybreak{\crcr\noalign{\allowbreak}}}
\def\displaybreak@{\def\displaybreak{\crcr\noalign{\break}}}
\def\intertext@{\def\intertext##1{\crcr\noalign{%
 \penalty\postdisplaypenalty \vskip\belowdisplayskip
 \vbox{\normalbaselines\noindent##1}%
 \penalty\predisplaypenalty \vskip\abovedisplayskip}}}
\newskip\centering@
\centering@\z@ plus\@m\p@
\def\align{\relax\ifingather@\DN@{\csname align (in
  \string\gather)\endcsname}\else
 \ifmmode\ifinner\DN@{\onlydmatherr@\align}\else
  \let\next@\align@\fi
 \else\DN@{\onlydmatherr@\align}\fi\fi\next@}
\newhelp\andhelp@
{An extra & here is so disastrous that you should probably exit^^J
and fix things up.}
\newif\iftag@
\newcount\and@
\def\align@{\inalign@true\inany@true
 \vspace@\allowdisplaybreak@\displaybreak@\intertext@
 \def\tag{\global\tag@true\ifnum\and@=\z@\DN@{&&}\else
          \DN@{&}\fi\next@}%
 \iftagsleft@\DN@{\csname align \endcsname}\else
  \DN@{\csname align \space\endcsname}\fi\next@}
\def\Tag@{\iftag@\else\errhelp\andhelp@\err@{Extra & on this line}\fi}
\newdimen\lwidth@
\newdimen\rwidth@
\newdimen\maxlwidth@
\newdimen\maxrwidth@
\newdimen\totwidth@
\def\measure@#1\endalign{\lwidth@\z@\rwidth@\z@\maxlwidth@\z@\maxrwidth@\z@
 \global\and@\z@                                                            
 \setbox@ne\vbox                                                            
  {\everycr{\noalign{\global\tag@false\global\and@\z@}}\Let@                
  \halign{\setboxz@h{$\m@th\displaystyle{\@lign##}$}
   \global\lwidth@\wdz@                                                     
   \ifdim\lwidth@>\maxlwidth@\global\maxlwidth@\lwidth@\fi                  
   \global\advance\and@\@ne                                                 
   &\setboxz@h{$\m@th\displaystyle{{}\@lign##}$}\global\rwidth@\wdz@        
   \ifdim\rwidth@>\maxrwidth@\global\maxrwidth@\rwidth@\fi                  
   \global\advance\and@\@ne                                                
   &\Tag@
   \eat@{##}\crcr#1\crcr}}
 \totwidth@\maxlwidth@\advance\totwidth@\maxrwidth@}                       
\def\displ@y@{\global\dt@ptrue\openup\jot
 \everycr{\noalign{\global\tag@false\global\and@\z@\ifdt@p\global\dt@pfalse
 \vskip-\lineskiplimit\vskip\normallineskiplimit\else
 \penalty\interdisplaylinepenalty\fi}}}
\def\black@#1{\noalign{\ifdim#1>\displaywidth
 \dimen@\prevdepth\nointerlineskip                                          
 \vskip-\ht\strutbox@\vskip-\dp\strutbox@                                   
 \vbox{\noindent\hbox to#1{\strut@\hfill}}
 \prevdepth\dimen@                                                          
 \fi}}
\expandafter\def\csname align \space\endcsname#1\endalign
 {\measure@#1\endalign\global\and@\z@                                       
 \ifingather@\everycr{\noalign{\global\and@\z@}}\else\displ@y@\fi           
 \Let@\tabskip\centering@                                                   
 \halign to\displaywidth
  {\hfil\strut@\setboxz@h{$\m@th\displaystyle{\@lign##}$}
  \global\lwidth@\wdz@\boxz@\global\advance\and@\@ne                        
  \tabskip\z@skip                                                           
  &\setboxz@h{$\m@th\displaystyle{{}\@lign##}$}
  \global\rwidth@\wdz@\boxz@\hfill\global\advance\and@\@ne                  
  \tabskip\centering@                                                       
  &\setboxz@h{\@lign\strut@\maketag@##\maketag@}
  \dimen@\displaywidth\advance\dimen@-\totwidth@
  \divide\dimen@\tw@\advance\dimen@\maxrwidth@\advance\dimen@-\rwidth@     
  \ifdim\dimen@<\tw@\wdz@\llap{\vtop{\normalbaselines\null\boxz@}}
  \else\llap{\boxz@}\fi                                                    
  \tabskip\z@skip                                                          
  \crcr#1\crcr                                                             
  \black@\totwidth@}}                                                      
\newdimen\lineht@
\expandafter\def\csname align \endcsname#1\endalign{\measure@#1\endalign
 \global\and@\z@
 \ifdim\totwidth@>\displaywidth\let\displaywidth@\totwidth@\else
  \let\displaywidth@\displaywidth\fi                                        
 \ifingather@\everycr{\noalign{\global\and@\z@}}\else\displ@y@\fi
 \Let@\tabskip\centering@\halign to\displaywidth
  {\hfil\strut@\setboxz@h{$\m@th\displaystyle{\@lign##}$}%
  \global\lwidth@\wdz@\global\lineht@\ht\z@                                 
  \boxz@\global\advance\and@\@ne
  \tabskip\z@skip&\setboxz@h{$\m@th\displaystyle{{}\@lign##}$}%
  \global\rwidth@\wdz@\ifdim\ht\z@>\lineht@\global\lineht@\ht\z@\fi         
  \boxz@\hfil\global\advance\and@\@ne
  \tabskip\centering@&\kern-\displaywidth@                                  
  \setboxz@h{\@lign\strut@\maketag@##\maketag@}%
  \dimen@\displaywidth\advance\dimen@-\totwidth@
  \divide\dimen@\tw@\advance\dimen@\maxlwidth@\advance\dimen@-\lwidth@
  \ifdim\dimen@<\tw@\wdz@
   \rlap{\vbox{\normalbaselines\boxz@\vbox to\lineht@{}}}\else
   \rlap{\boxz@}\fi
  \tabskip\displaywidth@\crcr#1\crcr\black@\totwidth@}}
\expandafter\def\csname align (in \string\gather)\endcsname
  #1\endalign{\vcenter{\align@#1\endalign}}
\Invalid@\endalign
\newif\ifxat@
\def\alignat{\RIfMIfI@\DN@{\onlydmatherr@\alignat}\else
 \DN@{\csname alignat \endcsname}\fi\else
 \DN@{\onlydmatherr@\alignat}\fi\next@}
\newif\ifmeasuring@
\newbox\savealignat@
\expandafter\def\csname alignat \endcsname#1#2\endalignat                   
 {\inany@true\xat@false
 \def\tag{\global\tag@true\count@#1\relax\multiply\count@\tw@
  \xdef\tag@{}\loop\ifnum\count@>\and@\xdef\tag@{&\tag@}\advance\count@\m@ne
  \repeat\tag@}%
 \vspace@\allowdisplaybreak@\displaybreak@\intertext@
 \displ@y@\measuring@true                                                   
 \setbox\savealignat@\hbox{$\m@th\displaystyle\Let@
  \attag@{#1}
  \vbox{\halign{\span\preamble@@\crcr#2\crcr}}$}%
 \measuring@false                                                           
 \Let@\attag@{#1}
 \tabskip\centering@\halign to\displaywidth
  {\span\preamble@@\crcr#2\crcr                                             
  \black@{\wd\savealignat@}}}                                               
\Invalid@\endalignat
\def\xalignat{\RIfMIfI@
 \DN@{\onlydmatherr@\xalignat}\else
 \DN@{\csname xalignat \endcsname}\fi\else
 \DN@{\onlydmatherr@\xalignat}\fi\next@}
\expandafter\def\csname xalignat \endcsname#1#2\endxalignat
 {\inany@true\xat@true
 \def\tag{\global\tag@true\def\tag@{}\count@#1\relax\multiply\count@\tw@
  \loop\ifnum\count@>\and@\xdef\tag@{&\tag@}\advance\count@\m@ne\repeat\tag@}%
 \vspace@\allowdisplaybreak@\displaybreak@\intertext@
 \displ@y@\measuring@true\setbox\savealignat@\hbox{$\m@th\displaystyle\Let@
 \attag@{#1}\vbox{\halign{\span\preamble@@\crcr#2\crcr}}$}%
 \measuring@false\Let@
 \attag@{#1}\tabskip\centering@\halign to\displaywidth
 {\span\preamble@@\crcr#2\crcr\black@{\wd\savealignat@}}}
\def\attag@#1{\let\Maketag@\maketag@\let\TAG@\Tag@                          
 \let\Tag@=0\let\maketag@=0
 \ifmeasuring@\def\llap@##1{\setboxz@h{##1}\hbox to\tw@\wdz@{}}%
  \def\rlap@##1{\setboxz@h{##1}\hbox to\tw@\wdz@{}}\else
  \let\llap@\llap\let\rlap@\rlap\fi                                         
 \toks@{\hfil\strut@$\m@th\displaystyle{\@lign\the\hashtoks@}$\tabskip\z@skip
  \global\advance\and@\@ne&$\m@th\displaystyle{{}\@lign\the\hashtoks@}$\hfil
  \ifxat@\tabskip\centering@\fi\global\advance\and@\@ne}
 \iftagsleft@
  \toks@@{\tabskip\centering@&\Tag@\kern-\displaywidth
   \rlap@{\@lign\maketag@\the\hashtoks@\maketag@}%
   \global\advance\and@\@ne\tabskip\displaywidth}\else
  \toks@@{\tabskip\centering@&\Tag@\llap@{\@lign\maketag@
   \the\hashtoks@\maketag@}\global\advance\and@\@ne\tabskip\z@skip}\fi      
 \atcount@#1\relax\advance\atcount@\m@ne
 \loop\ifnum\atcount@>\z@
 \toks@=\expandafter{\the\toks@&\hfil$\m@th\displaystyle{\@lign
  \the\hashtoks@}$\global\advance\and@\@ne
  \tabskip\z@skip&$\m@th\displaystyle{{}\@lign\the\hashtoks@}$\hfil\ifxat@
  \tabskip\centering@\fi\global\advance\and@\@ne}\advance\atcount@\m@ne
 \repeat                                                                    
 \xdef\preamble@{\the\toks@\the\toks@@}
 \xdef\preamble@@{\preamble@}
 \let\maketag@\Maketag@\let\Tag@\TAG@}                                      
\Invalid@\endxalignat
\def\xxalignat{\RIfMIfI@
 \DN@{\onlydmatherr@\xxalignat}\else\DN@{\csname xxalignat
  \endcsname}\fi\else
 \DN@{\onlydmatherr@\xxalignat}\fi\next@}
\expandafter\def\csname xxalignat \endcsname#1#2\endxxalignat{\inany@true
 \vspace@\allowdisplaybreak@\displaybreak@\intertext@
 \displ@y\setbox\savealignat@\hbox{$\m@th\displaystyle\Let@
 \xxattag@{#1}\vbox{\halign{\span\preamble@@\crcr#2\crcr}}$}%
 \Let@\xxattag@{#1}\tabskip\z@skip\halign to\displaywidth
 {\span\preamble@@\crcr#2\crcr\black@{\wd\savealignat@}}}
\def\xxattag@#1{\toks@{\tabskip\z@skip\hfil\strut@
 $\m@th\displaystyle{\the\hashtoks@}$&%
 $\m@th\displaystyle{{}\the\hashtoks@}$\hfil\tabskip\centering@&}%
 \atcount@#1\relax\advance\atcount@\m@ne\loop\ifnum\atcount@>\z@
 \toks@=\expandafter{\the\toks@&\hfil$\m@th\displaystyle{\the\hashtoks@}$%
  \tabskip\z@skip&$\m@th\displaystyle{{}\the\hashtoks@}$\hfil
  \tabskip\centering@}\advance\atcount@\m@ne\repeat
 \xdef\preamble@{\the\toks@\tabskip\z@skip}\xdef\preamble@@{\preamble@}}
\Invalid@\endxxalignat
\newdimen\gwidth@
\newdimen\gmaxwidth@
\def\gmeasure@#1\endgather{\gwidth@\z@\gmaxwidth@\z@\setbox@ne\vbox{\Let@
 \halign{\setboxz@h{$\m@th\displaystyle{##}$}\global\gwidth@\wdz@
 \ifdim\gwidth@>\gmaxwidth@\global\gmaxwidth@\gwidth@\fi
 &\eat@{##}\crcr#1\crcr}}}
\def\gather{\RIfMIfI@\DN@{\onlydmatherr@\gather}\else
 \ingather@true\inany@true\def\tag{&}%
 \vspace@\allowdisplaybreak@\displaybreak@\intertext@
 \displ@y\Let@
 \iftagsleft@\DN@{\csname gather \endcsname}\else
  \DN@{\csname gather \space\endcsname}\fi\fi
 \else\DN@{\onlydmatherr@\gather}\fi\next@}
\expandafter\def\csname gather \space\endcsname#1\endgather
 {\gmeasure@#1\endgather\tabskip\centering@
 \halign to\displaywidth{\hfil\strut@\setboxz@h{$\m@th\displaystyle{##}$}%
 \global\gwidth@\wdz@\boxz@\hfil&
 \setboxz@h{\strut@{\maketag@##\maketag@}}%
 \dimen@\displaywidth\advance\dimen@-\gwidth@
 \ifdim\dimen@>\tw@\wdz@\llap{\boxz@}\else
 \llap{\vtop{\normalbaselines\null\boxz@}}\fi
 \tabskip\z@skip\crcr#1\crcr\black@\gmaxwidth@}}
\newdimen\glineht@
\expandafter\def\csname gather \endcsname#1\endgather{\gmeasure@#1\endgather
 \ifdim\gmaxwidth@>\displaywidth\let\gdisplaywidth@\gmaxwidth@\else
 \let\gdisplaywidth@\displaywidth\fi\tabskip\centering@\halign to\displaywidth
 {\hfil\strut@\setboxz@h{$\m@th\displaystyle{##}$}%
 \global\gwidth@\wdz@\global\glineht@\ht\z@\boxz@\hfil&\kern-\gdisplaywidth@
 \setboxz@h{\strut@{\maketag@##\maketag@}}%
 \dimen@\displaywidth\advance\dimen@-\gwidth@
 \ifdim\dimen@>\tw@\wdz@\rlap{\boxz@}\else
 \rlap{\vbox{\normalbaselines\boxz@\vbox to\glineht@{}}}\fi
 \tabskip\gdisplaywidth@\crcr#1\crcr\black@\gmaxwidth@}}
\newif\ifctagsplit@
\def\CenteredTagsOnSplits{\global\ctagsplit@true}
\def\TopOrBottomTagsOnSplits{\global\ctagsplit@false}
\TopOrBottomTagsOnSplits
\def\split{\relax\ifinany@\let\next@\insplit@\else
 \ifmmode\ifinner\def\next@{\onlydmatherr@\split}\else
 \let\next@\outsplit@\fi\else
 \def\next@{\onlydmatherr@\split}\fi\fi\next@}
\def\insplit@{\global\setbox\z@\vbox\bgroup\vspace@\Let@\ialign\bgroup
 \hfil\strut@$\m@th\displaystyle{##}$&$\m@th\displaystyle{{}##}$\hfill\crcr}
\def\endsplit{\crcr\egroup\egroup\iftagsleft@\expandafter\lendsplit@\else
 \expandafter\rendsplit@\fi}
\def\rendsplit@{\global\setbox9 \vbox
 {\unvcopy\z@\global\setbox8 \lastbox\unskip}
 \setbox@ne\hbox{\unhcopy8 \unskip\global\setbox\tw@\lastbox
 \unskip\global\setbox\thr@@\lastbox}
 \global\setbox7 \hbox{\unhbox\tw@\unskip}
 \ifinalign@\ifctagsplit@                                                   
  \gdef\split@{\hbox to\wd\thr@@{}&
   \vcenter{\vbox{\moveleft\wd\thr@@\boxz@}}}
 \else\gdef\split@{&\vbox{\moveleft\wd\thr@@\box9}\crcr
  \box\thr@@&\box7}\fi                                                      
 \else                                                                      
  \ifctagsplit@\gdef\split@{\vcenter{\boxz@}}\else
  \gdef\split@{\box9\crcr\hbox{\box\thr@@\box7}}\fi
 \fi
 \split@}                                                                   
\def\lendsplit@{\global\setbox9\vtop{\unvcopy\z@}
 \setbox@ne\vbox{\unvcopy\z@\global\setbox8\lastbox}
 \setbox@ne\hbox{\unhcopy8\unskip\setbox\tw@\lastbox
  \unskip\global\setbox\thr@@\lastbox}
 \ifinalign@\ifctagsplit@                                                   
  \gdef\split@{\hbox to\wd\thr@@{}&
  \vcenter{\vbox{\moveleft\wd\thr@@\box9}}}
  \else                                                                     
  \gdef\split@{\hbox to\wd\thr@@{}&\vbox{\moveleft\wd\thr@@\box9}}\fi
 \else
  \ifctagsplit@\gdef\split@{\vcenter{\box9}}\else
  \gdef\split@{\box9}\fi
 \fi\split@}
\def\outsplit@#1$${\align\insplit@#1\endalign$$}
\newdimen\multlinegap@
\multlinegap@1em
\newdimen\multlinetaggap@
\multlinetaggap@1em
\def\MultlineGap#1{\global\multlinegap@#1\relax}
\def\multlinegap#1{\RIfMIfI@\onlydmatherr@\multlinegap\else
 \multlinegap@#1\relax\fi\else\onlydmatherr@\multlinegap\fi}
\def\nomultlinegap{\multlinegap{\z@}}
\def\multline{\RIfMIfI@
 \DN@{\onlydmatherr@\multline}\else
 \DN@{\multline@}\fi\else
 \DN@{\onlydmatherr@\multline}\fi\next@}
\newif\iftagin@
\def\tagin@#1{\tagin@false\in@\tag{#1}\ifin@\tagin@true\fi}
\def\multline@#1$${\inany@true\vspace@\allowdisplaybreak@\displaybreak@
 \tagin@{#1}\iftagsleft@\DN@{\multline@l#1$$}\else
 \DN@{\multline@r#1$$}\fi\next@}
\newdimen\mwidth@
\def\rmmeasure@#1\endmultline{%
 \def\shoveleft##1{##1}\def\shoveright##1{##1}
 \setbox@ne\vbox{\Let@\halign{\setboxz@h
  {$\m@th\@lign\displaystyle{}##$}\global\mwidth@\wdz@
  \crcr#1\crcr}}}
\newdimen\mlineht@
\newif\ifzerocr@
\newif\ifonecr@
\def\lmmeasure@#1\endmultline{\global\zerocr@true\global\onecr@false
 \everycr{\noalign{\ifonecr@\global\onecr@false\fi
  \ifzerocr@\global\zerocr@false\global\onecr@true\fi}}
  \def\shoveleft##1{##1}\def\shoveright##1{##1}%
 \setbox@ne\vbox{\Let@\halign{\setboxz@h
  {$\m@th\@lign\displaystyle{}##$}\ifonecr@\global\mwidth@\wdz@
  \global\mlineht@\ht\z@\fi\crcr#1\crcr}}}
\newbox\mtagbox@
\newdimen\ltwidth@
\newdimen\rtwidth@
\def\multline@l#1$${\iftagin@\DN@{\lmultline@@#1$$}\else
 \DN@{\setbox\mtagbox@\null\ltwidth@\z@\rtwidth@\z@
  \lmultline@@@#1$$}\fi\next@}
\def\lmultline@@#1\endmultline\tag#2$${%
 \setbox\mtagbox@\hbox{\maketag@#2\maketag@}
 \lmmeasure@#1\endmultline\dimen@\mwidth@\advance\dimen@\wd\mtagbox@
 \advance\dimen@\multlinetaggap@                                            
 \ifdim\dimen@>\displaywidth\ltwidth@\z@\else\ltwidth@\wd\mtagbox@\fi       
 \lmultline@@@#1\endmultline$$}
\def\lmultline@@@{\displ@y
 \def\shoveright##1{##1\hfilneg\hskip\multlinegap@}%
 \def\shoveleft##1{\setboxz@h{$\m@th\displaystyle{}##1$}%
  \setbox@ne\hbox{$\m@th\displaystyle##1$}%
  \hfilneg
  \iftagin@
   \ifdim\ltwidth@>\z@\hskip\ltwidth@\hskip\multlinetaggap@\fi
  \else\hskip\multlinegap@\fi\hskip.5\wd@ne\hskip-.5\wdz@##1}
  \halign\bgroup\Let@\hbox to\displaywidth
   {\strut@$\m@th\displaystyle\hfil{}##\hfil$}\crcr
   \hfilneg                                                                 
   \iftagin@                                                                
    \ifdim\ltwidth@>\z@                                                     
     \box\mtagbox@\hskip\multlinetaggap@                                    
    \else
     \rlap{\vbox{\normalbaselines\hbox{\strut@\box\mtagbox@}%
     \vbox to\mlineht@{}}}\fi                                               
   \else\hskip\multlinegap@\fi}                                             
\def\multline@r#1$${\iftagin@\DN@{\rmultline@@#1$$}\else
 \DN@{\setbox\mtagbox@\null\ltwidth@\z@\rtwidth@\z@
  \rmultline@@@#1$$}\fi\next@}
\def\rmultline@@#1\endmultline\tag#2$${\ltwidth@\z@
 \setbox\mtagbox@\hbox{\maketag@#2\maketag@}%
 \rmmeasure@#1\endmultline\dimen@\mwidth@\advance\dimen@\wd\mtagbox@
 \advance\dimen@\multlinetaggap@
 \ifdim\dimen@>\displaywidth\rtwidth@\z@\else\rtwidth@\wd\mtagbox@\fi
 \rmultline@@@#1\endmultline$$}
\def\rmultline@@@{\displ@y
 \def\shoveright##1{##1\hfilneg\iftagin@\ifdim\rtwidth@>\z@
  \hskip\rtwidth@\hskip\multlinetaggap@\fi\else\hskip\multlinegap@\fi}%
 \def\shoveleft##1{\setboxz@h{$\m@th\displaystyle{}##1$}%
  \setbox@ne\hbox{$\m@th\displaystyle##1$}%
  \hfilneg\hskip\multlinegap@\hskip.5\wd@ne\hskip-.5\wdz@##1}%
 \halign\bgroup\Let@\hbox to\displaywidth
  {\strut@$\m@th\displaystyle\hfil{}##\hfil$}\crcr
 \hfilneg\hskip\multlinegap@}
\def\endmultline{\iftagsleft@\expandafter\lendmultline@\else
 \expandafter\rendmultline@\fi}
\def\lendmultline@{\hfilneg\hskip\multlinegap@\crcr\egroup}
\def\rendmultline@{\iftagin@                                                
 \ifdim\rtwidth@>\z@                                                        
  \hskip\multlinetaggap@\box\mtagbox@                                       
 \else\llap{\vtop{\normalbaselines\null\hbox{\strut@\box\mtagbox@}}}\fi     
 \else\hskip\multlinegap@\fi                                                
 \hfilneg\crcr\egroup}
\def\bmod{\mskip-\medmuskip\mkern5mu\mathbin{\fam\z@ mod}\penalty900
 \mkern5mu\mskip-\medmuskip}
\def\pmod#1{\allowbreak\ifinner\mkern8mu\else\mkern18mu\fi
 ({\fam\z@ mod}\,\,#1)}
\def\pod#1{\allowbreak\ifinner\mkern8mu\else\mkern18mu\fi(#1)}
\def\mod#1{\allowbreak\ifinner\mkern12mu\else\mkern18mu\fi{\fam\z@ mod}\,\,#1}
\message{continued fractions,}
\newcount\cfraccount@
\def\cfrac{\bgroup\bgroup\advance\cfraccount@\@ne\strut
 \iffalse{\fi\def\\{\over\displaystyle}\iffalse}\fi}
\def\lcfrac{\bgroup\bgroup\advance\cfraccount@\@ne\strut
 \iffalse{\fi\def\\{\hfill\over\displaystyle}\iffalse}\fi}
\def\rcfrac{\bgroup\bgroup\advance\cfraccount@\@ne\strut\hfill
 \iffalse{\fi\def\\{\over\displaystyle}\iffalse}\fi}
\def\gloop@#1\repeat{\gdef\body{#1}\iterate}
\def\endcfrac{\gloop@\ifnum\cfraccount@>\z@\global\advance\cfraccount@\m@ne
 \egroup\hskip-\nulldelimiterspace\egroup\repeat}
\message{compound symbols,}
\def\binrel@#1{\setboxz@h{\thinmuskip0mu
  \medmuskip\m@ne mu\thickmuskip\@ne mu$#1\m@th$}%
 \setbox@ne\hbox{\thinmuskip0mu\medmuskip\m@ne mu\thickmuskip
  \@ne mu${}#1{}\m@th$}%
 \setbox\tw@\hbox{\hskip\wd@ne\hskip-\wdz@}}
\def\overset#1\to#2{\binrel@{#2}\ifdim\wd\tw@<\z@
 \mathbin{\mathop{\kern\z@#2}\limits^{#1}}\else\ifdim\wd\tw@>\z@
 \mathrel{\mathop{\kern\z@#2}\limits^{#1}}\else
 {\mathop{\kern\z@#2}\limits^{#1}}{}\fi\fi}
\def\underset#1\to#2{\binrel@{#2}\ifdim\wd\tw@<\z@
 \mathbin{\mathop{\kern\z@#2}\limits_{#1}}\else\ifdim\wd\tw@>\z@
 \mathrel{\mathop{\kern\z@#2}\limits_{#1}}\else
 {\mathop{\kern\z@#2}\limits_{#1}}{}\fi\fi}
\def\oversetbrace#1\to#2{\overbrace{#2}^{#1}}
\def\undersetbrace#1\to#2{\underbrace{#2}_{#1}}
\def\sideset#1\and#2\to#3{%
 \setbox@ne\hbox{$\dsize{\vphantom{#3}}#1{#3}\m@th$}%
 \setbox\tw@\hbox{$\dsize{#3}#2\m@th$}%
 \hskip\wd@ne\hskip-\wd\tw@\mathop{\hskip\wd\tw@\hskip-\wd@ne
  {\vphantom{#3}}#1{#3}#2}}
\def\rightarrowfill@#1{\setboxz@h{$#1-\m@th$}\ht\z@\z@
  $#1\m@th\copy\z@\mkern-6mu\cleaders
  \hbox{$#1\mkern-2mu\box\z@\mkern-2mu$}\hfill
  \mkern-6mu\mathord\rightarrow$}
\def\leftarrowfill@#1{\setboxz@h{$#1-\m@th$}\ht\z@\z@
  $#1\m@th\mathord\leftarrow\mkern-6mu\cleaders
  \hbox{$#1\mkern-2mu\copy\z@\mkern-2mu$}\hfill
  \mkern-6mu\box\z@$}
\def\leftrightarrowfill@#1{\setboxz@h{$#1-\m@th$}\ht\z@\z@
  $#1\m@th\mathord\leftarrow\mkern-6mu\cleaders
  \hbox{$#1\mkern-2mu\box\z@\mkern-2mu$}\hfill
  \mkern-6mu\mathord\rightarrow$}
\def\overrightarrow{\mathpalette\overrightarrow@}
\def\overrightarrow@#1#2{\vbox{\ialign{##\crcr\rightarrowfill@#1\crcr
 \noalign{\kern-\ex@\nointerlineskip}$\m@th\hfil#1#2\hfil$\crcr}}}

\def\overleftarrow{\mathpalette\overleftarrow@}
\def\overleftarrow@#1#2{\vbox{\ialign{##\crcr\leftarrowfill@#1\crcr
 \noalign{\kern-\ex@\nointerlineskip}$\m@th\hfil#1#2\hfil$\crcr}}}
\def\overleftrightarrow{\mathpalette\overleftrightarrow@}
\def\overleftrightarrow@#1#2{\vbox{\ialign{##\crcr\leftrightarrowfill@#1\crcr
 \noalign{\kern-\ex@\nointerlineskip}$\m@th\hfil#1#2\hfil$\crcr}}}
\def\underrightarrow{\mathpalette\underrightarrow@}
\def\underrightarrow@#1#2{\vtop{\ialign{##\crcr$\m@th\hfil#1#2\hfil$\crcr
 \noalign{\nointerlineskip}\rightarrowfill@#1\crcr}}}

\def\underleftarrow{\mathpalette\underleftarrow@}
\def\underleftarrow@#1#2{\vtop{\ialign{##\crcr$\m@th\hfil#1#2\hfil$\crcr
 \noalign{\nointerlineskip}\leftarrowfill@#1\crcr}}}
\def\underleftrightarrow{\mathpalette\underleftrightarrow@}
\def\underleftrightarrow@#1#2{\vtop{\ialign{##\crcr$\m@th\hfil#1#2\hfil$\crcr
 \noalign{\nointerlineskip}\leftrightarrowfill@#1\crcr}}}
\message{various kinds of dots,}
\let\DOTSI\relax
\let\DOTSB\relax

\newif\ifmath@
{\uccode`7=`\\ \uccode`8=`m \uccode`9=`a \uccode`0=`t \uccode`!=`h
 \uppercase{\gdef\math@#1#2#3#4#5#6\math@{\global\math@false\ifx 7#1\ifx 8#2%
 \ifx 9#3\ifx 0#4\ifx !#5\xdef\meaning@{#6}\global\math@true\fi\fi\fi\fi\fi}}}
\newif\ifmathch@
{\uccode`7=`c \uccode`8=`h \uccode`9=`\"
 \uppercase{\gdef\mathch@#1#2#3#4#5#6\mathch@{\global\mathch@false
  \ifx 7#1\ifx 8#2\ifx 9#5\global\mathch@true\xdef\meaning@{9#6}\fi\fi\fi}}}
\newcount\classnum@
\def\getmathch@#1.#2\getmathch@{\classnum@#1 \divide\classnum@4096
 \ifcase\number\classnum@\or\or\gdef\thedots@{\dotsb@}\or
 \gdef\thedots@{\dotsb@}\fi}
\newif\ifmathbin@
{\uccode`4=`b \uccode`5=`i \uccode`6=`n
 \uppercase{\gdef\mathbin@#1#2#3{\relaxnext@
  \DNii@##1\mathbin@{\ifx\space@\next\global\mathbin@true\fi}%
 \global\mathbin@false\DN@##1\mathbin@{}%
 \ifx 4#1\ifx 5#2\ifx 6#3\DN@{\FN@\nextii@}\fi\fi\fi\next@}}}
\newif\ifmathrel@
{\uccode`4=`r \uccode`5=`e \uccode`6=`l
 \uppercase{\gdef\mathrel@#1#2#3{\relaxnext@
  \DNii@##1\mathrel@{\ifx\space@\next\global\mathrel@true\fi}%
 \global\mathrel@false\DN@##1\mathrel@{}%
 \ifx 4#1\ifx 5#2\ifx 6#3\DN@{\FN@\nextii@}\fi\fi\fi\next@}}}
\newif\ifmacro@
{\uccode`5=`m \uccode`6=`a \uccode`7=`c
 \uppercase{\gdef\macro@#1#2#3#4\macro@{\global\macro@false
  \ifx 5#1\ifx 6#2\ifx 7#3\global\macro@true
  \xdef\meaning@{\macro@@#4\macro@@}\fi\fi\fi}}}
\def\macro@@#1->#2\macro@@{#2}
\newif\ifDOTS@
\newcount\DOTSCASE@
{\uccode`6=`\\ \uccode`7=`D \uccode`8=`O \uccode`9=`T \uccode`0=`S
 \uppercase{\gdef\DOTS@#1#2#3#4#5{\global\DOTS@false\DN@##1\DOTS@{}%
  \ifx 6#1\ifx 7#2\ifx 8#3\ifx 9#4\ifx 0#5\let\next@\DOTS@@\fi\fi\fi\fi\fi
  \next@}}}
{\uccode`3=`B \uccode`4=`I \uccode`5=`X
 \uppercase{\gdef\DOTS@@#1{\relaxnext@
  \DNii@##1\DOTS@{\ifx\space@\next\global\DOTS@true\fi}%
  \DN@{\FN@\nextii@}%
  \ifx 3#1\global\DOTSCASE@\z@\else
  \ifx 4#1\global\DOTSCASE@\@ne\else
  \ifx 5#1\global\DOTSCASE@\tw@\else\DN@##1\DOTS@{}%
  \fi\fi\fi\next@}}}
\newif\ifnot@
{\uccode`5=`\\ \uccode`6=`n \uccode`7=`o \uccode`8=`t
 \uppercase{\gdef\not@#1#2#3#4{\relaxnext@
  \DNii@##1\not@{\ifx\space@\next\global\not@true\fi}%
 \global\not@false\DN@##1\not@{}%
 \ifx 5#1\ifx 6#2\ifx 7#3\ifx 8#4\DN@{\FN@\nextii@}\fi\fi\fi
 \fi\next@}}}
\newif\ifkeybin@
\def\keybin@{\keybin@true
 \ifx\next+\else\ifx\next=\else\ifx\next<\else\ifx\next>\else\ifx\next-\else
 \ifx\next*\else\ifx\next:\else\keybin@false\fi\fi\fi\fi\fi\fi\fi}
\def\dots{\RIfM@\expandafter\mdots@\else\expandafter\tdots@\fi}
\def\tdots@{\unskip\relaxnext@
 \DN@{$\m@th\mathinner{\ldotp\ldotp\ldotp}\,
   \ifx\next,\,$\else\ifx\next.\,$\else\ifx\next;\,$\else\ifx\next:\,$\else
   \ifx\next?\,$\else\ifx\next!\,$\else$ \fi\fi\fi\fi\fi\fi}%
 \ \FN@\next@}
\def\mdots@{\FN@\mdots@@}
\def\mdots@@{\gdef\thedots@{\dotso@}
 \ifx\next\boldkey\gdef\thedots@\boldkey{\boldkeydots@}\else                
 \ifx\next\boldsymbol\gdef\thedots@\boldsymbol{\boldsymboldots@}\else       
 \ifx,\next\gdef\thedots@{\dotsc}
 \else\ifx\not\next\gdef\thedots@{\dotsb@}
 \else\keybin@
 \ifkeybin@\gdef\thedots@{\dotsb@}
 \else\xdef\meaning@{\meaning\next..........}\xdef\meaning@@{\meaning@}
  \expandafter\math@\meaning@\math@
  \ifmath@
   \expandafter\mathch@\meaning@\mathch@
   \ifmathch@\expandafter\getmathch@\meaning@\getmathch@\fi                 
  \else\expandafter\macro@\meaning@@\macro@                                 
  \ifmacro@                                                                
   \expandafter\not@\meaning@\not@\ifnot@\gdef\thedots@{\dotsb@}
  \else\expandafter\DOTS@\meaning@\DOTS@
  \ifDOTS@
   \ifcase\number\DOTSCASE@\gdef\thedots@{\dotsb@}%
    \or\gdef\thedots@{\dotsi}\else\fi                                      
  \else\expandafter\math@\meaning@\math@                                   
  \ifmath@\expandafter\mathbin@\meaning@\mathbin@
  \ifmathbin@\gdef\thedots@{\dotsb@}
  \else\expandafter\mathrel@\meaning@\mathrel@
  \ifmathrel@\gdef\thedots@{\dotsb@}
  \fi\fi\fi\fi\fi\fi\fi\fi\fi\fi\fi\fi
 \thedots@}
\def\plainldots@{\mathinner{\ldotp\ldotp\ldotp}}
\def\plaincdots@{\mathinner{\cdotp\cdotp\cdotp}}
\def\dotsi{\!\plaincdots@}
\let\dotsb@\plaincdots@
\newif\ifextra@
\newif\ifrightdelim@
\def\rightdelim@{\global\rightdelim@true                                    
 \ifx\next)\else                                                            
 \ifx\next]\else
 \ifx\next\rbrack\else
 \ifx\next\}\else
 \ifx\next\rbrace\else
 \ifx\next\rangle\else
 \ifx\next\rceil\else
 \ifx\next\rfloor\else
 \ifx\next\rgroup\else
 \ifx\next\rmoustache\else
 \ifx\next\right\else
 \ifx\next\bigr\else
 \ifx\next\biggr\else
 \ifx\next\Bigr\else                                                        
 \ifx\next\Biggr\else\global\rightdelim@false
 \fi\fi\fi\fi\fi\fi\fi\fi\fi\fi\fi\fi\fi\fi\fi}
\def\extra@{%
 \global\extra@false\rightdelim@\ifrightdelim@\global\extra@true            
 \else\ifx\next$\global\extra@true                                          
 \else\xdef\meaning@{\meaning\next..........}
 \expandafter\macro@\meaning@\macro@\ifmacro@                               
 \expandafter\DOTS@\meaning@\DOTS@
 \ifDOTS@
 \ifnum\DOTSCASE@=\tw@\global\extra@true                                    
 \fi\fi\fi\fi\fi}
\newif\ifbold@
\def\dotso@{\relaxnext@
 \ifbold@
  \let\next\delayed@
  \DNii@{\extra@\plainldots@\ifextra@\,\fi}%
 \else
  \DNii@{\DN@{\extra@\plainldots@\ifextra@\,\fi}\FN@\next@}%
 \fi
 \nextii@}
\def\extrap@#1{%
 \ifx\next,\DN@{#1\,}\else
 \ifx\next;\DN@{#1\,}\else
 \ifx\next.\DN@{#1\,}\else\extra@
 \ifextra@\DN@{#1\,}\else
 \let\next@#1\fi\fi\fi\fi\next@}
\def\ldots{\DN@{\extrap@\plainldots@}%
 \FN@\next@}
\def\cdots{\DN@{\extrap@\plaincdots@}%
 \FN@\next@}

\def\dotsc{\relaxnext@
 \DN@{\ifx\next;\plainldots@\,\else
  \ifx\next.\plainldots@\,\else\extra@\plainldots@
  \ifextra@\,\fi\fi\fi}%
 \FN@\next@}
\def\cdot{\mathchar"2201 }

\message{special superscripts,}
\def\dddot#1{{\mathop{#1}\limits^{\vbox to-1.4\ex@{\kern-\tw@\ex@
 \hbox{\rm...}\vss}}}}
\def\ddddot#1{{\mathop{#1}\limits^{\vbox to-1.4\ex@{\kern-\tw@\ex@
 \hbox{\rm....}\vss}}}}
\def\sphat{^{\mathchoice{}{}%
 {\,\,\botsmash{\hbox{\lower4\ex@\hbox{$\m@th\widehat{\null}$}}}}%
 {\,\botsmash{\hbox{\lower3\ex@\hbox{$\m@th\hat{\null}$}}}}}}

\def\spacute{^{\!\botsmash{\hbox{\lower\@ne ex\hbox{\'{}}}}}}
\def\spgrave{^{\mathchoice{}{}{}{\!}%
 \botsmash{\hbox{\lower\@ne ex\hbox{\`{}}}}}}
\def\spdot{^{\hbox{\raise\ex@\hbox{\rm.}}}}
\def\spddot{^{\hbox{\raise\ex@\hbox{\rm..}}}}
\def\spdddot{^{\hbox{\raise\ex@\hbox{\rm...}}}}
\def\spddddot{^{\hbox{\raise\ex@\hbox{\rm....}}}}
\def\spbreve{^{\!\botsmash{\hbox{\lower4\ex@\hbox{\u{}}}}}}

\message{\string\text,}
\def\textonlyfont@#1#2{\def#1{\RIfM@
 \Err@{Use \string#1\space only in text}\else#2\fi}}
\textonlyfont@\rm\tenrm
\textonlyfont@\it\tenit
\textonlyfont@\sl\tensl
\textonlyfont@\bf\tenbf
\def\oldnos#1{\RIfM@{\mathcode`\,="013B \fam\@ne#1}\else
 \leavevmode\hbox{$\m@th\mathcode`\,="013B \fam\@ne#1$}\fi}
\def\text{\RIfM@\expandafter\text@\else\expandafter\text@@\fi}
\def\text@@#1{\leavevmode\hbox{#1}}
\def\mathhexbox@#1#2#3{\text{$\m@th\mathchar"#1#2#3$}}
\def\dag{{\mathhexbox@279}}
\def\ddag{{\mathhexbox@27A}}
\def\S{{\mathhexbox@278}}
\def\P{{\mathhexbox@27B}}
\newif\iffirstchoice@
\firstchoice@true
\def\text@#1{\mathchoice
 {\hbox{\everymath{\displaystyle}\def\textfonti{\the\textfont\@ne}%
  \def\textfontii{\the\textfont\tw@}\textdef@@ T#1}}
 {\hbox{\firstchoice@false
  \everymath{\textstyle}\def\textfonti{\the\textfont\@ne}%
  \def\textfontii{\the\textfont\tw@}\textdef@@ T#1}}
 {\hbox{\firstchoice@false
  \everymath{\scriptstyle}\def\textfonti{\the\scriptfont\@ne}%
  \def\textfontii{\the\scriptfont\tw@}\textdef@@ S\rm#1}}
 {\hbox{\firstchoice@false
  \everymath{\scriptscriptstyle}\def\textfonti
  {\the\scriptscriptfont\@ne}%
  \def\textfontii{\the\scriptscriptfont\tw@}\textdef@@ s\rm#1}}}
\def\textdef@@#1{\textdef@#1\rm\textdef@#1\bf\textdef@#1\sl\textdef@#1\it}
\def\rmfam{0}
\def\textdef@#1#2{%
 \DN@{\csname\expandafter\eat@\string#2fam\endcsname}%
 \if S#1\edef#2{\the\scriptfont\next@\relax}%
 \else\if s#1\edef#2{\the\scriptscriptfont\next@\relax}%
 \else\edef#2{\the\textfont\next@\relax}\fi\fi}
\scriptfont\itfam\tenit \scriptscriptfont\itfam\tenit
\scriptfont\slfam\tensl \scriptscriptfont\slfam\tensl
\newif\iftopfolded@
\newif\ifbotfolded@
\def\topfoldedtext{\topfolded@true\botfolded@false\foldedtext@}
\def\botfoldedtext{\botfolded@true\topfolded@false\foldedtext@}
\def\foldedtext{\topfolded@false\botfolded@false\foldedtext@}
\Invalid@\foldedwidth
\def\foldedtext@{\relaxnext@
 \DN@{\ifx\next\foldedwidth\let\next@\nextii@\else
  \DN@{\nextii@\foldedwidth{.3\hsize}}\fi\next@}%
 \DNii@\foldedwidth##1##2{\setbox\z@\vbox
  {\normalbaselines\hsize##1\relax
  \tolerance1600 \noindent\ignorespaces##2}\ifbotfolded@\boxz@\else
  \iftopfolded@\vtop{\unvbox\z@}\else\vcenter{\boxz@}\fi\fi}%
 \FN@\next@}
\message{math font commands,}
\def\bold{\RIfM@\expandafter\bold@\else
 \expandafter\nonmatherr@\expandafter\bold\fi}
\def\bold@#1{{\bold@@{#1}}}
\def\bold@@#1{\fam\bffam\relax#1}
\def\slanted{\RIfM@\expandafter\slanted@\else
 \expandafter\nonmatherr@\expandafter\slanted\fi}
\def\slanted@#1{{\slanted@@{#1}}}
\def\slanted@@#1{\fam\slfam\relax#1}
\def\roman{\RIfM@\expandafter\roman@\else
 \expandafter\nonmatherr@\expandafter\roman\fi}
\def\roman@#1{{\roman@@{#1}}}
\def\roman@@#1{\fam\rmfam\relax#1}
\def\italic{\RIfM@\expandafter\italic@\else
 \expandafter\nonmatherr@\expandafter\italic\fi}
\def\italic@#1{{\italic@@{#1}}}
\def\italic@@#1{\fam\itfam\relax#1}
\def\Cal{\RIfM@\expandafter\Cal@\else
 \expandafter\nonmatherr@\expandafter\Cal\fi}
\def\Cal@#1{{\Cal@@{#1}}}
\def\Cal@@#1{\noaccents@\fam\tw@#1}
\mathchardef\Gamma="0000
\mathchardef\Delta="0001
\mathchardef\Theta="0002
\mathchardef\Lambda="0003
\mathchardef\Xi="0004
\mathchardef\Pi="0005
\mathchardef\Sigma="0006
\mathchardef\Upsilon="0007
\mathchardef\Phi="0008
\mathchardef\Psi="0009
\mathchardef\Omega="000A
\mathchardef\varGamma="0100
\mathchardef\varDelta="0101
\mathchardef\varTheta="0102
\mathchardef\varLambda="0103
\mathchardef\varXi="0104
\mathchardef\varPi="0105
\mathchardef\varSigma="0106
\mathchardef\varUpsilon="0107
\mathchardef\varPhi="0108
\mathchardef\varPsi="0109
\mathchardef\varOmega="010A
\let\alloc@@\alloc@
\def\hexnumber@#1{\ifcase#1 0\or 1\or 2\or 3\or 4\or 5\or 6\or 7\or 8\or
 9\or A\or B\or C\or D\or E\or F\fi}
\def\loadmsam{%
 \font@\tenmsa=msam10
 \font@\sevenmsa=msam7
 \font@\fivemsa=msam5
 \alloc@@8\fam\chardef\sixt@@n\msafam
 \textfont\msafam=\tenmsa
 \scriptfont\msafam=\sevenmsa
 \scriptscriptfont\msafam=\fivemsa
 \edef\next{\hexnumber@\msafam}%
 \mathchardef\dabar@"0\next39
 \edef\dashrightarrow{\mathrel{\dabar@\dabar@\mathchar"0\next4B}}%
 \edef\dashleftarrow{\mathrel{\mathchar"0\next4C\dabar@\dabar@}}%
 \let\dasharrow\dashrightarrow
 \edef\ulcorner{\delimiter"4\next70\next70 }%
 \edef\urcorner{\delimiter"5\next71\next71 }%
 \edef\llcorner{\delimiter"4\next78\next78 }%
 \edef\lrcorner{\delimiter"5\next79\next79 }%
 \edef\yen{{\noexpand\mathhexbox@\next55}}%
 \edef\checkmark{{\noexpand\mathhexbox@\next58}}%
 \edef\circledR{{\noexpand\mathhexbox@\next72}}%
 \edef\maltese{{\noexpand\mathhexbox@\next7A}}%
 \global\let\loadmsam\empty}%
\def\loadmsbm{%
 \font@\tenmsb=msbm10 \font@\sevenmsb=msbm7 \font@\fivemsb=msbm5
 \alloc@@8\fam\chardef\sixt@@n\msbfam
 \textfont\msbfam=\tenmsb
 \scriptfont\msbfam=\sevenmsb \scriptscriptfont\msbfam=\fivemsb
 \global\let\loadmsbm\empty
 }
\def\widehat#1{\ifx\undefined\msbfam \DN@{362}%
  \else \setboxz@h{$\m@th#1$}%
    \edef\next@{\ifdim\wdz@>\tw@ em%
        \hexnumber@\msbfam 5B%
      \else 362\fi}\fi
  \mathaccent"0\next@{#1}}
\def\widetilde#1{\ifx\undefined\msbfam \DN@{365}%
  \else \setboxz@h{$\m@th#1$}%
    \edef\next@{\ifdim\wdz@>\tw@ em%
        \hexnumber@\msbfam 5D%
      \else 365\fi}\fi
  \mathaccent"0\next@{#1}}
\message{\string\newsymbol,}
\def\newsymbol#1#2#3#4#5{\define#1{}%
  \count@#2\relax \advance\count@\m@ne 
 \ifcase\count@
   \ifx\undefined\msafam\loadmsam\fi \let\next@\msafam
 \or \ifx\undefined\msbfam\loadmsbm\fi \let\next@\msbfam
 \else  \Err@{\Invalid@@\string\newsymbol}\let\next@\tw@\fi
 \mathchardef#1="#3\hexnumber@\next@#4#5\space}

\def\Bbb{\RIfM@\expandafter\Bbb@\else
 \expandafter\nonmatherr@\expandafter\Bbb\fi}
\def\Bbb@#1{{\Bbb@@{#1}}}
\def\Bbb@@#1{\noaccents@\fam\msbfam\relax#1}
\message{bold Greek and bold symbols,}
\def\loadbold{%
 \font@\tencmmib=cmmib10 \font@\sevencmmib=cmmib7 \font@\fivecmmib=cmmib5
 \skewchar\tencmmib'177 \skewchar\sevencmmib'177 \skewchar\fivecmmib'177
 \alloc@@8\fam\chardef\sixt@@n\cmmibfam
 \textfont\cmmibfam\tencmmib
 \scriptfont\cmmibfam\sevencmmib \scriptscriptfont\cmmibfam\fivecmmib
 \font@\tencmbsy=cmbsy10 \font@\sevencmbsy=cmbsy7 \font@\fivecmbsy=cmbsy5
 \skewchar\tencmbsy'60 \skewchar\sevencmbsy'60 \skewchar\fivecmbsy'60
 \alloc@@8\fam\chardef\sixt@@n\cmbsyfam
 \textfont\cmbsyfam\tencmbsy
 \scriptfont\cmbsyfam\sevencmbsy \scriptscriptfont\cmbsyfam\fivecmbsy
 \let\loadbold\empty
}
\def\boldnotloaded#1{\Err@{\ifcase#1\or First\else Second\fi
       bold symbol font not loaded}}
\def\mathchari@#1#2#3{\ifx\undefined\cmmibfam
    \boldnotloaded@\@ne
  \else\mathchar"#1\hexnumber@\cmmibfam#2#3\space \fi}
\def\mathcharii@#1#2#3{\ifx\undefined\cmbsyfam
    \boldnotloaded\tw@
  \else \mathchar"#1\hexnumber@\cmbsyfam#2#3\space\fi}
\edef\bffam@{\hexnumber@\bffam}
\def\boldkey#1{\ifcat\noexpand#1A%
  \ifx\undefined\cmmibfam \boldnotloaded\@ne
  \else {\fam\cmmibfam#1}\fi
 \else
 \ifx#1!\mathchar"5\bffam@21 \else
 \ifx#1(\mathchar"4\bffam@28 \else\ifx#1)\mathchar"5\bffam@29 \else
 \ifx#1+\mathchar"2\bffam@2B \else\ifx#1:\mathchar"3\bffam@3A \else
 \ifx#1;\mathchar"6\bffam@3B \else\ifx#1=\mathchar"3\bffam@3D \else
 \ifx#1?\mathchar"5\bffam@3F \else\ifx#1[\mathchar"4\bffam@5B \else
 \ifx#1]\mathchar"5\bffam@5D \else
 \ifx#1,\mathchari@63B \else
 \ifx#1-\mathcharii@200 \else
 \ifx#1.\mathchari@03A \else
 \ifx#1/\mathchari@03D \else
 \ifx#1<\mathchari@33C \else
 \ifx#1>\mathchari@33E \else
 \ifx#1*\mathcharii@203 \else
 \ifx#1|\mathcharii@06A \else
 \ifx#10\bold0\else\ifx#11\bold1\else\ifx#12\bold2\else\ifx#13\bold3\else
 \ifx#14\bold4\else\ifx#15\bold5\else\ifx#16\bold6\else\ifx#17\bold7\else
 \ifx#18\bold8\else\ifx#19\bold9\else
  \Err@{\string\boldkey\space can't be used with #1}%
 \fi\fi\fi\fi\fi\fi\fi\fi\fi\fi\fi\fi\fi\fi\fi
 \fi\fi\fi\fi\fi\fi\fi\fi\fi\fi\fi\fi\fi\fi}
\def\boldsymbol#1{%
 \DN@{\Err@{You can't use \string\boldsymbol\space with \string#1}#1}%
 \ifcat\noexpand#1A%
   \let\next@\relax
   \ifx\undefined\cmmibfam \boldnotloaded\@ne
   \else {\fam\cmmibfam#1}\fi
 \else
  \xdef\meaning@{\meaning#1.........}%
  \expandafter\math@\meaning@\math@
  \ifmath@
   \expandafter\mathch@\meaning@\mathch@
   \ifmathch@
    \expandafter\boldsymbol@@\meaning@\boldsymbol@@
   \fi
  \else
   \expandafter\macro@\meaning@\macro@
   \expandafter\delim@\meaning@\delim@
   \ifdelim@
    \expandafter\delim@@\meaning@\delim@@
   \else
    \boldsymbol@{#1}%
   \fi
  \fi
 \fi
 \next@}
\def\mathhexboxii@#1#2{\ifx\undefined\cmbsyfam
    \boldnotloaded\tw@
  \else \mathhexbox@{\hexnumber@\cmbsyfam}{#1}{#2}\fi}
\def\boldsymbol@#1{\let\next@\relax\let\next#1%
 \ifx\next\cdot\mathcharii@201 \else
 \ifx\next\prime{{\null\mathcharii@030 \null}}\else
 \ifx\next\lbrack\mathchar"4\bffam@5B \else
 \ifx\next\rbrack\mathchar"5\bffam@5D \else
 \ifx\next\{\mathcharii@466 \else
 \ifx\next\lbrace\mathcharii@466 \else
 \ifx\next\}\mathcharii@567 \else
 \ifx\next\rbrace\mathcharii@567 \else
 \ifx\next\surd{{\mathcharii@170}}\else
 \ifx\next\S{{\mathhexboxii@78}}\else
 \ifx\next\P{{\mathhexboxii@7B}}\else
 \ifx\next\dag{{\mathhexboxii@79}}\else
 \ifx\next\ddag{{\mathhexboxii@7A}}\else
 \DN@{\Err@{You can't use \string\boldsymbol\space with \string#1}#1}%
 \fi\fi\fi\fi\fi\fi\fi\fi\fi\fi\fi\fi\fi}
\def\boldsymbol@@#1.#2\boldsymbol@@{\classnum@#1 \count@@@\classnum@        
 \divide\classnum@4096 \count@\classnum@                                    
 \multiply\count@4096 \advance\count@@@-\count@ \count@@\count@@@           
 \divide\count@@@\@cclvi \count@\count@@                                    
 \multiply\count@@@\@cclvi \advance\count@@-\count@@@                       
 \divide\count@@@\@cclvi                                                    
 \multiply\classnum@4096 \advance\classnum@\count@@                         
 \ifnum\count@@@=\z@                                                        
  \count@"\bffam@ \multiply\count@\@cclvi
  \advance\classnum@\count@
  \DN@{\mathchar\number\classnum@}%
 \else
  \ifnum\count@@@=\@ne                                                      
   \ifx\undefined\cmmibfam \DN@{\boldnotloaded\@ne}%
   \else \count@\cmmibfam \multiply\count@\@cclvi
     \advance\classnum@\count@
     \DN@{\mathchar\number\classnum@}\fi
  \else
   \ifnum\count@@@=\tw@                                                    
     \ifx\undefined\cmbsyfam
       \DN@{\boldnotloaded\tw@}%
     \else
       \count@\cmbsyfam \multiply\count@\@cclvi
       \advance\classnum@\count@
       \DN@{\mathchar\number\classnum@}%
     \fi
  \fi
 \fi
\fi}
\newif\ifdelim@
\newcount\delimcount@
{\uccode`6=`\\ \uccode`7=`d \uccode`8=`e \uccode`9=`l
 \uppercase{\gdef\delim@#1#2#3#4#5\delim@
  {\delim@false\ifx 6#1\ifx 7#2\ifx 8#3\ifx 9#4\delim@true
   \xdef\meaning@{#5}\fi\fi\fi\fi}}}
\def\delim@@#1"#2#3#4#5#6\delim@@{\if#32%
\let\next@\relax
 \ifx\undefined\cmbsyfam \boldnotloaded\@ne
 \else \mathcharii@#2#4#5\space \fi\fi}
\def\vert{\delimiter"026A30C }
\def\Vert{\delimiter"026B30D }
\let\|\Vert
\def\backslash{\delimiter"026E30F }
\def\boldkeydots@#1{\bold@true\let\next=#1\let\delayed@=#1\mdots@@
 \boldkey#1\bold@false}  
\def\boldsymboldots@#1{\bold@true\let\next#1\let\delayed@#1\mdots@@
 \boldsymbol#1\bold@false}
\message{Euler fonts,}

\def\frak{\mathfont@\frak}

\def\loadmathfont#1{%
   \expandafter\font@\csname ten#1\endcsname=#110
   \expandafter\font@\csname seven#1\endcsname=#17
   \expandafter\font@\csname five#1\endcsname=#15
   \edef\next{\noexpand\alloc@@8\fam\chardef\sixt@@n
     \expandafter\noexpand\csname#1fam\endcsname}%
   \next
   \textfont\csname#1fam\endcsname \csname ten#1\endcsname
   \scriptfont\csname#1fam\endcsname \csname seven#1\endcsname
   \scriptscriptfont\csname#1fam\endcsname \csname five#1\endcsname
   \expandafter\def\csname #1\expandafter\endcsname\expandafter{%
      \expandafter\mathfont@\csname#1\endcsname}%
 \expandafter\gdef\csname load#1\endcsname{}%
}
\def\mathfont@#1{\RIfM@\expandafter\mathfont@@\expandafter#1\else
  \expandafter\nonmatherr@\expandafter#1\fi}
\def\mathfont@@#1#2{{\mathfont@@@#1{#2}}}
\def\mathfont@@@#1#2{\noaccents@
   \fam\csname\expandafter\eat@\string#1fam\endcsname
   \relax#2}
\message{math accents,}
\def\accentclass@{7}
\def\noaccents@{\def\accentclass@{0}}
\def\makeacc@#1#2{\def#1{\mathaccent"\accentclass@#2 }}
\makeacc@\hat{05E}
\makeacc@\check{014}
\makeacc@\tilde{07E}
\makeacc@\acute{013}
\makeacc@\grave{012}
\makeacc@\dot{05F}
\makeacc@\ddot{07F}
\makeacc@\breve{015}
\makeacc@\bar{016}

\newcount\skewcharcount@
\newcount\familycount@
\def\theskewchar@{\familycount@\@ne
 \global\skewcharcount@\the\skewchar\textfont\@ne                           
 \ifnum\fam>\m@ne\ifnum\fam<16
  \global\familycount@\the\fam\relax
  \global\skewcharcount@\the\skewchar\textfont\the\fam\relax\fi\fi          
 \ifnum\skewcharcount@>\m@ne
  \ifnum\skewcharcount@<128
  \multiply\familycount@256
  \global\advance\skewcharcount@\familycount@
  \global\advance\skewcharcount@28672
  \mathchar\skewcharcount@\else
  \global\skewcharcount@\m@ne\fi\else
 \global\skewcharcount@\m@ne\fi}                                            
\newcount\pointcount@
\def\getpoints@#1.#2\getpoints@{\pointcount@#1 }
\newdimen\accentdimen@
\newcount\accentmu@
\def\dimentomu@{\multiply\accentdimen@ 100
 \expandafter\getpoints@\the\accentdimen@\getpoints@
 \multiply\pointcount@18
 \divide\pointcount@\@m
 \global\accentmu@\pointcount@}
\def\Makeacc@#1#2{\def#1{\RIfM@\DN@{\mathaccent@
 {"\accentclass@#2 }}\else\DN@{\nonmatherr@{#1}}\fi\next@}}
\def\unbracefonts@{\let\Cal@\Cal@@\let\roman@\roman@@\let\bold@\bold@@
 \let\slanted@\slanted@@}
\def\mathaccent@#1#2{\ifnum\fam=\m@ne\xdef\thefam@{1}\else
 \xdef\thefam@{\the\fam}\fi                                                 
 \accentdimen@\z@                                                           
 \setboxz@h{\unbracefonts@$\m@th\fam\thefam@\relax#2$}
 \ifdim\accentdimen@=\z@\DN@{\mathaccent#1{#2}}
  \setbox@ne\hbox{\unbracefonts@$\m@th\fam\thefam@\relax#2\theskewchar@$}
  \setbox\tw@\hbox{$\m@th\ifnum\skewcharcount@=\m@ne\else
   \mathchar\skewcharcount@\fi$}
  \global\accentdimen@\wd@ne\global\advance\accentdimen@-\wdz@
  \global\advance\accentdimen@-\wd\tw@                                     
  \global\multiply\accentdimen@\tw@
  \dimentomu@\global\advance\accentmu@\@ne                                 
 \else\DN@{{\mathaccent#1{#2\mkern\accentmu@ mu}%
    \mkern-\accentmu@ mu}{}}\fi                                             
 \next@}\Makeacc@\Hat{05E}
\Makeacc@\Check{014}
\Makeacc@\Tilde{07E}
\Makeacc@\Acute{013}
\Makeacc@\Grave{012}
\Makeacc@\Dot{05F}
\Makeacc@\Ddot{07F}
\Makeacc@\Breve{015}
\Makeacc@\Bar{016}
\def\Vec{\RIfM@\DN@{\mathaccent@{"017E }}\else
 \DN@{\nonmatherr@\Vec}\fi\next@}
\def\accentedsymbol#1#2{\csname newbox\expandafter\endcsname
  \csname\expandafter\eat@\string#1@box\endcsname
 \expandafter\setbox\csname\expandafter\eat@
  \string#1@box\endcsname\hbox{$\m@th#2$}\define
  #1{\copy\csname\expandafter\eat@\string#1@box\endcsname{}}}
\message{roots,}
\def\sqrt#1{\radical"270370 {#1}}
\let\underline@\underline
\let\overline@\overline
\def\underline#1{\underline@{#1}}
\def\overline#1{\overline@{#1}}
\Invalid@\leftroot
\Invalid@\uproot
\newcount\uproot@
\newcount\leftroot@
\def\root{\relaxnext@
  \DN@{\ifx\next\uproot\let\next@\nextii@\else
   \ifx\next\leftroot\let\next@\nextiii@\else
   \let\next@\plainroot@\fi\fi\next@}%
  \DNii@\uproot##1{\uproot@##1\relax\FN@\nextiv@}%
  \def\nextiv@{\ifx\next\space@\DN@. {\FN@\nextv@}\else
   \DN@.{\FN@\nextv@}\fi\next@.}%
  \def\nextv@{\ifx\next\leftroot\let\next@\nextvi@\else
   \let\next@\plainroot@\fi\next@}%
  \def\nextvi@\leftroot##1{\leftroot@##1\relax\plainroot@}%
   \def\nextiii@\leftroot##1{\leftroot@##1\relax\FN@\nextvii@}%
  \def\nextvii@{\ifx\next\space@
   \DN@. {\FN@\nextviii@}\else
   \DN@.{\FN@\nextviii@}\fi\next@.}%
  \def\nextviii@{\ifx\next\uproot\let\next@\nextix@\else
   \let\next@\plainroot@\fi\next@}%
  \def\nextix@\uproot##1{\uproot@##1\relax\plainroot@}%
  \bgroup\uproot@\z@\leftroot@\z@\FN@\next@}
\def\plainroot@#1\of#2{\setbox\rootbox\hbox{$\m@th\scriptscriptstyle{#1}$}%
 \mathchoice{\r@@t\displaystyle{#2}}{\r@@t\textstyle{#2}}
 {\r@@t\scriptstyle{#2}}{\r@@t\scriptscriptstyle{#2}}\egroup}
\def\r@@t#1#2{\setboxz@h{$\m@th#1\sqrt{#2}$}%
 \dimen@\ht\z@\advance\dimen@-\dp\z@
 \setbox@ne\hbox{$\m@th#1\mskip\uproot@ mu$}\advance\dimen@ 1.667\wd@ne
 \mkern-\leftroot@ mu\mkern5mu\raise.6\dimen@\copy\rootbox
 \mkern-10mu\mkern\leftroot@ mu\boxz@}
\def\boxed#1{\setboxz@h{$\m@th\displaystyle{#1}$}\dimen@.4\ex@
 \advance\dimen@3\ex@\advance\dimen@\dp\z@
 \hbox{\lower\dimen@\hbox{%
 \vbox{\hrule height.4\ex@
 \hbox{\vrule width.4\ex@\hskip3\ex@\vbox{\vskip3\ex@\boxz@\vskip3\ex@}%
 \hskip3\ex@\vrule width.4\ex@}\hrule height.4\ex@}%
 }}}
\message{commutative diagrams,}
\let\ampersand@\relax
\newdimen\minaw@
\minaw@11.11128\ex@
\newdimen\minCDaw@
\minCDaw@2.5pc
\def\minCDarrowwidth#1{\RIfMIfI@\onlydmatherr@\minCDarrowwidth
 \else\minCDaw@#1\relax\fi\else\onlydmatherr@\minCDarrowwidth\fi}
\newif\ifCD@
\def\CD{\bgroup\vspace@\relax\let\ampersand@&\iffalse}\fi
 \CD@true\vcenter\bgroup\Let@\tabskip\z@skip\baselineskip20\ex@
 \lineskip3\ex@\lineskiplimit3\ex@\halign\bgroup
 &\hfill$\m@th##$\hfill\crcr}
\def\endCD{\crcr\egroup\egroup\egroup}
\newdimen\bigaw@
\atdef@>#1>#2>{\ampersand@                                                  
 \setboxz@h{$\m@th\ssize\;{#1}\;\;$}
 \setbox@ne\hbox{$\m@th\ssize\;{#2}\;\;$}
 \setbox\tw@\hbox{$\m@th#2$}
 \ifCD@\global\bigaw@\minCDaw@\else\global\bigaw@\minaw@\fi                 
 \ifdim\wdz@>\bigaw@\global\bigaw@\wdz@\fi
 \ifdim\wd@ne>\bigaw@\global\bigaw@\wd@ne\fi                                
 \ifCD@\enskip\fi                                                           
 \ifdim\wd\tw@>\z@
  \mathrel{\mathop{\hbox to\bigaw@{\rightarrowfill@\displaystyle}}%
    \limits^{#1}_{#2}}
 \else\mathrel{\mathop{\hbox to\bigaw@{\rightarrowfill@\displaystyle}}%
    \limits^{#1}}\fi                                                        
 \ifCD@\enskip\fi                                                          
 \ampersand@}                                                              
\atdef@<#1<#2<{\ampersand@\setboxz@h{$\m@th\ssize\;\;{#1}\;$}%
 \setbox@ne\hbox{$\m@th\ssize\;\;{#2}\;$}\setbox\tw@\hbox{$\m@th#2$}%
 \ifCD@\global\bigaw@\minCDaw@\else\global\bigaw@\minaw@\fi
 \ifdim\wdz@>\bigaw@\global\bigaw@\wdz@\fi
 \ifdim\wd@ne>\bigaw@\global\bigaw@\wd@ne\fi
 \ifCD@\enskip\fi
 \ifdim\wd\tw@>\z@
  \mathrel{\mathop{\hbox to\bigaw@{\leftarrowfill@\displaystyle}}%
       \limits^{#1}_{#2}}\else
  \mathrel{\mathop{\hbox to\bigaw@{\leftarrowfill@\displaystyle}}%
       \limits^{#1}}\fi
 \ifCD@\enskip\fi\ampersand@}
\begingroup
 \catcode`\~=\active \lccode`\~=`\@
 \lowercase{%
  \global\atdef@)#1)#2){~>#1>#2>}
  \global\atdef@(#1(#2({~<#1<#2<}}
\endgroup
\atdef@ A#1A#2A{\llap{$\m@th\vcenter{\hbox
 {$\ssize#1$}}$}\Big\uparrow\rlap{$\m@th\vcenter{\hbox{$\ssize#2$}}$}&&}
\atdef@ V#1V#2V{\llap{$\m@th\vcenter{\hbox
 {$\ssize#1$}}$}\Big\downarrow\rlap{$\m@th\vcenter{\hbox{$\ssize#2$}}$}&&}
\atdef@={&\enskip\mathrel
 {\vbox{\hrule width\minCDaw@\vskip3\ex@\hrule width
 \minCDaw@}}\enskip&}
\atdef@|{\Big\Vert&&}
\atdef@\vert{\Big\Vert&&}
\def\pretend#1\haswidth#2{\setboxz@h{$\m@th\scriptstyle{#2}$}\hbox
 to\wdz@{\hfill$\m@th\scriptstyle{#1}$\hfill}}
\message{poor man's bold,}
\def\pmb{\RIfM@\expandafter\mathpalette\expandafter\pmb@\else
 \expandafter\pmb@@\fi}
\def\pmb@@#1{\leavevmode\setboxz@h{#1}%
   \dimen@-\wdz@
   \kern-.5\ex@\copy\z@
   \kern\dimen@\kern.25\ex@\raise.4\ex@\copy\z@
   \kern\dimen@\kern.25\ex@\box\z@
}
\def\binrel@@#1{\ifdim\wd2<\z@\mathbin{#1}\else\ifdim\wd\tw@>\z@
 \mathrel{#1}\else{#1}\fi\fi}
\newdimen\pmbraise@
\def\pmb@#1#2{\setbox\thr@@\hbox{$\m@th#1{#2}$}%
 \setbox4\hbox{$\m@th#1\mkern.5mu$}\pmbraise@\wd4\relax
 \binrel@{#2}%
 \dimen@-\wd\thr@@
   \binrel@@{%
   \mkern-.8mu\copy\thr@@
   \kern\dimen@\mkern.4mu\raise\pmbraise@\copy\thr@@
   \kern\dimen@\mkern.4mu\box\thr@@
}}
\def\documentstyle#1{\W@{}\input #1.sty\relax}
\message{syntax check,}
\font\dummyft@=dummy
\fontdimen1 \dummyft@=\z@
\fontdimen2 \dummyft@=\z@
\fontdimen3 \dummyft@=\z@
\fontdimen4 \dummyft@=\z@
\fontdimen5 \dummyft@=\z@
\fontdimen6 \dummyft@=\z@
\fontdimen7 \dummyft@=\z@
\fontdimen8 \dummyft@=\z@
\fontdimen9 \dummyft@=\z@
\fontdimen10 \dummyft@=\z@
\fontdimen11 \dummyft@=\z@
\fontdimen12 \dummyft@=\z@
\fontdimen13 \dummyft@=\z@
\fontdimen14 \dummyft@=\z@
\fontdimen15 \dummyft@=\z@
\fontdimen16 \dummyft@=\z@
\fontdimen17 \dummyft@=\z@
\fontdimen18 \dummyft@=\z@
\fontdimen19 \dummyft@=\z@
\fontdimen20 \dummyft@=\z@
\fontdimen21 \dummyft@=\z@
\fontdimen22 \dummyft@=\z@
\def\fontlist@{\\{\tenrm}\\{\sevenrm}\\{\fiverm}\\{\teni}\\{\seveni}%
 \\{\fivei}\\{\tensy}\\{\sevensy}\\{\fivesy}\\{\tenex}\\{\tenbf}\\{\sevenbf}%
 \\{\fivebf}\\{\tensl}\\{\tenit}}
\def\font@#1=#2 {\rightappend@#1\to\fontlist@\font#1=#2 }
\def\dodummy@{{\def\\##1{\global\let##1\dummyft@}\fontlist@}}
\def\nopages@{\output{\setbox\z@\box\@cclv \deadcycles\z@}%
 \alloc@5\toks\toksdef\@cclvi\output}
\let\galleys\nopages@
\newif\ifsyntax@
\newcount\countxviii@
\def\syntax{\syntax@true\dodummy@\countxviii@\count18
 \loop\ifnum\countxviii@>\m@ne\textfont\countxviii@=\dummyft@
 \scriptfont\countxviii@=\dummyft@\scriptscriptfont\countxviii@=\dummyft@
 \advance\countxviii@\m@ne\repeat                                           
 \dummyft@\tracinglostchars\z@\nopages@\frenchspacing\hbadness\@M}
\def\first@#1#2\end{#1}
\def\printoptions{\W@{Do you want S(yntax check),
  G(alleys) or P(ages)?}%
 \message{Type S, G or P, followed by <return>: }%
 \begingroup 
 \endlinechar\m@ne 
 \read\m@ne to\ans@
 \edef\ans@{\uppercase{\def\noexpand\ans@{%
   \expandafter\first@\ans@ P\end}}}%
 \expandafter\endgroup\ans@
 \if\ans@ P
 \else \if\ans@ S\syntax
 \else \if\ans@ G\galleys
 \else\message{? Unknown option: \ans@; using the `pages' option.}%
 \fi\fi\fi}
\def\alloc@#1#2#3#4#5{\global\advance\count1#1by\@ne
 \ch@ck#1#4#2\allocationnumber=\count1#1
 \global#3#5=\allocationnumber
 \ifalloc@\wlog{\string#5=\string#2\the\allocationnumber}\fi}
\def\document{\def\alloclist@{}\def\fontlist@{}}
\let\enddocument\bye

\let\proclaim\undefined
\let\footnote\undefined
\let\=\undefined
\let\>\undefined

\catcode`\@=\active
\message{... finished}

\expandafter\ifx\csname mathdefs.tex\endcsname\relax
  \expandafter\gdef\csname mathdefs.tex\endcsname{}
\else \message{Hey!  Apparently you were trying to
  \string\input{mathdefs.tex} twice.   This does not make sense.} 
\errmessage{Please edit your file (probably \jobname.tex) and remove
any duplicate ``\string\input'' lines}\endinput\fi




\catcode`\X=12\catcode`\@=11

\def\n@wcount{\alloc@0\count\countdef\insc@unt}
\def\n@wwrite{\alloc@7\write\chardef\sixt@@n}
\def\n@wread{\alloc@6\read\chardef\sixt@@n}
\def\r@s@t{\relax}\def\v@idline{\par}\def\@mputate#1/{#1}
\def\l@c@l#1X{\firstpart.#1}\def\gl@b@l#1X{#1}\def\t@d@l#1X{{}}

\def\crossrefs#1{\ifx\all#1\let\tr@ce=\all\else\def\tr@ce{#1,}\fi
   \n@wwrite\cit@tionsout\openout\cit@tionsout=\jobname.cit 
   \write\cit@tionsout{\tr@ce}\expandafter\setfl@gs\tr@ce,}
\def\setfl@gs#1,{\def\@{#1}\ifx\@\empty\let\next=\relax
   \else\let\next=\setfl@gs\expandafter\xdef
   \csname#1tr@cetrue\endcsname{}\fi\next}
\def\m@ketag#1#2{\expandafter\n@wcount\csname#2tagno\endcsname
     \csname#2tagno\endcsname=0\let\tail=\all\xdef\all{\tail#2,}
   \ifx#1\l@c@l\let\tail=\r@s@t\xdef\r@s@t{\csname#2tagno\endcsname=0\tail}\fi
   \expandafter\gdef\csname#2cite\endcsname##1{\expandafter
     \ifx\csname#2tag##1\endcsname\relax?\else\csname#2tag##1\endcsname\fi
     \expandafter\ifx\csname#2tr@cetrue\endcsname\relax\else
     \write\cit@tionsout{#2tag ##1 cited on page \folio.}\fi}
   \expandafter\gdef\csname#2page\endcsname##1{\expandafter
     \ifx\csname#2page##1\endcsname\relax?\else\csname#2page##1\endcsname\fi
     \expandafter\ifx\csname#2tr@cetrue\endcsname\relax\else
     \write\cit@tionsout{#2tag ##1 cited on page \folio.}\fi}
   \expandafter\gdef\csname#2tag\endcsname##1{\expandafter
      \ifx\csname#2check##1\endcsname\relax
      \expandafter\xdef\csname#2check##1\endcsname{}%
      \else\immediate\write16{Warning: #2tag ##1 used more than once.}\fi
      \multit@g{#1}{#2}##1/X%
      \write\t@gsout{#2tag ##1 assigned number \csname#2tag##1\endcsname\space
      on page \number\count0.}%
   \csname#2tag##1\endcsname}}

\def\multit@g#1#2#3/#4X{\def\t@mp{#4}\ifx\t@mp\empty%
      \global\advance\csname#2tagno\endcsname by 1 
      \expandafter\xdef\csname#2tag#3\endcsname
      {#1\number\csname#2tagno\endcsnameX}%
   \else\expandafter\ifx\csname#2last#3\endcsname\relax
      \expandafter\n@wcount\csname#2last#3\endcsname
      \global\advance\csname#2tagno\endcsname by 1 
      \expandafter\xdef\csname#2tag#3\endcsname
      {#1\number\csname#2tagno\endcsnameX}
      \write\t@gsout{#2tag #3 assigned number \csname#2tag#3\endcsname\space
      on page \number\count0.}\fi
   \global\advance\csname#2last#3\endcsname by 1
   \def\t@mp{\expandafter\xdef\csname#2tag#3/}%
   \expandafter\t@mp\@mputate#4\endcsname
   {\csname#2tag#3\endcsname\lastpart{\csname#2last#3\endcsname}}\fi}
\def\t@gs#1{\def\all{}\m@ketag#1e\m@ketag#1s\m@ketag\t@d@l p
\let\realscite\scite
\let\realstag\stag
   \m@ketag\gl@b@l r \n@wread\t@gsin
   \openin\t@gsin=\jobname.tgs \re@der \closein\t@gsin
   \n@wwrite\t@gsout\openout\t@gsout=\jobname.tgs }
\outer\def\localtags{\t@gs\l@c@l}
\outer\def\globaltags{\t@gs\gl@b@l}
\outer\def\newlocaltag#1{\m@ketag\l@c@l{#1}}
\outer\def\newglobaltag#1{\m@ketag\gl@b@l{#1}}

\newif\ifpr@ 
\def\m@kecs #1tag #2 assigned number #3 on page #4.%
   {\expandafter\gdef\csname#1tag#2\endcsname{#3}
   \expandafter\gdef\csname#1page#2\endcsname{#4}
   \ifpr@\expandafter\xdef\csname#1check#2\endcsname{}\fi}
\def\re@der{\ifeof\t@gsin\let\next=\relax\else
   \read\t@gsin to\t@gline\ifx\t@gline\v@idline\else
   \expandafter\m@kecs \t@gline\fi\let \next=\re@der\fi\next}
\def\pretags#1{\pr@true\pret@gs#1,,}
\def\pret@gs#1,{\def\@{#1}\ifx\@\empty\let\n@xtfile=\relax
   \else\let\n@xtfile=\pret@gs \openin\t@gsin=#1.tgs \message{#1} \re@der 
   \closein\t@gsin\fi \n@xtfile}

\newcount\sectno\sectno=0\newcount\subsectno\subsectno=0
\newif\ifultr@local \def\ultralocal{\ultr@localtrue}
\def\firstpart{\number\sectno}
\def\lastpart#1{\ifcase#1 \or a\or b\or c\or d\or e\or f\or g\or h\or 
   i\or k\or l\or m\or n\or o\or p\or q\or r\or s\or t\or u\or v\or w\or 
   x\or y\or z \fi}

\def\resetall{\global\advance\sectno by 1\subsectno=0
   \gdef\firstpart{\number\sectno}\r@s@t}
\def\resetsub{\global\advance\subsectno by 1
   \gdef\firstpart{\number\sectno.\number\subsectno}\r@s@t}
\def\newsection#1\par{\resetall\vskip0pt plus.3\vsize\penalty-250
   \vskip0pt plus-.3\vsize\bigskip\bigskip
   \message{#1}\leftline{\bf#1}\nobreak\bigskip}
\def\subsection#1\par{\ifultr@local\resetsub\fi
   \vskip0pt plus.2\vsize\penalty-250\vskip0pt plus-.2\vsize
   \bigskip\smallskip\message{#1}\leftline{\bf#1}\nobreak\medskip}


\newdimen\marginshift

\newdimen\margindelta
\newdimen\marginmax
\newdimen\marginmin

\def\margininit{       
\marginmax=3 true cm                  
				      
\margindelta=0.1 true cm              
\marginmin=0.1true cm                 
\marginshift=\marginmin
}    

\def\t@gsjj#1,{\def\@{#1}\ifx\@\empty\let\next=\relax\else\let\next=\t@gsjj
   \def\@@{p}\ifx\@\@@\else
   \expandafter\gdef\csname#1cite\endcsname##1{\citejj{##1}}
   \expandafter\gdef\csname#1page\endcsname##1{?}
   \expandafter\gdef\csname#1tag\endcsname##1{\tagjj{##1}}\fi\fi\next}
\newif\ifshowstuffinmargin
\showstuffinmarginfalse
\def\jjtags{\ifx\shlhetal\relax 
  \else
\ifx\shlhetal\undefinedcontrolseq
\else
\showstuffinmargintrue
\ifx\all\relax\else\expandafter\t@gsjj\all,\fi\fi \fi
}

\def\tagjj#1{\realstag{#1}\mginpar{\zeigen{#1}}}
\def\citejj#1{\rechnen{#1}\mginpar{\zeigen{#1}}}     

\def\rechnen#1{\expandafter\ifx\csname stag#1\endcsname\relax ??\else
                           \csname stag#1\endcsname\fi}

\newdimen\theight

\def\marginfont{\sevenrm}

\def\trymarginbox#1{\setbox0=\hbox{\marginfont\hskip\marginshift #1}%
		\global\marginshift\wd0 
		\global\advance\marginshift\margindelta}

\def \mginpar#1{%
\ifvmode\setbox0\hbox to \hsize{\hfill\rlap{\marginfont\quad#1}}%
\ht0 0cm
\dp0 0cm
\box0\vskip-\baselineskip
\else 
             \vadjust{\trymarginbox{#1}%
		\ifdim\marginshift>\marginmax \global\marginshift\marginmin
			\trymarginbox{#1}%
                \fi
             \theight=\ht0
             \advance\theight by \dp0    \advance\theight by \lineskip
             \kern -\theight \vbox to \theight{\rightline{\rlap{\box0}}%
\vss}}\fi}


\def\t@gsoff#1,{\def\@{#1}\ifx\@\empty\let\next=\relax\else\let\next=\t@gsoff
   \def\@@{p}\ifx\@\@@\else
   \expandafter\gdef\csname#1cite\endcsname##1{\zeigen{##1}}
   \expandafter\gdef\csname#1page\endcsname##1{?}
   \expandafter\gdef\csname#1tag\endcsname##1{\zeigen{##1}}\fi\fi\next}
\def\verbatimtags{\showstuffinmarginfalse
\ifx\all\relax\else\expandafter\t@gsoff\all,\fi}
\def\zeigen#1{\hbox{$\langle$}#1\hbox{$\rangle$}}

\def\margincite#1{\ifshowstuffinmargin\mginpar{\zeigen{#1}}\fi}

\def\margintag#1{\ifshowstuffinmargin\mginpar{\zeigen{#1}}\fi}

\def\(#1){\edef\dot@g{\ifmmode\ifinner(\hbox{\noexpand\etag{#1}})
   \else\noexpand\eqno(\hbox{\noexpand\etag{#1}})\fi
   \else(\noexpand\ecite{#1})\fi}\dot@g}

\newif\ifbr@ck
\def\eat#1{}
\def\[#1]{\br@cktrue[\br@cket#1'X]}
\def\br@cket#1'#2X{\def\temp{#2}\ifx\temp\empty\let\next\eat
   \else\let\next\br@cket\fi
   \ifbr@ck\br@ckfalse\br@ck@t#1,X\else\br@cktrue#1\fi\next#2X}
\def\br@ck@t#1,#2X{\def\temp{#2}\ifx\temp\empty\let\neext\eat
   \else\let\neext\br@ck@t\def\temp{,}\fi
   \def\teemp{#1}\ifx\teemp\empty\else\rcite{#1}\fi\temp\neext#2X}
\def\resetbr@cket{\gdef\[##1]{[\rtag{##1}]}}
\def\references{\resetbr@cket\newsection References\par}

\newtoks\symb@ls\newtoks\s@mb@ls\newtoks\p@gelist\n@wcount\ftn@mber
    \ftn@mber=1\newif\ifftn@mbers\ftn@mbersfalse\newif\ifbyp@ge\byp@gefalse
\def\defm@rk{\ifftn@mbers\n@mberm@rk\else\symb@lm@rk\fi}
\def\n@mberm@rk{\xdef\m@rk{{\the\ftn@mber}}%
    \global\advance\ftn@mber by 1 }
\def\rot@te#1{\let\temp=#1\global#1=\expandafter\r@t@te\the\temp,X}
\def\r@t@te#1,#2X{{#2#1}\xdef\m@rk{{#1}}}
\def\b@@st#1{{$^{#1}$}}\def\str@p#1{#1}
\def\symb@lm@rk{\ifbyp@ge\rot@te\p@gelist\ifnum\expandafter\str@p\m@rk=1 
    \s@mb@ls=\symb@ls\fi\write\f@nsout{\number\count0}\fi \rot@te\s@mb@ls}
\def\byp@ge{\byp@getrue\n@wwrite\f@nsin\openin\f@nsin=\jobname.fns 
    \n@wcount\currentp@ge\currentp@ge=0\p@gelist={0}
    \re@dfns\closein\f@nsin\rot@te\p@gelist
    \n@wread\f@nsout\openout\f@nsout=\jobname.fns }
\def\m@kelist#1X#2{{#1,#2}}
\def\re@dfns{\ifeof\f@nsin\let\next=\relax\else\read\f@nsin to \f@nline
    \ifx\f@nline\v@idline\else\let\t@mplist=\p@gelist
    \ifnum\currentp@ge=\f@nline
    \global\p@gelist=\expandafter\m@kelist\the\t@mplistX0
    \else\currentp@ge=\f@nline
    \global\p@gelist=\expandafter\m@kelist\the\t@mplistX1\fi\fi
    \let\next=\re@dfns\fi\next}
\def\symbols#1{\symb@ls={#1}\s@mb@ls=\symb@ls} 
\def\bigsymbol{\textstyle}
\symbols{\bigsymbol\ast,\dagger,\ddagger,\sharp,\flat,\natural,\star}
\def\ftnumbers{\ftn@mberstrue} \def\ftsymbols{\ftn@mbersfalse}
\def\paginal{\byp@ge} \def\resetftnumbers{\ftn@mber=1}
\def\ftnote#1{\defm@rk\expandafter\expandafter\expandafter\footnote
    \expandafter\b@@st\m@rk{#1}}

\long\def\jump#1\endjump{}
\def\ssum{\mathop{\lower .1em\hbox{$\textstyle\Sigma$}}\nolimits}

\def\qed{\nobreak\kern 1em \vrule height .5em width .5em depth 0em}
\def\newneq{\hbox{\rlap{\hbox to 1\wd9{\hss$=$\hss}}\raise .1em 
   \hbox to 1\wd9{\hss$\scriptscriptstyle/$\hss}}}
\def\subsetne{\setbox9 = \hbox{$\subset$}\mathrel{\hbox{\rlap
   {\lower .4em \newneq}\raise .13em \hbox{$\subset$}}}}
\def\supsetne{\setbox9 = \hbox{$\subset$}\mathrel{\hbox{\rlap
   {\lower .4em \newneq}\raise .13em \hbox{$\supset$}}}}

\def\vbar{\mathchoice{\vrule height6.3ptdepth-.5ptwidth.8pt\kern-.8pt}
   {\vrule height6.3ptdepth-.5ptwidth.8pt\kern-.8pt}
   {\vrule height4.1ptdepth-.35ptwidth.6pt\kern-.6pt}
   {\vrule height3.1ptdepth-.25ptwidth.5pt\kern-.5pt}}
\def\f@dge{\mathchoice{}{}{\mkern.5mu}{\mkern.8mu}}
\def\b@c#1#2{{\rm \mkern#2mu\vbar\mkern-#2mu#1}}
\def\b@b#1{{\rm I\mkern-3.5mu #1}}
\def\b@a#1#2{{\rm #1\mkern-#2mu\f@dge #1}}
\def\bb#1{{\count4=`#1 \advance\count4by-64 \ifcase\count4\or\b@a A{11.5}\or
   \b@b B\or\b@c C{5}\or\b@b D\or\b@b E\or\b@b F \or\b@c G{5}\or\b@b H\or
   \b@b I\or\b@c J{3}\or\b@b K\or\b@b L \or\b@b M\or\b@b N\or\b@c O{5} \or
   \b@b P\or\b@c Q{5}\or\b@b R\or\b@a S{8}\or\b@a T{10.5}\or\b@c U{5}\or
   \b@a V{12}\or\b@a W{16.5}\or\b@a X{11}\or\b@a Y{11.7}\or\b@a Z{7.5}\fi}}

\catcode`\X=11 \catcode`\@=12




\let\thischap\jobname

\def\partof#1{\csname returnthe#1part\endcsname}
\def\chapof#1{\csname returnthe#1chap\endcsname}

\def\setchapter#1,#2,#3.{%
  \expandafter\def\csname returnthe#1part\endcsname{#2}%
  \expandafter\def\csname returnthe#1chap\endcsname{#3}%
}

\setchapter 300a,A,I.
\setchapter 300b,A,II.
\setchapter 300c,A,III.
\setchapter 300d,A,IV.
\setchapter 300e,A,V.
\setchapter 300f,A,VI.
\setchapter 300g,A,VII.
\setchapter   88,B,I.
\setchapter  600,B,II.
\setchapter  705,B,III.

\def\cprefix#1{
\edef\theotherpart{\partof{#1}}\edef\theotherchap{\chapof{#1}}%
\ifx\theotherpart\thispart
   \ifx\theotherchap\thischap 
    \else 
     \theotherchap%
    \fi
   \else 
     \theotherpart.\theotherchap\fi}

\def\sectioncite[#1]#2{%
     \cprefix{#2}#1}

\edef\thispart{\partof{\thischap}}
\edef\thischap{\chapof{\thischap}}


\def\spuriousreset{}


\expandafter\ifx\csname citeadd.tex\endcsname\relax
\expandafter\gdef\csname citeadd.tex\endcsname{}
\else \message{Hey!  Apparently you were trying to
\string\input{citeadd.tex} twice.   This does not make sense.} 
\errmessage{Please edit your file (probably \jobname.tex) and remove
any duplicate ``\string\input'' lines}\endinput\fi

\sectno=-1   
\localtags
\jjtags
\NoBlackBoxes
\define\mr{\medskip\roster}
\define\sn{\smallskip\noindent}
\define\mn{\medskip\noindent}
\define\bn{\bigskip\noindent}
\define\ub{\underbar}
\define\wilog{\text{without loss of generality}}
\define\ermn{\endroster\medskip\noindent}

\define \nl{\newline}
\magnification=\magstep 1
\documentstyle{amsppt}

{    
\catcode`@11

\ifx\alicetwothousandloaded@\relax
  \endinput\else\global\let\alicetwothousandloaded@\relax\fi

\gdef\subjclass{\let\savedef@\subjclass
 \def\subjclass##1\endsubjclass{\let\subjclass\savedef@
   \toks@{\def\usualspace{{\rm\enspace}}\eightpoint}%
   \toks@@{##1\unskip.}%
   \edef\thesubjclass@{\the\toks@
     \frills@{{\noexpand\rm2000 {\noexpand\it Mathematics Subject
       Classification}.\noexpand\enspace}}%
     \the\toks@@}}%
  \nofrillscheck\subjclass}
} 


\expandafter\ifx\csname alice2jlem.tex\endcsname\relax
  \expandafter\xdef\csname alice2jlem.tex\endcsname{\the\catcode`@}
\else \message{Hey!  Apparently you were trying to
\string\input{alice2jlem.tex}  twice.   This does not make sense.}
\errmessage{Please edit your file (probably \jobname.tex) and remove
any duplicate ``\string\input'' lines}\endinput\fi

\expandafter\ifx\csname bib4plain.tex\endcsname\relax
  \expandafter\gdef\csname bib4plain.tex\endcsname{}
\else \message{Hey!  Apparently you were trying to \string\input
  bib4plain.tex twice.   This does not make sense.}
\errmessage{Please edit your file (probably \jobname.tex) and remove
any duplicate ``\string\input'' lines}\endinput\fi

\def\renewcommand{\newcommand}	       
\edef\cite{\the\catcode`@}%
\catcode`@ = 11
\let\@oldatcatcode = \cite
\chardef\@letter = 11
\chardef\@other = 12
%
%
%
%
\def\@innerdef#1#2{\edef#1{\expandafter\noexpand\csname #2\endcsname}}%
%
%
\@innerdef\@innernewcount{newcount}%
\@innerdef\@innernewdimen{newdimen}%
\@innerdef\@innernewif{newif}%
\@innerdef\@innernewwrite{newwrite}%
%
%
%
\def\@gobble#1{}%
%
%
%
\ifx\inputlineno\@undefined
   \let\@linenumber = \empty 
\else
   \def\@linenumber{\the\inputlineno:\space}%
\fi
%
%
%
\def\@futurenonspacelet#1{\def\cs{#1}%
   \afterassignment\@stepone\let\@nexttoken=
}%
\begingroup 
\def\\{\global\let\@stoken= }%
\\ 
\endgroup
\def\@stepone{\expandafter\futurelet\cs\@steptwo}%
\def\@steptwo{\expandafter\ifx\cs\@stoken\let\@@next=\@stepthree
   \else\let\@@next=\@nexttoken\fi \@@next}%
\def\@stepthree{\afterassignment\@stepone\let\@@next= }%
%
%
%
\def\@getoptionalarg#1{%
   \let\@optionaltemp = #1%
   \let\@optionalnext = \relax
   \@futurenonspacelet\@optionalnext\@bracketcheck
}%
%
%
\def\@bracketcheck{%
   \ifx [\@optionalnext
      \expandafter\@@getoptionalarg
   \else
      \let\@optionalarg = \empty
      \expandafter\@optionaltemp
   \fi
}%
\def\@@getoptionalarg[#1]{%
   \def\@optionalarg{#1}%
   \@optionaltemp
}%
%
%
%
\def\@nnil{\@nil}%
\def\@fornoop#1\@@#2#3{}%
\def\@for#1:=#2\do#3{%
   \edef\@fortmp{#2}%
   \ifx\@fortmp\empty \else
      \expandafter\@forloop#2,\@nil,\@nil\@@#1{#3}%
   \fi
}%
\def\@forloop#1,#2,#3\@@#4#5{\def#4{#1}\ifx #4\@nnil \else
       #5\def#4{#2}\ifx #4\@nnil \else#5\@iforloop #3\@@#4{#5}\fi\fi
}%
\def\@iforloop#1,#2\@@#3#4{\def#3{#1}\ifx #3\@nnil
       \let\@nextwhile=\@fornoop \else
      #4\relax\let\@nextwhile=\@iforloop\fi\@nextwhile#2\@@#3{#4}%
}%
%
%
%
\@innernewif\if@fileexists
\def\@testfileexistence{\@getoptionalarg\@finishtestfileexistence}%
\def\@finishtestfileexistence#1{%
   \begingroup
      \def\extension{#1}%
      \immediate\openin0 =
         \ifx\@optionalarg\empty\jobname\else\@optionalarg\fi
         \ifx\extension\empty \else .#1\fi
         \space
      \ifeof 0
         \global\@fileexistsfalse
      \else
         \global\@fileexiststrue
      \fi
      \immediate\closein0
   \endgroup
}%
%
%
%
%
\def\bibliographystyle#1{%
   \@readauxfile
   \@writeaux{\string\bibstyle{#1}}%
}%
\let\bibstyle = \@gobble
%
%
\let\bblfilebasename = \jobname
\def\bibliography#1{%
   \@readauxfile
   \@writeaux{\string\bibdata{#1}}%
   \@testfileexistence[\bblfilebasename]{bbl}%
   \if@fileexists
      \nobreak
      \@readbblfile
   \fi
}%
\let\bibdata = \@gobble
%
%
\def\nocite#1{%
   \@readauxfile
   \@writeaux{\string\citation{#1}}%
}%
\@innernewif\if@notfirstcitation
%
%
\def\cite{\@getoptionalarg\@cite}%
%
%
\def\@cite#1{%
   \let\@citenotetext = \@optionalarg
   \printcitestart
   \nocite{#1}%
   \@notfirstcitationfalse
   \@for \@citation :=#1\do
   {%
      \expandafter\@onecitation\@citation\@@
   }%
   \ifx\empty\@citenotetext\else
      \printcitenote{\@citenotetext}%
   \fi
   \printcitefinish
}%
\newif\ifweareinprivate
\weareinprivatetrue
\ifx\shlhetal\undefinedcontrolseq\weareinprivatefalse\fi
\ifx\shlhetal\relax\weareinprivatefalse\fi
\def\@onecitation#1\@@{%
   \if@notfirstcitation
      \printbetweencitations
   \fi
   \expandafter \ifx \csname\@citelabel{#1}\endcsname \relax
      \if@citewarning
         \message{\@linenumber Undefined citation `#1'.}%
      \fi
     \ifweareinprivate
      \expandafter\gdef\csname\@citelabel{#1}\endcsname{%
\strut 
\vadjust{\vskip-\dp\strutbox
\vbox to 0pt{\vss\parindent0cm \leftskip=\hsize 
\advance\leftskip3mm
\advance\hsize 4cm\strut\openup-4pt 
\rightskip 0cm plus 1cm minus 0.5cm ?  #1 ?\strut}}
         {\tt
            \escapechar = -1
            \nobreak\hskip0pt\pfeilsw
            \expandafter\string\csname#1\endcsname
             \pfeilso
            \nobreak\hskip0pt
         }%
      }%
     \else  
      \expandafter\gdef\csname\@citelabel{#1}\endcsname{%
            {\tt\expandafter\string\csname#1\endcsname}
      }%
     \fi  
   \fi
   \csname\@citelabel{#1}\endcsname
   \@notfirstcitationtrue
}%
%
%
\def\@citelabel#1{b@#1}%
%
%
\def\@citedef#1#2{\expandafter\gdef\csname\@citelabel{#1}\endcsname{#2}}%
%
%
%
\def\@readbblfile{%
   \ifx\@itemnum\@undefined
      \@innernewcount\@itemnum
   \fi
   \begingroup
      \def\begin##1##2{%
         \setbox0 = \hbox{\biblabelcontents{##2}}%
         \biblabelwidth = \wd0
      }%
      \def\end##1{}
      %
      %
      \@itemnum = 0
      \def\bibitem{\@getoptionalarg\@bibitem}%
      \def\@bibitem{%
         \ifx\@optionalarg\empty
            \expandafter\@numberedbibitem
         \else
            \expandafter\@alphabibitem
         \fi
      }%
      \def\@alphabibitem##1{%
         \expandafter \xdef\csname\@citelabel{##1}\endcsname {\@optionalarg}%
         \ifx\biblabelprecontents\@undefined
            \let\biblabelprecontents = \relax
         \fi
         \ifx\biblabelpostcontents\@undefined
            \let\biblabelpostcontents = \hss
         \fi
         \@finishbibitem{##1}%
      }%
      \def\@numberedbibitem##1{%
         \advance\@itemnum by 1
         \expandafter \xdef\csname\@citelabel{##1}\endcsname{\number\@itemnum}%
         \ifx\biblabelprecontents\@undefined
            \let\biblabelprecontents = \hss
         \fi
         \ifx\biblabelpostcontents\@undefined
            \let\biblabelpostcontents = \relax
         \fi
         \@finishbibitem{##1}%
      }%
      \def\@finishbibitem##1{%
         \biblabelprint{\csname\@citelabel{##1}\endcsname}%
         \@writeaux{\string\@citedef{##1}{\csname\@citelabel{##1}\endcsname}}%
         \ignorespaces
      }%
      %
      %
      \let\em = \bblem
      \let\newblock = \bblnewblock
      \let\sc = \bblsc
      \frenchspacing
      \clubpenalty = 4000 \widowpenalty = 4000
      \tolerance = 10000 \hfuzz = .5pt
      \everypar = {\hangindent = \biblabelwidth
                      \advance\hangindent by \biblabelextraspace}%
      \bblrm
      \parskip = 1.5ex plus .5ex minus .5ex
      \biblabelextraspace = .5em
      \bblhook
      \input \bblfilebasename.bbl
   \endgroup
}%
%
%
\@innernewdimen\biblabelwidth
\@innernewdimen\biblabelextraspace
%
%
%
\def\biblabelprint#1{%
   \noindent
   \hbox to \biblabelwidth{%
      \biblabelprecontents
      \biblabelcontents{#1}%
      \biblabelpostcontents
   }%
   \kern\biblabelextraspace
}%
%
%
%
\def\biblabelcontents#1{{\bblrm [#1]}}%
%
%
\def\bblrm{\rm}%
%
%
\def\bblem{\it}%
%
%
\def\bblsc{\ifx\@scfont\@undefined
              \font\@scfont = cmcsc10
           \fi
           \@scfont
}%
%
%
\def\bblnewblock{\hskip .11em plus .33em minus .07em }%
%
%
\let\bblhook = \empty
%
%
%
\def\printcitestart{[}
\def\printcitefinish{]}
\def\printbetweencitations{, }
\def\printcitenote#1{, #1}
%
%
%
\let\citation = \@gobble
%
%
%
\@innernewcount\@numparams
%
%
\def\newcommand#1{%
   \def\@commandname{#1}%
   \@getoptionalarg\@continuenewcommand
}%
%
%
\def\@continuenewcommand{%
   \@numparams = \ifx\@optionalarg\empty 0\else\@optionalarg \fi \relax
   \@newcommand
}%
%
%
\def\@newcommand#1{%
   \def\@startdef{\expandafter\edef\@commandname}%
   \ifnum\@numparams=0
      \let\@paramdef = \empty
   \else
      \ifnum\@numparams>9
         \errmessage{\the\@numparams\space is too many parameters}%
      \else
         \ifnum\@numparams<0
            \errmessage{\the\@numparams\space is too few parameters}%
         \else
            \edef\@paramdef{%
               \ifcase\@numparams
                  \empty  No arguments.
               \or ####1%
               \or ####1####2%
               \or ####1####2####3%
               \or ####1####2####3####4%
               \or ####1####2####3####4####5%
               \or ####1####2####3####4####5####6%
               \or ####1####2####3####4####5####6####7%
               \or ####1####2####3####4####5####6####7####8%
               \or ####1####2####3####4####5####6####7####8####9%
               \fi
            }%
         \fi
      \fi
   \fi
   \expandafter\@startdef\@paramdef{#1}%
}%
%
%
%
%
\def\@readauxfile{%
   \if@auxfiledone \else 
      \global\@auxfiledonetrue
      \@testfileexistence{aux}%
      \if@fileexists
         \begingroup
            \endlinechar = -1
            \catcode`@ = 11
            \input \jobname.aux
         \endgroup
      \else
         \message{\@undefinedmessage}%
         \global\@citewarningfalse
      \fi
      \immediate\openout\@auxfile = \jobname.aux
   \fi
}%
%
%
\newif\if@auxfiledone
\ifx\noauxfile\@undefined \else \@auxfiledonetrue\fi
%
%
%
%
\@innernewwrite\@auxfile
\def\@writeaux#1{\ifx\noauxfile\@undefined \write\@auxfile{#1}\fi}%
%
%
%
\ifx\@undefinedmessage\@undefined
   \def\@undefinedmessage{No .aux file; I won't give you warnings about
                          undefined citations.}%
\fi
%
%
\@innernewif\if@citewarning
\ifx\noauxfile\@undefined \@citewarningtrue\fi
%
%
%
\catcode`@ = \@oldatcatcode

\def\pfeilso{\leavevmode
            \vrule width 1pt height9pt depth 0pt\relax
           \vrule width 1pt height8.7pt depth 0pt\relax
           \vrule width 1pt height8.3pt depth 0pt\relax
           \vrule width 1pt height8.0pt depth 0pt\relax
           \vrule width 1pt height7.7pt depth 0pt\relax
            \vrule width 1pt height7.3pt depth 0pt\relax
            \vrule width 1pt height7.0pt depth 0pt\relax
            \vrule width 1pt height6.7pt depth 0pt\relax
            \vrule width 1pt height6.3pt depth 0pt\relax
            \vrule width 1pt height6.0pt depth 0pt\relax
            \vrule width 1pt height5.7pt depth 0pt\relax
            \vrule width 1pt height5.3pt depth 0pt\relax
            \vrule width 1pt height5.0pt depth 0pt\relax
            \vrule width 1pt height4.7pt depth 0pt\relax
            \vrule width 1pt height4.3pt depth 0pt\relax
            \vrule width 1pt height4.0pt depth 0pt\relax
            \vrule width 1pt height3.7pt depth 0pt\relax
            \vrule width 1pt height3.3pt depth 0pt\relax
            \vrule width 1pt height3.0pt depth 0pt\relax
            \vrule width 1pt height2.7pt depth 0pt\relax
            \vrule width 1pt height2.3pt depth 0pt\relax
            \vrule width 1pt height2.0pt depth 0pt\relax
            \vrule width 1pt height1.7pt depth 0pt\relax
            \vrule width 1pt height1.3pt depth 0pt\relax
            \vrule width 1pt height1.0pt depth 0pt\relax
            \vrule width 1pt height0.7pt depth 0pt\relax
            \vrule width 1pt height0.3pt depth 0pt\relax}

\def\pfeilsw{ \leavevmode 
            \vrule width 1pt height0.3pt depth 0pt\relax
            \vrule width 1pt height0.7pt depth 0pt\relax
            \vrule width 1pt height1.0pt depth 0pt\relax
            \vrule width 1pt height1.3pt depth 0pt\relax
            \vrule width 1pt height1.7pt depth 0pt\relax
            \vrule width 1pt height2.0pt depth 0pt\relax
            \vrule width 1pt height2.3pt depth 0pt\relax
            \vrule width 1pt height2.7pt depth 0pt\relax
            \vrule width 1pt height3.0pt depth 0pt\relax
            \vrule width 1pt height3.3pt depth 0pt\relax
            \vrule width 1pt height3.7pt depth 0pt\relax
            \vrule width 1pt height4.0pt depth 0pt\relax
            \vrule width 1pt height4.3pt depth 0pt\relax
            \vrule width 1pt height4.7pt depth 0pt\relax
            \vrule width 1pt height5.0pt depth 0pt\relax
            \vrule width 1pt height5.3pt depth 0pt\relax
            \vrule width 1pt height5.7pt depth 0pt\relax
            \vrule width 1pt height6.0pt depth 0pt\relax
            \vrule width 1pt height6.3pt depth 0pt\relax
            \vrule width 1pt height6.7pt depth 0pt\relax
            \vrule width 1pt height7.0pt depth 0pt\relax
            \vrule width 1pt height7.3pt depth 0pt\relax
            \vrule width 1pt height7.7pt depth 0pt\relax
            \vrule width 1pt height8.0pt depth 0pt\relax
            \vrule width 1pt height8.3pt depth 0pt\relax
            \vrule width 1pt height8.7pt depth 0pt\relax
            \vrule width 1pt height9pt depth 0pt\relax
      }


\def\widestnumber#1#2{}

\def\citewarning#1{\ifx\shlhetal\relax 
    \else
    \par{#1}\par
    \fi
}

\def\rm{\fam0 \tenrm}

\def\fakesubhead#1\endsubhead{\bigskip\noindent{\bf#1}\par}



%
%
%

%

\font\textrsfs=rsfs10
\font\scriptrsfs=rsfs7
\font\scriptscriptrsfs=rsfs5

\newfam\rsfsfam
\textfont\rsfsfam=\textrsfs
\scriptfont\rsfsfam=\scriptrsfs
\scriptscriptfont\rsfsfam=\scriptscriptrsfs

\edef\oldcatcodeofat{\the\catcode`\@}
\catcode`\@11

\def\Cal@@#1{\noaccents@ \fam \rsfsfam #1}

\catcode`\@\oldcatcodeofat


\expandafter\ifx \csname margininit\endcsname \relax\else\margininit\fi

\long\def\red#1\endred{}
\long\def\green#1\endgreen{}
\long\def\blue#1\endblue{}

\def\endred{ \unmatched endred! }
\def\endgreen{ \unmatched endgreen! }
\def\endblue{ \unmatched endblue! }

\ifx\latexcolors\undefinedcs\def\latexcolors{}\fi

\def\emptycs{}
\def\evaluatelatexcolors{%
        \ifx\latexcolors\emptycs\else
        \expandafter\xxevaluate\latexcolors\xxfertig\evaluatelatexcolors\fi}
\def\xxevaluate#1,#2\xxfertig{\setupthiscolor{#1}%
        \def\latexcolors{#2}}

\font\smallfont=cmsl7
\def\rutgerscolor{\ifmmode\else\endgraf\fi\smallfont
\advance\leftskip0.5cm\relax}
\def\setupthiscolor#1{\edef\tmptmpcs{\noexpand\bgroup\noexpand\rutgerscolor
\noexpand\def\noexpand\currentcolor{#1}%
\noexpand}%
\expandafter\let\csname#1\endcsname\tmptmpcs
\def\tmptmpcs{\checkColorUnmatched{#1}\popthecolor}
\expandafter\let\csname end#1\endcsname\tmptmpcs}

\def\checkColorUnmatched#1{\def\expectcolor{#1}%
    \ifx\expectcolor\currentcolor   
    \else \edef\failhere{\noexpand\tryingToClose '\currentcolor' with end\expectcolor}\failhere\fi}

\def\currentcolor{???}

\def\popthecolor{\ifmmode\else\endgraf\fi\egroup}

\expandafter\def\csname#1\endcsname{}

\evaluatelatexcolors

 \let\outerhead\head
 \def\head{\innerhead}
 \let\innerhead\outerhead

 \let\outersubhead\subhead
 \def\subhead{\innersubhead}
 \let\innersubhead\outersubhead

 \let\outersubsubhead\subsubhead
 \def\subsubhead{\innersubsubhead}
 \let\innersubsubhead\outersubsubhead

 \def\proclaim{\innerproclaim}
 \let\innerproclaim\outerproclaim

 %
 %
 %
 %

\def\demo#1{\medskip\noindent{\it #1.\/}}
\def\enddemo{\smallskip}

\def\remark#1{\medskip\noindent{\it #1.\/}}
\def\endremark{\smallskip}

\pageheight{8.5truein}
\topmatter
\title{Two cardinals models with gap one revisited} \endtitle
\rightheadtext{Two cardinal models}
\author {Saharon Shelah \thanks {\null\newline 
The author would like to thank the Israel Science Foundation for
partial support of this research (Grant No. 242/03). Publication 824.
\null\newline
I would like to thank Alice Leonhardt for the beautiful typing.}
\endthanks} \endauthor  


\affil{Institute of Mathematics\\
 The Hebrew University\\
 Jerusalem, Israel
 \medskip
 Rutgers University\\
 Mathematics Department\\
 New Brunswick, NJ  USA} \endaffil

\abstract   We succeed to say something on the identities of
$(\mu^+,\mu)$ when $\mu > \theta > \text{ cf}(\mu),\mu$ strong limit
$\theta$-compact or even $\mu$ limit of compact cardinals.  
This hopefully will help to prove that
\mr
\item "{$(a)$}"  the pair $(\mu^+,\mu)$ is compact and
\sn
\item "{$(b)$}"  the consistency
of ``some pair $(\mu^+,\mu)$ is not compact", however, this has 
not been proved. 
\endroster
\endabstract
\endtopmatter
\document

\newpage

\head {Annotated Content} \endhead  \resetall 
 \spuriousreset
\bn
\S0 Introduction
\mr
\item "{${{}}$}"  [We give the basic definitions.]
\endroster
\bn
\S1 2-simplicity for gap one
\mr
\item "{${{}}$}"  [We prove that if $\mu = 2^{< \mu}$ then the family
of identities of $(\mu^+,\mu)$ is 2-simple.  So this applies to $\mu$
singular strong limit but also, e.g., to triples
$(\mu^+,\mu,\kappa),\mu = 2^{< \mu} > \kappa$.]
\endroster
\bn
\S2 Successor of strong limit above supercompact:2-identities 
\mr
\item "{${{}}$}"  [Consider a pair $(\mu^+,\mu)$ with $\mu$ strong
limit singular $> \theta > \text{ cf}(\mu),\theta$ a compact
cardinal. We point out quite simply 2-identities which belong to
ID$_2(\mu^+,\mu)$ but not to ID$_2(\aleph_1,\aleph_0)$.]
\endroster
\newpage

\head {\S0 Introduction} \endhead  \resetall \sectno=0
 \spuriousreset
\bigskip

There has been much work on $\kappa$-compactness of pairs
$(\lambda,\mu)$ of cardinals, i.e., when: \ub{if} $T$ is a set of
first order sentences of cardinality $\le \kappa$ and every finite
subset has a $(\lambda,\mu)$ mod $M$ (i.e., $\|M\| = \lambda,|P^M| =
\mu$ for a fixed unary $P$).  \ub{Then} $T$ has a
$(\lambda,\mu)$-model.  

A particularly important case is $\lambda = \mu^+$ in which case this
can be represented as a problem on the $\kappa$-compactness of the
logic $\Bbb L(\bold Q^{\text{card}}_\lambda)$, i.e., 
$(\bold Q^{\text{card}}_{\ge \lambda} x)\varphi$ says that there are at least
$\lambda$ element $x$ satisfying $\varphi_i$.  We deal here only with
this case.  See Furkhen \cite{Fu65}, Morley and Vaught \cite{MoVa62},
Keisler \cite{Ke70}, Mitchel \cite{M1}; for more history see
\cite{Sh:604}. 

Now two cardinal theorems can be translated to partition problems: see
\cite{Sh:8}, \cite{Sh:E17}, lately Shelah and Vaananan
\cite{ShVa:790}.

Restricting ourselves to pairs $(\mu^+,\mu)$, the identities of
$(\aleph_1,\aleph_0)$ were sorted out in \cite{Sh:74}, but we do not
know of the identities of any really different pair $(\mu^+,\mu)$,
i.e., one for which $(\aleph_1,\aleph_0) \nrightarrow (\mu^+,\mu)$.
We know of some such pairs is suitable set theory.  By Mitchel
$(\aleph_2,\aleph_1)$ after suitably collapsing of a Mahlo strongly
inaccessible to $\aleph_2$.  The other, when there is a compact
cardinal in $(\text{cf}(\mu),\mu)$ by Litman and Shelah.  So it would
be nice to know (taking the extreme case).
\bn
\margintag{0.0}\ub{\stag{0.0} Question}:  Assume $\mu$ is a singular cardinal the
limit of compact and even supercompact cardinals.
\nl
1) What are the identities of $(\mu^+,\mu)$?
\nl
2) Is $(\mu^+,\mu) \, \aleph_0$-compact (equivalently $\mu$-compact)?

Note that though we already know that there are some identities of
$(\mu^+,\mu)$ which are not identities of $(\aleph_1,\aleph_0)$ we
have no explicit example.  We give here a partial solution to
\scite{0.0}(1) by finding families of such identities.

Another problem is consistency of failure of compactness.
\nl
In \cite{Sh:604} we have dealt with the simplest case for pairs
$(\lambda,\mu)$ by a reasonable criterion: including no use of 
large cardinals.  From another
perspective the simplest case is the consistency of non compactness of
$\Bbb L(\bold Q),\bold Q$ one cardinality quantifier, and the simplest
one is $\bold Q = \exists^{\ge \mu^+}$.  So we are again drawn to pairs
$(\mu^+,\mu)$, that is gap one instead of gap 2 as in \cite{Sh:604}, 
so necessarily we need to use large cardinals as if, e.g., $\neg 0^\#$
then every such pair is compact.  
\bigskip

\definition{\stag{0.1} Definition}  1) A partial \ub{identity}
\footnote{identification in the terminology of \cite{Sh:8}} 
$\bold s$ is a pair $(a,e) = (\text{Dom}_{\bold s},e_{\bold s})$ where
$a$ is a finite set and $e$ is an equivalence relation on a subfamily
of the family of the finite subsets of $a$, having the property

$$
b\, e\, c \Rightarrow |b| = |c|.
$$
\mn
The equivalence class of $b$ with respect to $e$ will be denoted
$b/e$. \nl
1A) We say $\bold s$ is a full identity or identity if Dom$(e) = {\Cal
P}(a)$. \nl
1B) We say that partial identities $\bold s_1 = (a_1,e_1),\bold s_2 =
(a_2,e_2)$ are isomorphic if there is an isomorphism $h$ from $\bold
s_1$ onto $\bold s_2$ which mean that $h$ is a one-to-one function
from $a_1$ onto $a_2$ such that for every $b_1,c_1 \subseteq a_1$ we
have $(b_1 e_1 c_1) \equiv h(b_1) e_2 h(b_2)$ (so $h$ maps Dom$(e_1)$
onto Dom$(e_2)$).  We define similarly ``$h$ is an embedding of $\bold
s_1$ into $\bold s_2$. \nl
2)  We say that $\lambda \rightarrow (a,e)_\mu$, \ub{if} $(a,e)$ is an
identity or a partial identity and for every function
$f:\,[\lambda]^{<\aleph_0}\to\mu$, there is a one-to-one function
$h:a \rightarrow \lambda$ such that

$$
b\, e\, c \Rightarrow f(h''(b)) = f(h''(c)).
$$
\mn
(Instead Rang$(f) \subseteq \mu$ we may just require $|\text{Rang}(f)|
\le \mu$, this is equivalent). \nl
3)  We define

$$
\text{ID}(\lambda,\mu) =: \{(n,e):\,n < \omega \and (n,e)
\text{ is an identity and } \lambda \rightarrow (n,e)_\mu\}
$$
\mn
and for $f:\,[\lambda]^{<\aleph_0} \rightarrow X$ we let

$$
\align
\text{ID}(f) =: \{(n,e):&(n,e) \text{ is an identity such that for some
one-to-one function} \\
  &h \text{ from } n = \{0,\dotsc,n-1\} \text{ to } \lambda \text{ we
have} \\
  &(\forall b,c \subseteq n)(b\, e\, c \Rightarrow f(h''(b)) = f(h''(c)))\}.
\endalign
$$
\enddefinition
\bn
Clearly two-place functions are easier to understand; this motivates:
\definition{\stag{0.2} Definition}  1) A two-identity or $2$-identity
\footnote{it is not an identity as $e$ is an equivalence relation on
too small set but it is a partial identity} 
is a pair $(a,e)$ where $a$ is a finite set and $e$ 
is an equivalence relation on $[a]^2$. 
Let $\lambda \rightarrow (a,e)_\mu$ mean $\lambda \rightarrow
(a,e^+)_\mu$ where $be^+ c \leftrightarrow [(bec) \vee (b = c \subseteq
a)]$ for any $b,c \subseteq a$.
\nl
2) We defined

$$
\text{ID}_2(\lambda,\mu) =: 
\{(n,e):\,(n,e) \text{ is a } 2\text{-identity and }
\lambda \rightarrow(n,e)_\mu\}
$$
\mn
we define ID$_2(f)$ when $f:[\lambda]^2 \rightarrow X$ as

$$
\align
\biggl\{(n,e):&(n,e) \text{ is a two-identity such that for some } h,
\\
  & \text{ a one-to-one function from } 
\{0,\dotsc,n-1\} \text{ into } \lambda \\
  &\text{ we have } \{\ell_1,\ell_2\} e \{k_1,k_2\} \text{ implies that }
\ell_1 \ne \ell_2 \in \{0,\dotsc,n-1\}, \\
  &k_1 \ne k_2 \in \{0,\dotsc,n-1\} \text{ and }
f(\{h(\ell_1),h(\ell_2)\}) = f(\{h(k_1),h(k_2)\}) \biggr\}.
\endalign
$$
\mn
3) Let us define

$$
\align
\text{ID}^\circledast_2 =: 
\{&({}^n 2,e):\,({}^n 2,e) \text{ is a two-identity and if} \\
  &\{\eta_1,\eta_2\} \ne \{\nu_1,\nu_2\} \text{ are } \subseteq {}^n
2, \text{ then} \\
  &\{\eta_1,\eta_2\}e\{\nu_1,\nu_2\} \Rightarrow \eta_1 \cap \eta_2
= \nu_1 \cap \nu_2\}.
\endalign
$$
\mn
4) In parts (1) and (2) we may replace 2 by $k<\omega$ (only $k <
|a_{\bold s}|$ is interesting) and by $(\le k)$.
\enddefinition
\bn
\margintag{0.3}\ub{\stag{0.3} Discussion}:  By \cite{Sh:49}, 
under the assumption $\aleph_\omega < 2^{\aleph_0}$,
the families ID$_2(\aleph_\omega,\aleph_0)$ and $ID^\circledast_2$ 
coincide (up
to an isomorphism of identities). 
In Gilchrist and Shelah \cite{GcSh:491} and \cite{GcSh:583}
we considered the question of the equality between these
ID$_2(2^{\aleph_0},\aleph_0)$ and ID$^\circledast_2$
under the assumption $2^{\aleph_0} = \aleph_2$. We showed that
consistently the answer may be ``yes" and may be ``no".

Note that $(\aleph_n,\aleph_0) \nrightarrow (\aleph_\omega,\aleph_0)$ 
so ID$(\aleph_2,\aleph_0) \ne \text{ ID}(\aleph_\omega,\aleph_0)$, but for
identities for pairs (i.e. ID$_2$) the question is meaningful.
\bn
We can look more at ordered identities 
\definition{\stag{0.4} Definition}  1) An ord-identity or order
identity is an identity
$\bold s$ such that $a_s \subseteq$ Ord or just: $a$ is an ordered set. \nl
2) $\lambda \rightarrow_{or} (\bold s)_\mu$ if $\bold s$ is an
ord-identity and for every $\bold c:[\lambda]^{< \aleph_0} \rightarrow
\mu$ we have $\bold s \in \text{ OID}(\bold c)$, see below (equivalently
Dom$(\bold c) = [\lambda]^{< \aleph_0},|\text{Rang}(\bold c)| \le \mu$). \nl
3) For $\bold c:[\lambda]^{< \aleph_0} < \mu$ let
OID$(\bold c) = \{(a,e):a$ is a set of ordinals and there is an order
preserving function $f:a \rightarrow \lambda$ such that $b_1 e b_2
\Rightarrow \bold c(f''(b_1)) = \bold c(f''(b_2))\}$. \nl
4) OID$(\lambda,\mu) = \{(n,e):(n,e) \in \text{ OID}(\bold c)$ for every
$\bold c:[\lambda]^{< \aleph_0} \rightarrow \mu \text{ we say } (n,e)
\in \text{ OID}(\bold c)\}$. \nl 
5) Similarly OID$_2$, OID$_k$, OID$_{\le k}$.
\enddefinition
\bn
Of course,
\proclaim{\stag{0.5A} Claim}  1)  {\rm ID}$(\lambda,\mu)$ can be computed from
{\rm OID}$(\lambda,\mu)$. \nl
2) Let $a$ be a finite set of ordinals and $e$ a function.
If $(a,e)$ is an identity, $a$ a set of ordinals and $\lambda >
\mu$, \ub{then} $(a,e) \in { \text{\rm ID\/}}(\lambda,\mu)$ \ub{iff} for some
permutation $\pi$ of $a$ we have $(a,e^\pi) \in { \text{\rm OID\/}}
(\lambda,\mu)$ where $e^\pi = \{(b,c):(\pi''(b),\pi''(c)) \in e\}$. \nl
3) Let $A$ be a set of ordinals, $(a,e)$ an ord-identity and $\bold c$
a function with domain $[A]^{< \aleph_0}$.  \ub{Then} 
$(a,e) \in { \text{\rm ID\/}}(\bold c)$ \ub{iff} for 
some permutation $\pi$ of $a,(a,e^\pi) \in { \text{\rm OID\/}}(\bold c)$. \nl
4) Similarly for 2-identities and $k$-identities and $(\le
k)$-identities and partial identities.
\endproclaim
\bigskip

\proclaim{\stag{0.6} Claim}  For $n \in [1,\omega)$ and $\bold s$ an
ordered partial identity \ub{then} 
there is a first order sentence $\psi_{\bold s}$ such
that: $\psi_{\bold s}$ has a $(\mu^{+n},\mu)$-model \ub{iff} $\bold s
\notin { \text{\rm OID\/}}(\mu^{+n},\mu)$.
\endproclaim
\bigskip

\demo{Proof}  Easy as for some first order $\psi$ sentence if $M$ is a
$(\mu^{+n},\mu)$-model of $\psi$ \ub{then} $<^M$ is a linear order of
$M$ (of cardinality $\mu^{+n}$) which is $\mu^{+n}$-like (i.e. every
initial segment has cardinality).  \hfill$\square_{\scite{0.6}}$\margincite{0.6}
\enddemo
\bn
We define simplicity:
\definition{\stag{0.7} Definition}  1) For $k \le \aleph_0$, we say
$(\lambda,\mu)$ has $k$-simple identities \ub{when} $(a,e) \in
\text{ ID}(\lambda,\mu) \Rightarrow (a,e') \in \text{
ID}(\lambda,\mu)$ whenever:
\mr
\item "{$(*)_k$}"  $a \subseteq \omega,(a,e)$ is an identity of
$(\lambda,\mu)$ and $e'$ is defined by
$$
b e' c \text{ iff } |b|=|c| \and (\forall b'c')[b' \subseteq b \and
|b'| \le k \and c' = \text{ OP}_{c,b}(b') \rightarrow b' ec];
$$
recall OP$_{A,B}(\alpha) = \beta$ iff $\alpha \in A \and \beta \in B
\and \text{ otp}(\alpha \cap A) = \text{ otp}(\beta \cap B)$.
\ermn
2) We define ``$(\lambda,\mu)$ for $k$-simple ordered identities".
\enddefinition
\bn
We can ask \nl
\margintag{0.8}\ub{\stag{0.8} Question}:  1) Define reasonably a pair $(\lambda,\mu)$
such that consistently
\mr
\item "{$\circledast$}"  ID$(\lambda,\mu)$ is not recursive
\sn
\item "{$\circledast'$}"  ID$(\lambda,\mu)$ is not, in a reasonable way,
finitely generated.
\ermn
2) Similarly for ID$_2(\lambda,\mu)$. \nl
3) Restrict yourself to $(\mu^+,\mu)$.
\newpage

\head {\S1 2-simplicity for gap one} \endhead  \resetall \sectno=1
 \spuriousreset
\bigskip

\proclaim{\stag{t.6} Claim}  1) If $\mu$ is strong limit singular 
\ub{then} {\rm ID}$_2(\mu^+,\mu)$ is 2-simple. \nl
2) If $\mu = 2^{< \mu}$ and $c_0:[\mu^+]^{< \aleph_0}
\rightarrow \mu$ \ub{then} we can find $c^*:[\mu^+]^2 \rightarrow \mu$ such
that:
\mr
\item "{$(\alpha)$}"  if $n \in [2,\omega)$ and
$\alpha_0,\dotsc,\alpha_{n-1} < \mu^+$ are with no repetitions and
$\beta_0,\dotsc,\beta_{n-1} < \mu^+$ are with no repetitions and $\ell
< k < n \Rightarrow c^*\{\alpha_\ell,\alpha_k\} =
c^*\{\beta_\ell,\beta_k\}$ \ub{then} $c_0\{\alpha_0,\dotsc,\alpha_{n-1}\} 
= c_0\{\beta_0,\dotsc,\beta_{n-1}\}$ and even
$c^*\{\alpha_0,\dotsc,\alpha_{n-1}\} = c^*\{\beta_0,\dotsc,\beta_{n-1}\}$
\sn
\item "{$(\beta)$}"  if in addition $\alpha_0 < \alpha_1 < \ldots <
\alpha_{n-1}$ \ub{then} $\beta_0 < \beta_1 < \ldots < \beta_{n-3} <
\beta_{n-2},\beta_{n-3} < \beta_{n-1}$.
\endroster
\endproclaim
\bigskip

\remark{\stag{t.6A} Remark}  1) We may wonder what is the gain in
\scite{t.6}(2) as compared to \scite{t.6}(1), as if $\mu = 2^{< \mu}$
is regular then we know all relevant theory on $(\mu^+,\mu)$?  The
answer is that it clarifies identities of triples $(\mu^+,\mu,\kappa)$, e.g.
\mr
\item "{$(a)$}"  $(\mu^+,\mu,\kappa),\mu$ strong limit singular $>
\kappa \ge \text{ cf}(\mu)$
\sn
\item "{$(b)$}"  $(\mu^+,\mu,\kappa),\mu = \mu^{\beth_\omega(\kappa)}$.
\ermn
2) Replacing $\mu^+,2$ by $\mu^{+k},k+1 \ge 2$ is similar and easier.
\endremark
\bigskip

\demo{Proof}  1) By part (2). \nl
2)  By $\boxdot_1 - \boxdot_5$ below the claim is easy (see details in
the end).
\mr
\item "{$\boxdot_1$}"  There is $c_1:[\mu^+]^2 \rightarrow \mu$ such
that if $\alpha_0 < \alpha_1 <\alpha_2 < \mu^+$ and $\beta_0,\beta_1,\beta_2
< \mu^+$ are with no repetitions and $c_1\{\beta_\ell,\beta_k\} =
c_1\{\alpha_\ell,\alpha_k\}$ for $\ell < k < 3$ \ub{then} at least two of
the following holds $\beta_0 < \beta_1,\beta_0 < \beta_2,\beta_1 <
\beta_2$.
\ermn
[Why?  Let $\eta_\alpha \in {}^\mu 2$ for $\alpha < \mu^+$ be pairwise
distinct and for $\alpha \ne \beta < \mu^+$ let
$\varepsilon\{\alpha,\beta\} = \text{ Min}\{\varepsilon:\eta_\alpha
\restriction \varepsilon \ne \eta_\beta \restriction \varepsilon\}$
and define the function $c'_1$ with domain $[\mu^+]^2$ by
$c'_1\{\alpha,\beta\} = \{\eta_\alpha \restriction
\varepsilon\{\alpha,\beta\},\eta_\beta \restriction
\varepsilon\{\alpha,\beta\}\}$, now $|\text{Rang}(c'_1)| \le \mu$
holds because $\mu = 2^{< \mu}$.  For $\alpha \ne \beta$, let
$c''_1\{\alpha,\beta\}$ be 1 if $(\eta_\alpha <_{\text{lex}}
\eta_\beta) \equiv (\alpha < \beta)$ and 0 otherwise (the Sierpinski
colouring).  Lastly, define $c_1$ by
$c_1,c_1\{\alpha,\beta\} =
(c'_1\{\alpha,\beta\},c''_1\{\alpha,\beta\})$, it 
is a function with domain $[\mu^+]^2$ and
range of cardinality $\le \mu$ and easily it is as required.]
\mr
\item "{$\boxdot_2$}"  for every $c:[\mu^+]^{< \aleph_0} \rightarrow \mu$
there is $c_2:[\mu^+]^2 \rightarrow \mu$ such that: if $n \ge 2,\alpha_0 <
\alpha_1 < \ldots < \alpha_{n-1} < \mu^+,\beta_0 < \beta_1 < \ldots <
\beta_{n-1} < \mu^+$ and $\ell < k < n \Rightarrow
c_2\{\alpha_\ell,\alpha_k\} = c_2\{\beta_\ell,\beta_k\}$ then
$c\{\alpha_0,\dotsc,\alpha_{n-1}\} = c\{\beta_0,\dotsc,\beta_{n-1}\}$.
\ermn
[Why?  We are given $c:[\mu^+]^{< \aleph_0} \rightarrow \mu$ and for
each $\alpha < \mu^+$ let $f_\alpha$ be a one-to-one function 
from $\alpha$ onto the ordinal $|\alpha| \le \mu$ and we shall use
those $f_\alpha$'s also later.
\nl
We define an equivalence relation $E$ on $[\mu^+]^2$
\mr
\item "{$(*)$}"  for $\alpha_1 < \beta_1 < \mu^+$ and $\alpha_2 <
\beta_2 < \mu^+$ we have $\{\alpha_1,\beta_1\} E \{\alpha_2,\beta_2\}$
\ub{iff}
{\roster
\itemitem{ $(a)$ }   $f_{\beta_1}(\alpha_1) = 
f_{\beta_2}(\alpha_2)$ and
\sn
\itemitem{ $(b)$ }   for any $n < \omega$ and 
$\gamma_0 < \ldots < \gamma_{n-1} < f_{\beta_1}(\alpha_1)$ we have
\endroster}
$$
c\{\alpha_1,\beta_1,
f^{-1}_{\beta_1}(\gamma_0),\dotsc,f^{-1}_{\beta_1}(\gamma_{n-1})\}
= c\{\alpha_2,\beta_2,
f^{-1}_{\beta_2}(\gamma_0),\dotsc,f^{-1}_{\beta_1}(\gamma_{n-1})\}
$$
and similarly if we omit $\alpha_1,\alpha_2$ and/or $\beta_1,\beta_2$.
\ermn
So $[\mu^+]^2/E$ has cardinality $\le {}^{\mu >} 2 = \mu$ 
and let $c_2:[\mu^+] \rightarrow \mu$
be such that $c_2\{\alpha_1,\beta_1\} = c_2\{\alpha_2,\beta_2\}$ iff
$\{\alpha_1,\beta_1\}/E = \{\alpha_2,\beta_2\}/E$.  We now check that 
it is as required in $\boxdot_2$.
Let $n,\langle \alpha_\ell:\ell < n \rangle,\langle \beta_\ell:\ell <
n \rangle$ be as in $\boxdot_2$; so $\ell < k < n \Rightarrow
c_2\{\alpha_\ell,\alpha_n\} = c_2\{\beta_\ell,\beta_n\}$, hence by
$(*)(a)$ above (for $k=n-1$) we have $\ell < n-1 \Rightarrow
f_{\alpha_{n-1}}(\alpha_\ell) = f_{\beta_{n-1}}(\beta_\ell)$, call it
$\gamma_\ell$.  Let $\ell(*) <n(*)$ be such that $\gamma_\ell$ is maximal.
Now apply $(*)(b)$ with
$\alpha_{\ell(*)},\alpha_{n-1},\beta_{\ell(*)},
\beta_{n-2}$ here standing for
$\alpha_1,\beta_1,\alpha_2,\beta_2$ there and we get the desired result.]
\mr
\item "{$\boxdot_3$}"   In $\boxdot_2$, using $f_\alpha:\alpha
\rightarrow \mu$ as in its proof, we have
$c\{\alpha_0,\dotsc,\alpha_{n-1}\} = c\{\beta_0,\dotsc,\beta_{n-2}\}$
also when
{\roster
\itemitem{ $(*)$ }  $n \ge 2,\alpha_0 < \alpha_1 < \ldots <
\alpha_{n-3} < \alpha_{n-2} < \alpha_{n-1} < \mu^+,\beta_0 < \beta_1 <
\ldots < \beta_{n-3} < \beta_{n-1} < \beta_{n-2}$ and $\ell < n-2
\Rightarrow f_{\alpha_{n-1}}(\alpha_\ell) =
f_{\alpha_{n-2}}(\alpha_\ell)$ and $\ell < k < n \Rightarrow
c_2\{\alpha_\ell,\alpha_k\} = c_2\{\beta_\ell,\beta_k\}$.
\endroster}
\ermn
[Why?  Just the same proof.]
\mr
\item "{$\boxdot_4$}"  there is $c_4:[\mu^+] \rightarrow \mu$ such
that if $\alpha_0 < \alpha_1 < \alpha_2 < \mu^+$ and
$\beta_0,\beta_1,\beta_2 < \mu^+$ with no repetitions,
$c_4\{\beta_\ell,\beta_k\} = c_4\{\alpha_\ell,\alpha_k\}$ for $\ell <
k < 3$ then $\beta_0 < \beta_1 \and \beta_0 < \beta_2$.
\ermn
[Why?  For $\alpha < \beta < \mu^+$ we let $c'\{\alpha,\beta\} =
\{f_\beta(\gamma):\gamma < \alpha \and f_\beta(\gamma) <
f_\beta(\beta)\}$ and let $c_4\{\alpha,\beta\} =
(c'\{\alpha,\beta\},c_1\{\alpha,\beta\},f_\beta(\alpha))$ where $c_1$
is from $\boxdot_1$ and $\langle f_\gamma:\gamma < \mu^+ \rangle$ is
from the proof of $\boxdot_2$. 
Clearly $|\text{Rang}(c')| \le \dsize \sum_{\zeta < \mu} 2^{|\zeta|} =
\mu$ hence $|\text{Rang}(c_4)| \le \mu^3 = \mu$.
If $\alpha_\ell,\beta_\ell(\ell < 3)$ form a counterexample, then 
$c_1\{\alpha_\ell,\alpha_k\} = c_1\{\beta_\ell,\beta_k\}$ for $\ell <
k < 3$ hence by $\boxdot_1$ we have four cases according to which one of the
inequalities $\beta_\ell < \beta_k,\ell < k < 3$ fail.  So the proof of
$\boxdot_4$ splits to three cases.
\bn
\ub{Case 0}:  $\beta_0 < \beta_1 < \beta_2$.

Trivial: the desired conclusion holds.
\bn
\ub{Case 1}:  $\beta_1 < \beta_0$ so $\beta_1 < \beta_0 < \beta_2$.

Let $\zeta_\ell = f_{\alpha_2}(\alpha_\ell)$ for $\ell = 0,1$ hence
$\zeta_0 \ne \zeta_1$ as $f_{\alpha_2}$ is one to one and $\zeta_\ell
= f_{\beta_2}(\beta_\ell)$.  Now on the one hand if
$\zeta_0 < \zeta_1$ then $c'\{\alpha_1,\alpha_2\} 
\ne c'\{\beta_1,\beta_2\}$ (as 
$\zeta_0 \in c'\{\alpha_1,\alpha_2\},\zeta_0 \notin
c'\{\beta_1,\beta_2\})$, contradiction.  On the other hand 
if $\zeta_1 < \zeta_0$ then $c'\{\alpha_0,\alpha_2\} \ne
c'\{\beta_0,\beta_2\}$ (as $\zeta_1 \in c'\{\beta_0,\beta_2\},
\zeta_1 \notin c'\{\alpha_0,\alpha_2\}$), a contradiction, too.
\bn
\ub{Case 2}:  $\beta_2 < \beta_0$.

Then at least one of $\beta_1 < \beta_0,\beta_2 < \beta_1$ hold
contradicting $\boxdot_1$, (i.e., the case we are in).
\bn
\ub{Case 3}:  $\beta_2 < \beta_1$.

By $\boxdot_1$ we have $\beta_0 < \beta_2 < \beta_1$.

This is O.K. for $\boxdot_4$.]
\mr
\item "{$\boxdot_5$}"  for every $c:[\mu^+]^2 \rightarrow \mu$ there
is $c_5:[\mu^+]^2 \rightarrow \mu$ such that
{\roster
\itemitem{ $(a)$ }  $c_5\{\alpha_1,\beta_1\} = c_5\{\alpha_2,\beta_2\}
\Rightarrow c_2\{\alpha_1,\beta_1\} = c_2\{\alpha_2,\beta_2\}$ where
$c_2$ is from $\boxdot_2$ (so also $\boxdot_3$)
\sn
\itemitem{ $(b)$ }  there are no $\alpha_0 < \alpha_1 < \alpha_2 <
\mu^+$ and $\beta_0 < \beta_1 < \beta_2 < \mu^+$ such that
$f_{\alpha_2}(\alpha_0) \ne
f_{\alpha_1}(\alpha_0),c_5\{\alpha_0,\alpha_1\} =
c_5\{\beta_0,\beta_2\},c_5\{\alpha_0,\alpha_2\} = c_5\{\beta_0,\beta_1\}$
and $c_5\{\alpha_1,\alpha_2\} = c_5\{\beta_1,\beta_2\}$
\sn
\itemitem{ $(c)$ }  $c_5\{\alpha_1,\beta_1\} =
c_5\{\alpha_2,\beta_2\} \Rightarrow c_4\{\alpha_1,\beta_1\} =
c_4\{\alpha_2,\beta_2\}$ where $c_4$ is from $\boxdot_4$.
\sn
[Why?  Let $\kappa = \text{ cf}(\mu) \le \mu$ and $\mu = \dsize \sum_{i <
\kappa} \lambda_i$ be such that if $\mu$ is a limit cardinal then 
$\lambda_i$ is (strictly) increasing continuous and if $\mu$ is a successor
cardinal then $\mu = \lambda^+$ and $\lambda_i = \lambda$ for $i <
\kappa$.   We can find $d:[\mu^+]^2 \rightarrow \kappa$ and $\bar g$ such that
\sn
\itemitem{ $\circledast_0$ }  $\qquad (i) \quad$ for 
$\beta < \mu^+,i < \kappa$ the set
$A_{\beta,i} =: \{\alpha < \beta:d\{\alpha,\beta\} \le i\}$ has \nl

\hskip55pt cardinality $\le \lambda_i$ and
\sn
\itemitem{ ${{}}$ } $\qquad (ii) \quad$ if 
$\alpha < \beta < \gamma < \mu^+$ then
$d\{\alpha,\gamma\} \le \text{
max}\{d\{\alpha,\beta\},d\{\beta,\gamma\}\}$
\sn
\itemitem{ ${{}}$ } $\qquad (iii) \quad \bar g$ is a sequence 
$\langle g_\alpha:\alpha < \mu^+ \rangle$ 
\sn
\itemitem{ ${{}}$ }  $\qquad (iv) \quad g_\alpha:\alpha \rightarrow
\mu$ is one to one and \nl

\hskip60pt $\lambda^+_i < \mu \and i < \kappa \and \alpha < 
\beta \Rightarrow ((g_\beta(\alpha) < \lambda^+_i) \equiv$ \nl

\hskip55pt $(d\{\alpha,\beta\} \le i))$
\sn
\itemitem{ ${{}}$ }  $\qquad (v) \quad$ if $\alpha <
\beta,d\{\alpha,\beta\} = i$ and $\lambda^+_i = \mu$ then
$g_\beta(\alpha) < d\{\alpha,\beta\}$.
\endroster}
\ermn
[Why we can find them?  By induction on $\beta < \mu^+$
by induction on $i < \mu$ for $\alpha=f^{-1}_\beta(i)$ 
we choose $d\{\alpha,\beta\}$ and $g_\beta(\alpha)$ as required.]
\sn
Define the functions $c'_6$ and $c'_7$ 
with domain $[\mu^+]^2$ as follows: if $\alpha < \beta$ then \nl
$c'_6\{\alpha,\beta\} = \{(t,\zeta_0,\zeta_1):\zeta_0,\zeta_1 \le
g_\beta(\alpha),t<2$ and $t=0 \Rightarrow g^{-1}_\beta(\zeta_1) <
g_\beta(\zeta_2),t=1 \Rightarrow g_\beta(\zeta_1) > g_\beta(\zeta_2)\}$ and 
$c'_7\{\alpha,\beta\} = \{(t,\zeta,\xi):\zeta \in 
\lambda^+_{d\{\alpha,\beta\}} \cap \text{ Rang}(g_\alpha)$ and $\xi \in
\lambda^+_{d\{\alpha,\beta\}} \cap \text{ Rang}(g_\beta)$ and
$[\lambda^+_{d\{\alpha,\beta\}} = \mu \Rightarrow \zeta <
d\{\alpha,\beta\} \and \xi < d\{\alpha,\beta\}]$  
and $g^{-1}_\alpha(\zeta) < g^{-1}_\beta(\xi) \and t=0$ or
$g^{-1}_\alpha(\zeta) = g^{-1}_\beta(\xi)
\and t=1$ or $g^{-1}_\alpha(\zeta) > g^{-1}_\beta(\xi) \and t=2\}$.

Now for $\alpha < \beta < \mu^+$ we define $c'_5\{\alpha,\beta\} \in
\Pi\{\lambda^+_j:j \le d\{\alpha,\beta\}\}$, we do this by induction
on $\beta$ and for a fixed $\beta$ by induction $i=d\{\alpha,\beta\}$
and for a fixed $\beta$ and $i$ by induction on $\alpha$. \nl
Arriving to $\alpha < \beta$ so $\zeta < 
\lambda^+_{d\{\alpha,\beta\}}$, for each
$j \le d\{\alpha,\beta\}$, let $(c'_5\{\alpha,\beta\})(j)$ be the first
ordinal $\xi < \lambda^+_j$ such that:
\mr
\item "{$\circledast_1$}"  if $\gamma < \beta \and 
d\{\gamma,\beta\} \le j \and (d\{\gamma,\beta\} = d\{\alpha,\beta\} 
\Rightarrow \gamma < \alpha)$ then
$$
(c'_5\{\alpha,\gamma\})(j) < \xi.
$$
\ermn
Clearly possible.  The colouring we use is $c_5$ where for $\alpha <
\beta < \mu^+$ we let
$c_5\{\alpha,\beta\} = (d\{\alpha,\beta\},g_\beta(\alpha),f_\beta(\alpha),
c_2\{\alpha,\beta\},c'_5\{\alpha,\beta\},
c'_6\{\alpha,\beta\},c'_7\{\alpha,\beta\},c_4\{\alpha,\beta\})$,
recalling $c_4$ is from $\boxdot_4$ and $c_2$ is from $\boxdot_2$.
Obviously, $|\text{Rang}(c_5)| \le \mu$ and clauses (a) + (c) of $\boxdot_5$
holds.
So assume $\alpha_0 < \alpha_1 < \alpha_2,\beta_0 < \beta_1 < \beta_2$
form a counterexample to clause (b) of $\boxdot_5$ and we shall
eventually derive a contradiction.
\sn
Clearly
\mr
\item "{$\circledast_2$}"  $(i) \quad d\{\alpha_0,\alpha_2\} =
d\{\beta_0,\beta_1\},d\{\alpha_0,\alpha_1\} =
d\{\beta_0,\beta_2\},d\{\alpha_1,\alpha_2\} = d\{\beta_1,\beta_2\}$
\sn
\item "{${{}}$}"  $(ii) \quad$ similarly for $c',c'_0,c'_1,c_4$.
\ermn
By clause (ii) above we have $d\{\alpha_0,\alpha_2\} \le 
\text{ max}\{d\{\alpha_0,\alpha_1\},d\{\alpha_1,\alpha_2\}\}$, and applying
clause (ii) to $\beta_0 < \beta_1 < \beta_2$ and using $\circledast_2$ we have
$d\{\alpha_0,\alpha_1\} \le \text{ max}\{d\{\alpha_0,\alpha_2\},
d\{\alpha_1,\alpha_2\}$.  Hence
$d\{\alpha_0,\alpha_1\} =
d\{\alpha_0,\alpha_2\} > d\{\alpha_1,\alpha_2\}$ or $\dsize
\bigwedge^2_{\ell=1} [d\{\alpha_0,\alpha_\ell\} \le
d\{\alpha_1,\alpha_2\}]$;
we deal with those two cases separately.
\bn
\ub{Case 1}:  $\varepsilon = d\{\alpha_0,\alpha_1\} = d\{\alpha_0,\alpha_2\} >
d\{\alpha_1,\alpha_2\}$. \nl
So (see the definition of $c'_5$, with
$\alpha_0,\alpha_2,\alpha_1,\varepsilon$ here standing for
$\alpha,\beta,\gamma,j$ there recalling that
$\alpha_0 < \alpha_1 < \alpha_2$) we have $\lambda^+_\varepsilon >
(c'_5\{\alpha_0,\alpha_2\})(\varepsilon) > (c'_5\{\alpha_0,\alpha_1\})
(\varepsilon)$.
Similarly, $\lambda^+_\varepsilon > (c'_5\{\beta_0,\beta_2\})(\varepsilon) >
(c'_5\{\beta_0,\beta_1\})(\varepsilon)$.  
This contradicts $c'_5\{\alpha_0,\alpha_\ell\}
= c'_5\{\beta_0,\beta_{3 - \ell}\}$ for $\ell=1,2$.
\bn
\ub{Case 2}:  $d\{\alpha_0,\alpha_\ell\} \le
d\{\alpha_1,\alpha_2\}$ for $\ell=1,2$.

Let $\varepsilon = d\{\alpha_1,\alpha_2\}$.  Let $\zeta_\ell =
g_{\alpha_\ell}(\alpha_0)$ for $\ell = 1,2$ so $\zeta_\ell =
g_{\beta_{3 - \ell}}(\beta_0)$ for $\ell = 1,2$.
By the assumption toward contradiction, i.e., by a demand
in clause (b) of $\boxdot_5$ we have $\zeta_1 \ne
\zeta_2$.  Clearly $\zeta_\ell < \lambda^+_{d\{\alpha_0,\alpha_\ell\}}
\le \lambda^+_{d\{\alpha_1,\alpha_2\}} = \lambda^+_\varepsilon$ and
$\lambda^+_\varepsilon = \mu \Rightarrow \zeta_\ell <
d\{\alpha_0,\alpha_\ell\} \le d\{\alpha_1,\alpha_2\} \le \varepsilon$.

As $c'_7\{\alpha_1,\alpha_2\} = c'_7\{\beta_1,\beta_2\}$ and
$g^{-1}_{\alpha_1}(\zeta_1) = g^{-1}_{\alpha_2}(\zeta_2)$ clearly
$g^{-1}_{\beta_1}(\zeta_1) = g^{-1}_{\beta_2}(\zeta_2)$ and they are
well defined.

For $\ell=1,2$ as $c_5\{\alpha_0,\alpha_\ell\} = 
c_5\{\beta_0,\beta_{3-\ell}\}$ by the 
choice of $\zeta_\ell$ (that is $\zeta_\ell =
g_{\alpha_\ell}(\alpha_0))$ we have 
$g_{\beta_\ell}(\beta_0) = \zeta_{3-\ell}$ so 
$g^{-1}_{\beta_\ell}(\zeta_{3-\ell}) = \beta_0$ for $\ell=1,2$ hence
$g^{-1}_{\beta_1}(\zeta_2) = g^{-1}_{\beta_2}(\zeta_1)$.
As $c_5\{\alpha_1,\alpha_2\} = c_5\{\beta_1,\beta_2\}$ we have
$c'_7\{\alpha_1,\alpha_2\} = c'_7\{\beta_1,\beta_2\}$ but $\zeta_1,\zeta_2
\le g_{\alpha_2}(\alpha_1)$ hence
\mr
\item "{$\circledast_3$}"  $(g^{-1}_{\alpha_\ell}(\zeta_1) < 
g^{-1}_{\alpha_\ell}(\zeta_2)) \equiv
(g^{-1}_{\beta_\ell}(\zeta_1) < g^{-1}_{\beta_\ell}(\zeta_2))$ 
for $\ell = 1,2$.
\ermn
As $\zeta_1 \ne \zeta_2$ we have $g^{-1}_{\alpha_1}(\zeta_1) \ne
g^{-1}_{\alpha_1}(\zeta_2)$. 

By symmetry \wilog \, $\zeta_1 > \zeta_2$ so 
$g^{-1}_{\beta_1}(\zeta_1) < g^{-1}_{\beta_1}(\zeta_2)$ \ub{iff} (by
equalities above) $g^{-1}_{\beta_2}(\zeta_2) < g^{-1}_{\beta_2}(\zeta_1)$
\ub{iff} (the equivalence in $\circledast_3$) 
$g^{-1}_{\alpha_2}(\zeta_2) < g^{-1}_{\alpha_2}(\zeta_1)$ \ub{iff} by
the choice of $\zeta_1,g^{-1}_{\alpha_2}(\zeta_1) =
\alpha_0),g^{-1}_{\alpha_2}(\zeta_2) < \alpha_0$ \ub{iff} (as 
$c'_5\{\alpha_0,\alpha_2\} = c'_5\{\beta_0,\beta_1\}$ and $\zeta_2 <
\zeta_1 = g_{\alpha_1}(\beta_0)),
g^{-1}_{\beta_1}(\zeta_2) < \beta_0$ \ub{iff} 
(as $\beta_0 = g^{-1}_{\beta_1}(\zeta_1)),g^{-1}_{\beta_1}(\zeta_2) <
g^{-1}_{\beta_1}(\zeta_1)$, clear contradiction. \nl
So we have proved $\boxdot_5$.
\enddemo
\bn
We can now sum up, i.e.:
\demo{Proof of \scite{t.6}(2) from $\boxdot_1 - \boxdot_5$}   We are given
$c_0:[\mu^+]^{< \aleph_0} \rightarrow \mu$.  First we apply $\boxdot_2$
for $c=c_0$ and get $c_2:[\mu^+]^2 \rightarrow \mu$ as there. \nl
Second, we apply $\boxdot_5$ for $c=c_2$ and get $c_5$ as there.  Let
us check that $c_5$ is as required on $c^*$ in \scite{t.6}(2).  So assume
$(*)_0 + (*)_1$ below and (as the case $n=2$ is trivial)
assume $n \ge 3$ where
\mr
\item "{$(*)_0$}"  $\{\alpha_0,\dotsc,\alpha_{n-1}\} \in [\mu^+]^n$
and $\{\beta_0,\dotsc,\beta_{n-1}\} \in [\mu^+]^n$ and
\sn
\item "{$(*)_1$}"  $\ell < k < n \Rightarrow
c_5\{\alpha_\ell,\alpha_k\} = c_5\{\beta_\ell,\beta_k\}$. 
\ermn
Without loss of generality (by renaming)
\mr
\item "{$(*)_2$}"   $\alpha_0 < \ldots < \alpha_{n-1}$. 
\ermn
and it is enough to prove that $c_0\{\alpha_0,\dotsc,\alpha_{n-1}\} =
c_0\{\beta_0,\dotsc,\beta_{n-1}\}$. 
By clause (a) of $\boxdot_5$ we have
\mr
\item "{$(*)_3$}"   $\ell < k < n \Rightarrow
c_2\{\alpha_\ell,\alpha_k\} = c_2\{\beta_\ell,\beta_k\}$. 
\ermn
By clause (c) of $\boxdot_5$ we have
\mr
\item "{$(*)_4$}"  $\ell < k < n \Rightarrow
c_4\{\alpha_\ell,\alpha_k\} = c_4\{\beta_\ell,\beta_k\}$. 
\ermn
Hence by $\boxdot_4$ we have
\mr
\item "{$(*)_5$}"   if $\ell < k < n$ and $\ell < n-2$ then
$\beta_\ell < \beta_k$.
\ermn
[Why?  Apply $\boxdot_4$  to $\alpha_\ell,\alpha_{\ell
+1},\alpha_k;\beta_\ell,\beta_{\ell +1},\beta_k$ if $\ell +1 < k$, and
apply $\boxdot_4$ to \nl
$\alpha_\ell,\alpha_{\ell +1},\alpha_{\ell
+2};\beta_\ell,\beta_{\ell +1},\beta_{\ell +2}$ if $\ell+1=k$.] \nl
So
\mr
\item "{$(*)_6(i)$}"  $\beta_0 < \beta_1 < \ldots < \beta_{n-3} <
\beta_{n-2} < \beta_{n-1}$ or
\sn
\item "{${{}}(ii)$}"  $\beta_0 < \beta_1 < \ldots < \beta_{n-3} <
\beta_{n-1} < \beta_{n-2}$.
\ermn
So clause $(\beta)$ of \scite{t.6} holds.

If $(i)$ of $(*)_6$ holds, then the choice of $c_2$, i.e., by $\boxdot_2$
and $(*)_3$ above we get $c_0\{\alpha_0,\dotsc,\alpha_{n-1}\} =
c_0\{\beta_0,\dotsc,\beta_{n-1}\}$ so we are done.  Otherwise we have
$(ii)$ of $(*)_6$ so by clause (b) of $\boxdot_5$ we have
\mr
\item "{$(*)_7$}"  if $\ell < n-2$ then $f_{\alpha_{n-1}}(\alpha_\ell)
= f_{\beta_{n-2}}(\beta_\ell)$.
\ermn
[Why?  Apply $\boxdot_5(b)$ to $\alpha_\ell,\alpha_{n-2},\alpha_{n-1};
\beta_\ell,\beta_{n-2},\beta_{n-1}$.] \nl
So by $\boxdot_3$ we get $c_0\{\alpha_0,\dotsc,\alpha_{n-1}\} = 
c_0\{\beta_0,\dotsc,\beta_{n-1}\}$ finishing.
\hfill$\square_{\scite{t.6}}$\margincite{t.6}
\enddemo
\bigskip

\proclaim{\stag{t.7} Claim}  Defining {\rm ID}$(\lambda,\mu)$, we
can restrict ourselves to $c:[\lambda]^{< \aleph_0} \rightarrow \mu$
such that $c \restriction [\lambda]^1$ is constant \ub{if} 
{\rm cf}$(\lambda) > \mu$.
\endproclaim
\bigskip

\proclaim{\stag{t.7A} Claim}  1) Assume $\mu = \mu^{< \mu}$ and $n \in
[1,\omega)$.  The identities of {\rm ID}$(\mu^{+n},\mu)$ are $(n+1)$-simple
(and also {\rm OID}$(\mu^+,\mu))$.
\endproclaim
\bigskip

\demo{Proof}  As in \scite{t.6}, only easier in the additional cases.
\hfill$\square_{\scite{t.8}}$\margincite{t.8}
\enddemo
\newpage

\head {\S2 Successor of strong limit above supercompact: 2-identities} \endhead  \resetall \sectno=2
 \spuriousreset
\bigskip

So we know that if $\mu$ is strong limit singular and there is a compact
cardinal in (cf$(\mu),\mu$) then ID$_2(\mu^+,\mu) 
\ne \text{ ID}_2(\aleph_1,\aleph_0)$.  It seems desirable to find
explicitly such 2-identity.
\bn
The proof of the following does much more.
\proclaim{\stag{t.8} Claim}  Assume
\mr
\item "{$(a)$}"  $\bold s_k = (k + \binom k2,e_{{\bold s}_k})$ where the
non-singleton $e_{\bold s_k}$-equivalence classes are the set sets
here $\binom 12 = 0$) \nl

$\{\{\ell_0,\ell_2\}:\ell_0 < k$ and for some $\ell_1 \in 
\{\ell_0+1,\dotsc,k-1\}$ we have 
$\ell_2 = k + {\binom{\ell_1}2} + \ell_0\}$ and \nl

$\{\{\ell_1,\ell_2\}:\ell_1 < k$ and for some $\ell_0 < \ell_1$ we
have $\ell_2 = k + \binom{\ell_1}2 + \ell_0\}$
\sn
\item "{$(b)$}"  $\mu$ is strong limit, $\theta$ a compact cardinal and
\, ${\text{\rm cf\/}}(\mu) < \theta < \mu$.
\ermn
1) $\bold s_k \in { \text{\rm ID\/}}_2(\mu^+,\mu)$, moreover $\bold s_k \in
{ \text{\rm OID\/}}_2(\mu^+,\mu)$. \nl
2) $\bold s_k \notin { \text{\rm ID\/}}_2(\aleph_1,\aleph_0)$ for $k \ge 3$ so
for $k=3$ we have 
$\bold s_k = (6,e_{\bold s})$ and the non-singleton equivalence
classes, after permuting $\{3,5\}$ are $\{\{1,3\},\{0,4\},\{0,5\}\}$ and
$\{\{1,5\},\{2,3\},\{2,4\}\}$. 
\endproclaim
\bigskip

\demo{Proof}  Part (1) 
follows from subclaim \scite{t.8A}(3) below and part (2) follows from
\scite{t.8B} below.  \hfill$\square_{\scite{t.8}}$\margincite{t.8}
\enddemo
\bigskip

\proclaim{\stag{t.8A} Claim}  Assume  
\mr
\item "{$(a)$}"  $\mu$ is strong limit,
\sn
\item "{$(b)$}"   $\theta$ is compact and \, ${\text{\rm cf\/}}(\mu) 
< \theta < \mu$
\sn
\item "{$(c)$}"  $\kappa = \,{\text{\rm cf\/}}(\mu),
\langle \lambda_i:i < \kappa
\rangle$ is increasing with limit $\mu$
\sn
\item "{$(d)$}"  $c:[\mu^+]^2 \rightarrow \mu$
\sn
\item "{$(e)$}"  $d\{\alpha,\beta\} = \,{\text{\rm Min\/}}\{i:c\{\alpha,\beta\}
< \lambda_i\}$.
\ermn
1) We can find $i(*),A,f$ such that
\mr
\item "{$(*)(i)$}"  $i(*) < \kappa,A \in [\mu^+]^{\mu^+}$ and $f:A \rightarrow
\lambda_{i(*)}$ 
\sn
\item "{$(ii)$}"  for every set $B \subseteq A$ of cardinality $<
\theta$ there are $\mu^+$ ordinals $\gamma \in A$ satisfying $(\forall
\alpha \in B)[d\{\alpha,\gamma\} = f(\alpha)]$.
\ermn
2) In part (1) we also have: if $A_1 \subseteq A,|A_1| \ge
\beth_n(\lambda)^+$ and $\lambda_{i(*)} \le \lambda < \mu$, \ub{then} for
some $\langle \gamma_\ell:\ell < n \rangle \in {}^n(\lambda_{i(*)})$ 
and $B \in [A_1]^\lambda$ for every $\alpha_0 <
\ldots < \alpha_{n-1}$ from $B$ for arbitrarily large $\beta <
\lambda$ we have $\ell < n \Rightarrow c\{\alpha_\ell,\beta\} =
\gamma_\ell$. \nl
3) $\bold s_k \in { \text{\rm ID\/}}_2(c)$ 
where $\bold s_k$ is from clause (a) of \scite{t.8}.
\endproclaim
\bigskip

\demo{Proof}  1) Let $D$ be a uniform $\theta$-complete ultrafilter on
$\mu^+$.

Define $f:\mu^+ \rightarrow \kappa$ by $f(\alpha) = i \Leftrightarrow
\{\gamma < \mu^+:d\{\alpha,\gamma\} = i\} \in D$, note that the
function $f$ is well defined as $D$
is a $\theta$-complete ultrafilter on $\mu^+$ and 
$\theta > \kappa$.  So for some $i(*)$, the
set $A =: \{\alpha < \mu^+:f(\alpha) = i(*)\}$ belongs to $D$ and 
check that $(*)$ holds, that is (i) + (ii) hold. \nl
2) Define $c^*:[A]^n \rightarrow {}^n(\lambda_{i(*)})$ such that
\mr
\item "{$\circledast$}"  if $\alpha_0 < \ldots < \alpha_{n-1}$ are
from $A$ then for $\mu^+$ ordinals $\beta < \mu^+$ we have $\langle
c\{\alpha_\ell,\beta\}:\ell < n\}\rangle =
c^*\{\alpha_0,\dotsc,\alpha_{n-1}\}$.
\ermn
So Rang$(c^*)$ has cardinality $\le (\lambda_{i(*)})^n =
\lambda_{i(*)}$ hence by the Erd\"os-Rado theorem there is $B \subseteq A_1$
infinite (even of any pregiven cardinality $< \lambda$) such that $c^*
\restriction [B]^n$ is constant.
 \nl
3) Straight: in part (2) use $n=2,A_1=A$ and get $B$ and $\langle
\gamma_0,\gamma_1 \rangle \in {}^2(\lambda_{i(*)})$ as there and choose
$\alpha_0 < \ldots < \alpha_{k-1}$ from $B$.  
Next choose $\alpha_\ell$ for $\ell = 0,1,\dotsc,\binom k2 -1$,
choosing $\beta_\ell$ by induction on $\ell$.  If 
$\ell = \binom{\ell_1}2 + \ell_0$ and $\ell_0 < \ell_1 < k$
choose $\beta_\ell \in A$ satisfying $\beta_\ell > \alpha_{k-1}$ and 
$\beta_\ell > \beta_m$ for $m < \ell$ such that
$c\{\alpha_{\ell_0},\beta_\ell\} =
\gamma_0,c\{\alpha_{\ell_1},\beta_\ell\} = \gamma_1$. \nl
Now let $\alpha_{k + \ell} = \beta_\ell$ for $\ell < \binom k2$, and
clearly $\langle \alpha_\ell:\ell < k + \binom k2 \rangle$ realize the
identity $\bold s_k$.  \nl
${{}}$  \hfill$\square_{\scite{t.8A}}$\margincite{t.8A}
\enddemo
\bigskip

\proclaim{\stag{t.8B} Subclaim}  1) If $\bold s \in 
{ \text{\rm ID\/}}_2(\aleph_1,\aleph_0)$, \ub{then} we can find a function
$h:[{\text{\rm Dom\/}}_{\bold s}]^2/\bold s \rightarrow \omega$
respecting $e_{\bold s}$ (i.e. $\{\ell_1,\ell_2\} e_{\bold
s}\{\ell_3,\ell_4\} \Rightarrow h\{\ell_1,\ell_2\} =
h\{\ell_3,\ell_4\}$) and there is a linear order $<$ of 
{\rm Dom}$_{\bold s}$ satisfying
\mr
\item "{$\circledast$}"  for any equivalence class $\bold a$ of $e$
there are $a_0,a_1$ such that
{\roster
\itemitem{ $(i)$ }  $a_0,a_1$ are disjoint finite subsets of
{\rm Dom}$_{\bold s}$
\sn
\itemitem{ $(ii)$ }  if $\{\ell_0,\ell_1\} \in \bold a$ and $\ell_0 <
\ell_1$ then $\ell_0 \in a_0 \and \ell_1 \in a_1$
\sn
\itemitem{ $(iii)$ }  if $\ell_0 \ne \ell_1$ are from $a_0 \cup a_1$
and $\{\ell_0,\ell_1\} \notin \bold a$ \ub{then} $h(\{\ell_0,\ell_1\})
> h(\bold a)$.
\endroster}
\ermn
2) We can add in $\circledast$
\mr
\item "{$(iv)$}"   if $\bold a_0,\bold a_1$ are distinct $\bold
e_{\bold s}$-equivalence classes \ub{then} for some $m \in \{0,1\}$ we have
$[\cup \bold a_m]^2 \backslash \bold a_m$ is disjoint to $\bold
a_{1-m}$
\sn
\item "{$(v)$}"  in $\circledast$ above $a_0,a_1$ can be defined as
$\{\ell_0:\{\ell_0,\ell_1\} \in \bold a,\ell_0 <
\ell_1\},\{\ell_1:\{\ell_0,\ell_1\} \in \bold a,\ell_0 < \ell_1\}$
respectively. 
\ermn
3) If $k \ge 3,\bold s_k$ from \scite{t.8} clause (a) \ub{then} $\bold
s_k$ does not belong to {\rm ID}$_2(\aleph_1,\aleph_0)$.
\endproclaim
\bigskip

\demo{Proof}  1) Remember that by \scite{0.5A} 
we can deal with OID$(\aleph_1,\aleph_0)$.
By \cite{Sh:74} we know what is
ID$(\aleph_1,\aleph_0)$, i.e., the  family of identities in
OID$(\aleph_1,\aleph_0)$ is generated by two operations; one is called
duplication and the other of restriction (see below) from the trivial identity
(i.e. $|\text{dom}_{\bold s}| = 1$) and 
we prove $\circledast$ by induction on $n$, the number of times we
need to apply the operations.
\nl
Recall that $(a,e)$ is gotten by duplication if 
we can find sets $a_0,a_1,a_2$ and a function $g$ such that
\mr
\item "{$\circledast_1(a)$}"  $a_0 < a_1 
< a_2$ (i.e. $\ell_0 \in a_0,\ell_1 \in
a_1,\ell_2 \in a_2 \Rightarrow \ell_0 < \ell_1 < \ell_2$)
\sn
\item "{$(b)$}"  $a = a_0 \cup a_1 \cup a_2$
\sn
\item "{$(c)$}"  $g$ a one-to-one order preserving function from $a_0
\cup a_1$ onto $a_0 \cup a_1$ (so $g \restriction a_0 = 
\text{ id}_{a_0}$; let $g_1 = g,g_2 = g^{-1}$
\sn
\item "{$(d)$}"  for $\ell_0 \ne \ell_1 \in (a_0 \cup a_1)$ we have
$\{\ell_0,\ell_1\} e \{g(\ell_0),g(\ell_1)\}$
\sn
\item "{$(e)$}"  if $\ell_1 \in a_1,\ell_2 \in a_2$ then
$\{\ell_1,\ell_2\}/e$ is a singleton
\sn
\item "{$(f)$}"  $\bold s_\ell = (a_0 \cup a_\ell,e \restriction [a_0
\cup a_\ell]^2)$ is from a lower level (up to isomorphism).
\ermn
Recall that $(a,e)$ is gotten by restriction from $(a',e')$ if $a
\subseteq a',e = e' \restriction [a]^2$.

Now we prove the existence of $h$ as required by induction on the
level.  If $|\text{Dom}_{\bold s}| =1$ this is trivial.  If $\bold s$
is gotten by restriction it is trivial too, (as if $\bold s = (a,e),s'
= (a',e'),a' \subseteq a,e' = e \restriction a'$ and $h:[a]^2/e$ is as
guaranteed then we let $h'(\{\ell_0,\ell_1\}/e') = h(\{\ell_0,\ell_1\}/e)$ for
$\ell_0 < \ell_1$.  Easily $h'$ is as required).  So assume $\bold s =
(a,e)$ is gotten by duplication, so let $a_0,a_1,a_2,g_1,g_2$ be as in
$\circledast_1$ and let $h_1$ be as required for $\bold s_1 = (a_0 \cup a_1,e
\restriction [a_0 \cup a_1)^2)$ and similarly define $h_2$ by 
$h_2\{\alpha,\beta\} =
h_1\{g_2(\alpha),g_2(\beta)\}$.  Let $n^* = \text{ max Rang}(h_1)$ and
define $h:[a_0 \cup a_1 \cup a_2]^2 \Rightarrow \omega$ by $h
\supseteq h_1,h \supseteq h_2$ and if $k \in a_1,\ell \in a_2$ then we
let $h\{k,\ell\} = n^*+1$.  Now check.
\nl
2) By symmetry, \wilog \, $h(\bold a_0) < h(\bold a_1)$ and now $m=1$
satisfies the requirement by applying $\circledast_1$ to the equivalence
class $\bold a = \bold a_1$. \nl
3) It is enough to deal with $\bold s_3$.   
By direct checking the criterion in part (2) fails. \nl
${{}}$  \hfill$\square_{\scite{t.8B}}$\margincite{t.8B}
\enddemo
\bn
The following is like \scite{t.8} with $\mu$ just limit (not
necessarily a strong limit cardinal) so
\proclaim{\stag{t.9} Claim}   Assume
\mr
\item "{$(a)$}"  $\bold s'_n \in { \text{\rm OID\/}}_2$ is $(2n+n^2,
e_{{\bold s}'_n})$ where the non-singleton $e_{{\bold s}'_n}$-equivalence
classes are \nl
$\{\{\ell_0,2n+n \ell_0 + \ell_1\}:\ell_0,\ell_1 < n\}$ and \nl
$\{\{n+ \ell_1,2n+n \ell_0 + \ell_1\}:\ell_0,\ell_1 < n\}$
\sn
\item "{$(b)$}"  $\mu$ is a limit cardinal, 
$\mu > \theta > { \text{\rm cf\/}}(\mu)$ and $\theta$ is a compact cardinal
\sn
\item "{$(c)$}"  $s''_n \in { \text{\rm OID\/}}_n$ is $(2^n +
2^{2n},e_{{\bold s}''_s})$ where the non-singleton 
$e_{{\bold s}''_n}$-equivalence classes are: for 
$m<n,\eta \in {}^m 2,i=0,1$ let
$\bold a^i_\eta = \{\{\ell_i,2^n + \binom{2^n}{\ell_0} +
\ell_1\}:\ell_0,\ell_1 < 2^n$ and for some $\nu_0,\nu_1 \in {}^n 2$ we
have $\eta \char 94 \langle 0 \rangle \trianglelefteq \nu_0,\eta \char
94 \langle 1 \rangle \trianglelefteq \nu_1$ and $\ell_0 = \Sigma\{\nu_0(j)
2^j:j < n\}$ and $\ell_1 = \Sigma\{\nu_1(j)2^j:j <n\}\}$.
\ermn
1) $\bold s'_n \in { \text{\rm ID\/}}_2(\mu^+,\mu)$, moreover $\bold s'_n \in
{ \text{\rm OID\/}}_2(\mu^+,\mu)$ similarly for $\bold s'_n$. \nl
2) $\bold s'_n \notin { \text{\rm ID\/}}_2(\aleph_1,\aleph_0)$ for $n
\ge 2$, similarly for $\bold s''_n$.
\endproclaim
\bigskip

\demo{Proof}  1) Like the proof of \scite{t.8A} using \cite{Sh:49} (or
just \cite[\S5]{Sh:604}) instead of the Erd\"os-Rado theorem. \nl
2) Otherwise there is $(a,e) \in \text{ ID}_2(\aleph_1,\aleph_0)$ and
an embedding $h$ of $\bold s'_n$ into $(a,e)$ and by \scite{0.5A}
without loss of generality $h$ is order preserving and $(a,e) \in
\text{ OID}_2(\aleph_1,\aleph_0)$.  Now
\mr
\item "{$(*)_1$}"    if $\ell_0 < n,\ell_1 < n$ and $\ell = 2n + n
\ell_0 + \ell_1$ \ub{then} $h(\ell_0) < h(\ell)$. \nl
[Why?  Choose $\ell'_1 < n,\ell'_1 \ne \ell_1$ and $\ell' = 2n + n
\ell_0 + \ell'_1$, so $\ell \ne \ell'$ and 
$\{\ell_0,\ell\}e_{{\bold s}'_n} \{\ell_0,n + \ell'\}$ hence
$\{h(\ell_0),h(\ell)\},\{h(\ell_0),h(\ell')\}$ are
$e$-equivalent and $h(\ell) \ne h(\ell')$.  But on $(a,e)$
we know that if $\{m_0,m_1\} e \{m_0,m_2\}$ then $m_2 < m_1 < m_0$ and
$m_2 < m_0 < m_1$ are impossible (see \scite{t.10}(2) below) so we are done.]
\sn
\item "{$(*)_2$}"  if $\ell_0 < n,\ell_1 < n$ and $\ell = 2n +
n\ell_0 + \ell_1$ then $h(\ell_1) < h(\ell)$. \nl
[Why?  Like $(*)_1$.]
\ermn
Now we apply \scite{t.8B}(1) + (2) above so $\bold s'_n \notin \text{
ID}_2(\aleph_2,\aleph_1)$.  The conclusion about $\bold s''_n$
follows.  \hfill$\square_{\scite{t.9}}$\margincite{t.9}
\enddemo
\bigskip

\demo{\stag{t.10} Observation}  1) If $k \ge 2,\bold s = 
(n,e) \in \text{ OID}_2(\mu^+,\mu)$ \ub{then} we can find 
$\bold s' = (n',e')$ in fact $n' = 2n-1$ such that:
\mr
\widestnumber\item{$(iii)$}
\item "{$(i)$}"  $e' \restriction [n]^2 = e$
\sn
\item "{$(ii)$}"  $\bold s' \in \text{ ID}(\mu^+,\mu)$
\sn
\item "{$(iii)$}"   for every $c:[\mu^+]^{< \aleph_0} \rightarrow \mu$
there is $c':[\mu^+]^{< \aleph_0} \rightarrow \mu$ refining $c$
(i.e. $c'(u_1) = c'(u_2) \Rightarrow c(u_1) = c(u_2))$ such that: if
$h:\{0,\dotsc,2n-2\} \rightarrow \mu^+$ is one to one and
satisfies $u_1 e' u_2 \Rightarrow c'(h''(u_1)) 
= c'(h''(u_2))$ then $h \restriction \{0,\dotsc,n-1\}$ is increasing.
\ermn
2) There is $c:[\mu^+]^2 \rightarrow \mu$ such that:

if $\alpha,\beta,\gamma$ are distinct and $c\{\alpha,\beta\} =
c\{\alpha,\gamma\}$ then $\alpha < \beta \and \alpha < \gamma$. \nl
3) We can replace in (1), $(\mu^+,\mu)$ by $(\lambda,\mu)$ if there is
$\bold s = (n,e) \in \text{ ID}(\lambda,\mu)$ such that for some
$c:[\lambda]^{< \aleph_0} \rightarrow \mu$ such that
\mr
\item "{$\circledast$}"   if $h:n \rightarrow \lambda$ induces 
$e_{\bold s}$ then $h(0) < h(1)$.
\endroster
\enddemo
\bigskip

\demo{Proof}  1) Define $e':u_1 e' u_2$ \ub{iff} $u_1
e u_2 \vee u_1 = u_2 \vee \dsize \bigvee_{\ell < n-1}
(u_1 = \{\ell,n+ \ell+1\} \and u_2 e \{\ell,\ell+1\}) \vee \dsize
\bigvee_{\ell <n} (u_2 = \{\ell,n + \ell +1\} \and u_1 e \{\ell,\ell +1\})$.  
Now use (2). 
\nl
2) Let $f_\alpha:\alpha \rightarrow \mu$ be one to one and let $<^*$ a dense
linear order on $\mu^+$ with $\{\alpha:\alpha < \mu\}$ a dense
subset.  Now choose $c_1:[\mu^+]^2 \rightarrow \mu$ such that $\alpha < \beta
\Rightarrow \alpha \le^* c_1\{\alpha,\beta\} <^* \beta$ and
$c:[\mu^+]^2 \rightarrow \mu$ be $\alpha < \beta \Rightarrow
c\{\alpha,\beta\} = \text{ pr}(f_\beta(\alpha),c_1\{\alpha,\beta\})$
for some pairing function pr. \nl
3) Similar to part (1) only $|\text{Dom}_{\bold s'}|$ is larger.
\hfill$\square_{\scite{t.10}}$\margincite{t.10}
\enddemo
\newpage

     \shlhetal 

\nocite{ignore-this-bibtex-warning} 
\newpage
    
REFERENCES.  
\bibliographystyle{lit-plain}
\bibliography{lista,listb,listx,listf,liste}

\enddocument